%% file: mainRev.tex
\documentclass[10pt]{article}

\input{defs.tex}
\usepackage{todonotes}

  \title{Energy-conserving explicit and implicit time integration methods for the multi-dimensional Hermite-\DG{} 
    discretization of the Vlasov-Maxwell equations}
\author[]{%
C. Pagliantini\thanks{\textit{Corresponding author}. Centre for Analysis, Scientific computing and Applications,
Eindhoven University of Technology,
5600 MB Eindhoven, The Netherlands.
Email: \texttt{c.pagliantini@tue.nl}}\;,
G. Manzini\thanks{T-5 Applied Mathematics and Plasma Physics Group,
Los Alamos National Laboratory, Los Alamos, NM 87545, USA.
Email: \texttt{gmanzini@lanl.gov}}\;,
O. Koshkarov\thanks{T-5 Applied Mathematics and Plasma Physics Group,
Los Alamos National Laboratory, Los Alamos, NM 87545, USA.
Email: \texttt{koshkarov@lanl.gov}}\;,
G. L. Delzanno\thanks{T-5 Applied Mathematics and Plasma Physics Group,
Los Alamos National Laboratory, Los Alamos, NM 87545, USA.
Email: \texttt{delzanno@lanl.gov}}\;, and
V. Roytershteyn\thanks{Space Science Institute, 4765 Walnut St, Suite B,
    Boulder, CO 80301, USA.
Email: \texttt{vroytershteyn@spacescience.org}}
}
\date{}

\newcommand{\rev}[1]{\textcolor{black}{#1}}
\newcommand{\revtwo}[1]{\textcolor{black}{#1}}

\begin{document}
\maketitle

\input{abstract.tex}
  

\renewcommand{\arraystretch}{1.}
\raggedbottom
\input{sec1_intro.tex}
\input{sec2_model.tex}
\input{sec3_discrete_time.tex}
\input{sec4_discrete_space.tex}
\input{sec5_invariants.tex}
\input{sec6_cons_GL_RK.tex}

\input{sec7_results.tex}
\input{conclusion.tex}


\section*{Declarations}
\textbf{Funding.}
The work of GLD, OK, GM was supported by the Laboratory Directed
Research and Development - Exploratory and Research (LDRD-ER) Program
of Los Alamos National Laboratory under project number 20170207ER.
Los Alamos National Laboratory is operated by Triad National Security,
LLC, for the National Nuclear Security Administration of
U.S. Department of Energy (Contract No. 89233218CNA000001).
Computational resources for the SPS-\DG{} simulations were provided by
the Los Alamos National Laboratory Institutional Computing Program.

VR's contributions were supported by DOE grant
DE-SC0019315.

\medskip
\noindent
\textbf{Data availability statement.}
Data sharing not applicable to this article as no datasets were generated or analysed during the current study.


\bibliographystyle{plain}
\bibliography{ref}


\end{document}

%% file: defs.tex
\usepackage{etex}
\usepackage[utf8]{inputenc} 
\usepackage[american]{babel}
\usepackage[T1]{fontenc} 
\usepackage{amsmath,amssymb,amsfonts,amsthm}
\usepackage{caption}
\usepackage{authblk}

\usepackage[text={6.25in,9in},centering]{geometry}
\usepackage{graphicx}
\usepackage{bm,bbm}
\usepackage{comment}
\usepackage{rotating}
\usepackage{color}
\usepackage{multirow}
\usepackage{multicol}
\usepackage{numprint}
\usepackage{overpic}
\usepackage[numbers,sort]{natbib}
\usepackage{float}
\usepackage{mdframed}
\usepackage{outlines}
\usepackage{savesym}
\savesymbol{AND}
\usepackage{algpseudocode}
\usepackage{algorithm}
\usepackage[
	pdftex,
	bookmarks,
	colorlinks=true, 
	pdftitle={Energy-preserving RK},
	pdfauthor={Cecilia Pagliantini},
	pdfsubject={},
	pdfkeywords={},
	urlcolor= blue,
	linkcolor= blue,
	citecolor= blue
]{hyperref}
\usepackage{url}

\theoremstyle{plain}
\newtheorem{theorem}{Theorem}[section]
\newtheorem{proposition}{Proposition}[section]
\newtheorem{lemma}{Lemma}[section]

\theoremstyle{definition}
  \newtheorem{remark}{Remark}[section]
  \newtheorem{definition}{Definition}[section]

\newcommand{\DG}{DG}

\definecolor{MyDarkGreen}{rgb}{0,0.45,0}
\newcommand{\CP}[1]{\noindent{\color{blue}{[CP:~#1]}}}

\newcommand{\gv}{\bm{g}}

\newcommand{\nv}{\bm{n}}

\newcommand{\uv}{\bm{u}}
\newcommand{\vv}{\bm{v}}

\newcommand{\xv}{\bm{x}}

\newcommand{\Bv}{\bm{B}}

\newcommand{\Ev}{\bm{E}}

\newcommand{\Jv}{\bm{J}}

\newcommand{\Ov}{\bm{O}}

\newcommand{\Sv}{\bm{S}}

\newcommand{\Uv}{\bm{U}}

\newcommand{\UvC}{\bm{u}}

\newcommand{\as}{a}
\newcommand{\bs}{b}
\newcommand{\cs}{c}
\newcommand{\ds}{d}

\newcommand{\fs}{f}

\newcommand{\ks}{k}

\newcommand{\ms}{m}
\newcommand{\ns}{n}

\newcommand{\qs}{q}

\renewcommand{\ss}{s}
\newcommand{\ts}{t}
\newcommand{\us}{u}
\newcommand{\vs}{v}

\newcommand{\xs}{x}
\newcommand{\ys}{y}
\newcommand{\zs}{z}

\newcommand{\As}{A}
\newcommand{\Bs}{B}

\newcommand{\Fs}{F}

\newcommand{\Hs}{H}
\newcommand{\Is}{I}

\newcommand{\Ls}{L}

\newcommand{\Ns}{N}

\newcommand{\calA}{\mathcal{A}}
\newcommand{\calB}{\mathcal{B}}
\newcommand{\calC}{\mathcal{C}}
\newcommand{\calD}{\mathcal{D}}
\newcommand{\calE}{\mathcal{E}}

\newcommand{\calH}{\mathcal{H}}
\newcommand{\calI}{\mathcal{I}}

\newcommand{\calL}{\mathcal{L}}

\newcommand{\calV}{\mathcal{V}}

\newcommand{\calY}{\mathcal{Y}}

\renewcommand{\fs}{f^s}

\renewcommand{\As}{\calA}
\renewcommand{\Bs}{\calB}
\renewcommand{\Ls}{\calL}

\newcommand{\Yv}{\ensuremath{\mathbf{Y}}}

\newcommand{\xiv}{{\bm\xi}}

\newcommand{\C}{C}

\newcommand{\PT}[1]{\ensuremath{\dfrac{\partial #1}{\partial t}}}

\renewcommand{\as}{\alpha^s}
\renewcommand{\us}{u^s}

\newcommand{\GRADX}{\nabla_{\xv}}
\newcommand{\GRADV}{\nabla_{\vv}}

\newcommand{\matFx}{\mathbbm{F}_x}

\newcommand{\matF}{\mathbbm{F}}

\newcommand{\dxv}{d\xv}
\newcommand{\dvv}{d\vv}

\newcommand{\dS}{dS}

\newcommand{\Nn}{\ensuremath{N_{v_x}}}
\newcommand{\Nm}{\ensuremath{N_{v_y}}}
\newcommand{\Np}{\ensuremath{N_{v_z}}}

\newcommand{\REAL}{\mathbbm{R}}

\newcommand{\vphi}{\bm\varphi}
\newcommand{\varphiIl}{\varphi^{I,l}}

\renewcommand{\ms}{m^s}
\renewcommand{\qs}{q^s}

\newcommand{\dy}{\,dy}
\newcommand{\dz}{\,dz}

\newcommand{\dt}{dt}
\renewcommand{\ds}{d\sigma}

\newcommand{\PS}[1]{\mathbb{P}_{#1}}

\newcommand{\ABS}  [1]{\left|#1\right|}

\def\trait #1 #2 #3 {\vrule width #1pt height #2pt depth #3pt}
\def\fin{\hfill
        \trait .3 5 0
        \trait 5 .3 0
        \kern-5pt
        \trait 5 5 -4.7
        \trait 0.3 5 0
\medskip}
\newcommand{\ENDPROOF}{\fin}
\newcommand{\BEGINPROOF}{\emph{Proof}.~~}

\usepackage{scalerel,stackengine}
\stackMath
\newcommand\reallywidehat[1]{%
\savestack{\tmpbox}{\stretchto{%
  \scaleto{%
    \scalerel*[\widthof{\ensuremath{#1}}]{\kern-.6pt\bigwedge\kern-.6pt}%
    {\rule[-\textheight/2]{1ex}{\textheight}}
  }{\textheight}%
}{0.5ex}}%
\stackon[1pt]{#1}{\tmpbox}%
}

\newcommand{\NFLX}[1]{\reallywidehat{#1}}        
\newcommand{\JUMP}[1]{\big[\hspace{-1mm}\big[#1\big]\hspace{-1mm}\big]}

\newcommand{\JMPF}{\mathcal{J}_{\F}}

\renewcommand{\Hs} {\mathcal{H}}

\renewcommand{\Fs} {\mathcal{F}}
\newcommand{\FNs}{\mathcal{F}^N}
\newcommand{\FN}[1]{\mathcal{F}^N_{#1}}

\newcommand{\Ox}{\Omega_x}
\renewcommand{\Ov}{\Omega_\vs}

\newcommand{\fsN}{f^{s,N}}
\newcommand{\fsNz}{f^{s,N}_{0}}

\newcounter{numbs}
\newcounter{numbi}
\newcounter{numbii}

\newcommand{\BegBoxQuote}[1]{
   \vspace{\baselineskip}
   \begin{center}
   \begin{minipage}{#1\textwidth}
   \begin{mdframed}
}
\newcommand{\EndBoxQuote}{
   \end{mdframed}
   \end{minipage}
   \end{center}
}

\newcommand{\ocs}{\omega_{cs}}
\newcommand{\oce}{\omega_{ce}}
\newcommand{\ope}{\omega_{pe}}

\newcommand{\Ex}{E_x}
\newcommand{\Ey}{E_y}
\newcommand{\Ez}{E_z}
\newcommand{\ExN}{E_x^N}
\newcommand{\EyN}{E_y^N}
\newcommand{\EzN}{E_z^N}

\newcommand{\Bx}{B_x}
\newcommand{\By}{B_y}
\newcommand{\Bz}{B_z}
\newcommand{\BxN}{B_x^N}
\newcommand{\ByN}{B_y^N}
\newcommand{\BzN}{B_z^N}

\newcommand{\vpsi}{\psi}

\newcommand{\cN}{\mathcal{N}}

\newcommand{\Il} {I,l}

\newcommand{\F}{\mathsf{f}}
\newcommand{\Fip}{\mathsf{f}_{i+\frac{1}{2},j,k}}

\newcommand{\SPAN}{\mathsf{span}}

\newcommand{\vecE}  {\Ev}
\newcommand{\vecEN} {{\Ev}^N}
\newcommand{\vecENz}{{\Ev}^N_0}
\newcommand{\vecEIl}{{\Ev}^{\Il}}
\newcommand{\ExIl}  {\Ex^{\Il}}
\newcommand{\EyIl}  {\Ey^{\Il}}
\newcommand{\EzIl}  {\Ez^{\Il}}

\newcommand{\vecB}  {\Bv}
\newcommand{\vecBN} {{\Bv}^N}
\newcommand{\vecBNz}{{\Bv}^N_0}
\newcommand{\vecBIl}{{\Bv}^{\Il}}
\newcommand{\BxIl}  {\Bx^{\Il}}
\newcommand{\ByIl}  {\By^{\Il}}
\newcommand{\BzIl}  {\Bz^{\Il}}

\newcommand{\vecJ}{\Jv}
\newcommand{\vecJN}{\Jv^N}

\newcommand{\FF}{\mathbb{F}}
\newcommand{\Fvx}{\mathbf{F}_x}
\newcommand{\Fvy}{\mathbf{F}_y}
\newcommand{\Fvz}{\mathbf{F}_z}

\newcommand{\LTWO}{L^2}

\newcommand{\Nl}{N_l}

\newcommand{\FFE}{\FF_{\Ev}}
\newcommand{\FFB}{\FF_{\Bv}}

\newcommand{\csVM} {\frac{\qs}{\ms}\frac{\oce}{\ope}}

\newcommand{\fse}{f^s}
\newcommand{\Bve}{\Bv}
\newcommand{\Eve}{\Ev}

\newcommand{\calVN} {\calV^N} \newcommand{\calHN} {\calH^N}
\newcommand{\calHpN}{\mathcal{\widetilde{H}}^N}

\newcommand{\fsNt}{\ensuremath{f^{s,N,\tau}}}
\newcommand{\fsNtp}{\ensuremath{f^{s,N,\tau+1}}}

\newcommand{\vecENt}{\ensuremath{\Ev^{N,\tau}}}

\newcommand{\vecBNt}{\ensuremath{\Bv^{N,\tau}}}

\newcommand{\BsE}{\Bs_{\Ev}}
\newcommand{\BsB}{\Bs_{\Bv}}

\newcommand{\KIN}{{\bf kin}}
\newcommand{\TOT}{{\bf tot}}

\newcommand{\MOD}{\textbf{mod}}
\newcommand{\ELE}{\Ev,\Bv}
\newcommand{\JMP}{{\bf jump}}
\newcommand{\BND}{{\bf bnd}}

\newcommand{\phiRK} {\varphi_{RK}}
\newcommand{\phiGL} {\varphi_{GL}}
\newcommand{\phiRKh}{\widehat{\varphi}_{RK}}
\newcommand{\bhati}{\widehat{b}_i}

\newcommand{\ysb}{\overline{\ys}}

\newcommand{\NRK}{{N_{RK}}}
\newcommand{\NGL}{{N_{GL}}}
\newcommand{\Ng}{{N_{T}}}

\newcommand{\species}{\text{s}}

%% file: abstract.tex
\begin{abstract}
  We study the conservation properties of the Hermite-discontinuous Galerkin (Hermite-DG) approximation of the Vlasov-Maxwell equations. In this semi-discrete formulation, the total mass is preserved independently for every plasma species. Further, an energy invariant exists if central numerical fluxes are used in the DG approximation of Maxwell's equations, while a dissipative term is present when upwind fluxes are employed. In general, traditional temporal integrators might fail to preserve invariants associated with conservation laws during the time evolution. Hence, we analyze the capability of explicit and implicit Runge-Kutta (RK) temporal integrators to preserve such invariants. Since explicit RK methods can only ensure preservation of linear invariants but do not provide any control on the system energy, we consider modified explicit RK methods in the family of relaxation Runge-Kutta methods (RRK).
 These methods can be tuned to preserve the energy invariant at the continuous or semi-discrete level, a distinction that is important when upwind fluxes are used in the discretization of Maxwell's equations since upwind provides a numerical source of energy dissipation that is not present when central fluxes are used. We prove that the proposed methods are able to preserve the energy invariant and to maintain the semi-discrete energy dissipation (if present) according to the discretization of Maxwell's equations. An extensive set of numerical experiments corroborates the theoretical findings. It also suggests that maintaining the semi-discrete energy dissipation when upwind fluxes are used leads to an overall better accuracy of the method relative to using upwind fluxes while forcing exact energy conservation.
\end{abstract}

%% file: sec1_intro.tex
\section{Introduction}

The Vlasov-Maxwell equations for modeling collisionless plasmas in
presence of self-consistent electromagnetic fields possess an infinite
number of invariants, i.e., quantities that do not change while the
physical system evolves in
time~\cite{Glassey:1996,Bittencourt:2004,Boyd-Sandserson:2003}.
These quantities are related to the moments of
the particle distribution functions with respect to the independent
velocity variable and include terms in their definition that may
depend on the electromagnetic fields.
Importantly, a few of them express fundamental properties of the
physical system such as the conservation of the number of particles,
which also entails mass and charge conservation, the conservation of
the total momentum and the conservation of the total energy.

Reproducing some (or all) of the conservation properties
mentioned above in the discrete setting makes a numerical method very
compelling, but at the same time it is a nontrivial
task~\cite{Juno-Hakim-TenBarge-Shi-Dorland:2018}. 
Indeed, it turns out that for only very few numerical methods we can
define discrete invariants that can be calculated from the numerical
solution and correspond to the number of particles, the total
momentum or the total energy of the physical model.

\revtwo{
The majority of the literature on numerical schemes able to preserve key invariants of motion at the fully discrete level has focused on particle-based discretizations.
Variational algorithms were introduced in \cite{Le70,Le72,East91} based on least action principles. Energy-conserving time integrators for finite element particle-in-cell discretizations of the Vlasov–Maxwell equations have been recently proposed in \cite{KS21,CPKS22}.
Using the Hamiltonian formulation of some kinetic plasma models, Refs.~\cite{He15,Tao16,GEMPIC} have focused on the preservation of the phase-space structure, via e.g. symplectic integrators, while ensuring global bounds of the numerical errors of conserved quantities.
}

\medskip
In a series of previous works~\cite{Delzanno:2015,Manzini-Delzanno-Markidis-Vencels:2016,koshkarov2021}, it was shown that the spectral
expansion of the plasma distribution function, \rev{see also \cite{Schumer-Holloway:1998}}, makes it possible to formulate numerical methods for the solution of the kinetic equations that show excellent conservation properties and posses many other desirable properties. Such methods are particularly attractive for problems that involve fluid-kinetic coupling, i.e., the
coupling between macroscopic and microscale system dynamics. That is because a suitable choice of the spectral basis makes the fluid-kinetic
coupling an intrinsic property of the numerical approximation: the low-order terms of the expansion provide a fluid description of
the plasma while additional higher-order terms are able to capture
features from the underlying kinetics.


\medskip
In the present work, we investigate the conservation properties for one such method: the numerical approximation of the system of Vlasov-Maxwell equations
provided by the spectral-discontinuous Galerkin method proposed in
\cite{koshkarov2021}.
The method is based on a spectral expansion in velocity space using
asymmetrically-weighted Hermite basis functions and an approximation
in physical space using piecewise discontinuous polynomials of a given order.
The discontinuous Galerkin (DG) method is also applied to the
discretization in space of the time-dependent Maxwell equations. Unlike the work presented in~\cite{koshkarov2021}, here we treat the fully discrete case and focus on the conservation of the total number of particles and the total energy of the system.

The conservation of the number of particles, which implies mass and
charge conservation, is a consequence of the conservative formulation
of the method, and holds regardless of the discontinuity of the
approximate solution fields.
However, the discontinuous nature of the approximate distribution
functions and electromagnetic fields plays a major role in the
conservation of the total momentum and the total energy.
Indeed, the momentum is not conserved and a dissipative term
proportional to the square of the jumps of the electromagnetics fields
affects the energy conservation when upwind numerical fluxes are used in the discretization
of Maxwell's equations.
The total energy is instead conserved when central
numerical fluxes are adopted for Maxwell's equations.

\rev{We remark that the work of Ref. \cite{koshkarov2021} does not deal with temporal integration. On the contrary, the present manuscript studies the numerical temporal discretization of the system of ordinary differential equations resulting from the aforementioned DG-Hermite discretization. In particular, we study the effect of time integration on the conservation of mass and energy, i.e., the invariants of the semi-discrete problem.}
The availability of
time discretizations able to satisfy the conservation laws of
the continuous or semi-discrete problem is of fundamental importance to ensure
numerical solutions that reproduce key physical properties of kinetic
models.
With this objective in mind, we consider two families of one-step
temporal integrators: explicit Runge-Kutta methods and (implicit)
Gauss-Legendre methods.
First, we show that
the implicit Gauss-Legendre schemes are able to
preserve both the conserved quantities and the dissipation of the
total energy, meaning that the amount of energy dissipated by the
fully-discrete scheme equals the temporal approximation of the
dissipative term ensuing from the Hermite-DG semi-discrete formulation.
Secondly, we study the properties of explicit RK schemes:
they ensure the preservation of the number of particles (linear
invariant) but do not guarantee any control on the evolution of the
total energy.
However, in practical implementation of plasma solvers, the explicit
schemes are often preferred to the implicit ones despite the superior
stability and conservation properties of the latter.
The reason is that explicit schemes are computationally cheap,
unlike implicit schemes whose computational cost is dominated by the solution of a
nonlinear system of equations at each time step.
This remark also suggests that the availability of \emph{explicit} RK
methods of arbitrary order able to preserve the conservation and
dissipation laws of the continuous or semi-discrete problem could have a significant
impact on the efficient solution of kinetic plasma models.
We address this issue
by resorting to
explicit RK methods based on projection.
Projection methods, \emph{cf.} for instance \cite[Section 5.3]{ESF98} and \cite{GQ05}, belong to the class of extrinsic numerical time integrators and are based on the projection of the numerical solution onto the manifold that is constrained by a given conservation property.
The projection methods considered in this work consist in modifying the RK step by suitable scalar factors that depend on the approximate solution and are determined by the given conservation/dissipation law.
This approach was originally proposed in
\cite[pp. 265--266]{DV84} to derive conservative RK methods.
The idea was extended in \cite{DBM02} to a restricted class of fourth-order methods and further developed in \cite{Calvo-HernandezAbreu-Montijano-Randez:2006} via the so-called incremental direction technique (IDT). More recently, an extension of these methods to preserve inner-product norms of the solution has led to the so-called relaxation Runge--Kutta (RRK) methods \cite{Ketcheson19}.
We analyze this family of modified RK methods in the context of the Hermite-DG discretization of the Vlasov-Maxwell equations.
For the semi-discrete problem considered, this correction guarantees the preservation of the energy invariant and can be tuned to preserve it either at the continuous level, or at the semi-discrete level (this distinction is only relevant when upwind numerical fluxes are used in Maxwell's equations, since for central fluxes the total energy is conserved at both continuum and semi-discrete levels). With numerical experiments, we show that the latter choice that maintains the semi-discrete energy dissipation when upwind numerical fluxes are used in the discretization of Maxwell's equations leads to a better overall accuracy of the method relative to a method that uses upwind fluxes but forces exact energy conservation.
  

\medskip 
The remainder of the paper is organized as follows.
In Section~\ref{sec:Vlasov-Maxwell:eqs}, we introduce the
Vlasov-Maxwell equations as the mathematical model that describes
the transport phenomena of different charged particle species in a
collisionless plasma under the action of the self-consistent
electromagnetic field.
Section~\ref{sec:Runge-Kutta:methods} is devoted to the temporal
discretization of the Vlasov-Maxwell equations by using explicit
Runge-Kutta methods and implicit Gauss-Legendre methods that are
reformulated as a single time-evolution problem.
In Section~\ref{sec:variational:formulation}, we summarize the
space-velocity discretization of the Vlasov-Maxwell equations
performed with the Hermite-DG methods introduced in
\cite{koshkarov2021}.
Next, Section~\ref{sec:invariants} pertains to the study of the
conservation laws associated with the total number of particles and
the total energy and to their characterization as conserved quantities
of the ODE resulting from the semi-discrete system of
Vlasov-Maxwell equations.
In Section~\ref{sec:cons-general}, we analyze the capability of
Runge-Kutta temporal integrators to preserve the conservative and
dissipative relations associated with the different spatial
discretizations, and we study modified explicit schemes that preserve
mass and total energy of the semi-discrete Vlasov-Maxwell problem.
In Section~\ref{sec:numerical:results}, we experimentally assess the
performance of the proposed methods on the whistler instability test
case, on the Orzsag-Tang vortex test case, and in resolving dispersion
properties of high frequency waves.
Finally, in Section~\ref{sec:conclusions} we present some conclusions
and final considerations.

\medskip
\emph{Notation and normalization.} 
We normalize the model equations as follows.
Time $t$ is normalized to the electron plasma frequency
$\ope=\sqrt{e^2 n^{e}_0\slash{}\varepsilon_0m^e}$, where $e$ is the
elementary charge, $m^e$ is the electron mass, $\varepsilon_0$ is the
permittivity of vacuum, and $n^e_0$ is a reference electron density.
The velocity coordinate $\vs$ is normalized to the speed of light
$\cs$; the spatial coordinate $\xs$ is normalized to the electron
inertial length $d_e=\cs\slash{\ope}$; the magnetic field $\vecB$ is
normalized to a reference magnetic field $\Bs_0$, and the electric
field is normalized to $\cs\Bs_0$.
We denote the quantities associated with a given plasma species by the
superscript $s$, which may take the specific values $s=e$ (electrons)
and $s=i$ (ions), etc.
Accordingly, we denote the mass of the particles of species $s$ by
$\ms$ and their charge by $\qs=\pm e$.  
We normalize the charge $\qs$ and the mass $\ms$ to the elementary charge $e$ and the
mass $m^e$, respectively.
Finally, we define the cyclotron frequency of species $s$ as
$\ocs=e\Bs_0/\ms$. For simplicity, we maintain the same symbol for each normalized variable and from now on we will only consider normalized quantities.


%% file: sec2_model.tex

\section{Vlasov-Maxwell equations}
\label{sec:Vlasov-Maxwell:eqs}
The Vlasov-Maxwell equations provide a model for
collisionless plasmas~\cite{goldston1995}.
At any time instant $t>0$,
the behavior of the particles of species $s$ in the plasma
is described by the non-negative distribution function
$\fs(\xv,\vv,t)$, where
$\xv$ denotes the position in the physical space $\Ox$ and $\vv$
the position in the velocity space $\Ov$.
The plasma evolves self-consistently
under the action of
the electric and magnetic fields
$\Eve(\xv,t)$ and $\Bve(\xv,t)$ generated by particles' motion (and external sources): the
distribution function of species $s$ satisfies the (normalized) Vlasov
equation
\begin{equation}
  \frac{\partial\fse}{\partial t} + \vv\cdot\GRADX\fse + 
  \csVM\big(\Eve+\vv\times\Bve\big)\cdot\GRADV\fse = 0.
  \label{eq:Vlasov:fs}
\end{equation}
Let $\vecJ$ and $\rho$ denote the self-consistent electric current and
charge density, respectively, induced by the plasma particles, namely
\begin{equation*}
  \rho (\xv,t) = \sum_{s=1}^{\Ns_{\species}}\qs\int_{\Ov}    \fse(\xv,\vv,t)\dvv,\qquad
  \vecJ(\xv,t) = \sum_{s=1}^{\Ns_{\species}}\qs\int_{\Ov} \vv\fse(\xv,\vv,t)\dvv,
\end{equation*}
where $\Ns_{\species}$ denotes the number of plasma species.
The electric and magnetic fields $\Eve$ and $\Bve$
satisfy the time-dependent wave propagation equations
\begin{align}
  &\frac{\partial\Eve}{\partial t} - \GRADX\times\Bve =  - \frac{\ope}{\oce}\vecJ,\label{eq:dEdt}\\[0.5em]
  &\frac{\partial\Bve}{\partial t} + \GRADX\times\Eve =  0,    \label{eq:dBdt}
\end{align}
and the divergence equations
\begin{align*}
  &\GRADX\cdot\Eve = \frac{\ope}{\oce}\rho,\\[0.5em]
  &\GRADX\cdot\Bve = 0.
\end{align*}

We consider the unbounded velocity space $\Ov=\REAL^{3}$ and we assume
that, for $\ABS{\vv}\to\infty$,
each distribution function $\fse$ decays sufficiently fast, e.g. as
$\exp(-\ABS{\vv}^2)$ ~\cite{Glassey:1996}.
For example, this assumption is physically consistent with near-Maxwellian velocity
distribution of a plasma close to thermodynamic
equilibrium~\cite{Grad:1949}.
Similarly, we consider the closed bounded subset $\Ox\subset\REAL^3$
with boundary $\partial\Ox$, and we assume that suitable
problem-dependent boundary conditions for $\fse$, $\Eve$, and $\Bve$
are prescribed at $\partial\Ox$ for any time $t$ and any $\vv$ in $\Ov$.
Moreover, we consider the plasma evolution in the temporal interval $[0,T]$ and assume that physically meaningful initial conditions are provided 
for the unknown fields $\fse$, $\Eve$, $\Bve$ at the initial time
$t=0$.
Since we pursue a numerical approximation of the Vlasov-Maxwell equations
based on a spatial discontinuous Galerkin method,
we reformulate Eqs.~\eqref{eq:dEdt}-\eqref{eq:dBdt} in
conservative form as follows.
Let us introduce, for any $t>0$ and $\xv\in\Ox$, the vector of
conservative unknowns $\UvC$ and the source term $\Sv$,
\begin{align}
  \UvC(\xv,t) := 
  \left(
    \begin{array}{c}
      \vecE(\xv,t) \\[0.5em]
      \vecB(\xv,t)
    \end{array}
  \right),
  \qquad
  \Sv(\xv,t) = \left(
    \begin{array}{c}
      \Sv_{\vecE} \\[0.5em]
      \Sv_{\vecB}
    \end{array}
  \right):= -\frac{\ope}{\oce} \left(
    \begin{array}{c}
      \vecJ(\xv,t) \\[0.5em]
      \mathbf{0}
    \end{array}
  \right),
    \label{eq:Maxwell:conservative:05}
\end{align}
and the fluxes,
\begin{align}
  \qquad
  \FF(\UvC) := 
  \left(
    \begin{array}{c}
      \FFE(\UvC)\\[0.5em]
      \FFB(\UvC)
    \end{array}
  \right),
  \qquad
  \nabla_{\xv}\cdot\FF(\UvC) =
  \left(
    \begin{array}{c}
      \nabla_{\xv}\cdot\FFE(\UvC)\\[0.5em]
      \nabla_{\xv}\cdot\FFB(\UvC)
    \end{array}
  \right) =
  \left(
    \begin{array}{c}
      -\nabla_{\xv}\times\vecB\\[0.5em]
      \nabla_{\xv}\times\vecE
    \end{array}
  \right),
  \label{eq:divg:FF}
\end{align}
where the operator $\nabla_{\xv}\,\cdot\,$ denotes here the row-wise
divergence in physical space. An explicit form for fluxes $\FF$ is given in Section~\ref{ssec:semi-disc-VF}.
With these definitions, Maxwell's equations \eqref{eq:dEdt} and
\eqref{eq:dBdt} can be written in conservative form as
\begin{equation}\label{eq:Maxwell:conservative:00}
  \frac{\partial\UvC}{\partial t} + \nabla_{\xv}\cdot\FF(\UvC) = \rev{\Sv(\xv,t)}.
\end{equation}


%% file: sec3_discrete_time.tex

\section{Temporal discretization of the Vlasov-Maxwell equations}
\label{sec:Runge-Kutta:methods}

For the numerical approximation of the Vlasov-Maxwell equations, we
pursue a method of line approach with Runge-Kutta temporal integrators
coupled with the DG-spectral approximation in space-velocity
introduced in \cite{koshkarov2021}, and summarized in
Section~\ref{sec:variational:formulation}.

Since the focus of this work is on numerical time integrators, we
start by reformulating the Vlasov equation~\eqref{eq:Vlasov:fs} for
every plasma species $s$, and the wave propagation equations
\eqref{eq:dEdt}-\eqref{eq:dBdt} as a single time-evolution problem and
consider its temporal discretization.
To this aim, consider the vector-valued function collecting all the unknowns
\begin{align}\label{eq:yvecC}
  \Yv(\xv,\vv,t) = 
  \left(
    \begin{array}{c}
      \big(\fs(\xv,\vv,t)\big)_{s}\\[0.5em]
      \Eve(\xv,t)\\[0.5em]
      \Bve(\xv,t)
    \end{array}
  \right),\qquad t>0,\; \xv\in\Ox,\; \vv\in\Ov,\; s=1,\ldots,\Ns_{\species}\,.
\end{align}
The function $\Yv$ satisfies the partial differential equation
\begin{align}\label{eq:PDE}
  \dfrac{\partial \Yv}{\partial t} = \Fs(\Yv;\xv,\vv,t)
  \qquad\mbox{in}\quad\Ox\times\Ov\times (0,T],
\end{align}
with suitable initial conditions $\Yv(\xv,\vv,0)=\Yv_{0}(\xv,\vv)$,
and boundary conditions, as described
in~Section~\ref{sec:Vlasov-Maxwell:eqs}.
%
The operator $\Fs$ in \eqref{eq:PDE} is defined as
\begin{align}\label{eq:rhs}
  \Fs(\Yv;\xv,\vv,t) = 
  \left(
    \begin{array}{c}
      \big(
      - \vv\cdot\GRADX\fse - \csVM\big(\Eve+\vv\times\Bve\big)\cdot\GRADV\fse
      \big)_{s}
      \\[1em]
    \rev{\Sv(\xv,t)}-\nabla_{\xv}\cdot\FF(\UvC)
    \end{array}
  \right),
\end{align}
where the first term corresponds to the Vlasov equation
\eqref{eq:Vlasov:fs}, while the second term is associated with
Maxwell's equations in conservative form
\eqref{eq:Maxwell:conservative:00}.
The boundary conditions can be introduced in the above setting through
additional conditions on the unknown functions $\fs$ and $\UvC$ or by
suitably including them in the discretization in space and velocity
that is introduced in Section~\ref{sec:variational:formulation}.

\medskip
For the temporal discretization of problem~\eqref{eq:PDE}, we split
the time integration domain $(0,T]$ into the union of intervals
$(t^{\tau},t^{\tau+1}]$, $\tau\in\mathbb{N}_0$, with time step
$\Delta\ts^{\tau}=t^{\tau+1}-t^{\tau}$.
Let us denote by $\Yv^\tau(\xv,\vv)$ the approximation of the
function~\eqref{eq:yvecC} at time $t^\tau$, and assume that the
initial condition $\Yv^0(\xv,\vv)$ is given, for all $\xv\in\Ox$ and
$\vv\in\Ov$.
In the forthcoming formulation of the temporal schemes
for~\eqref{eq:PDE} we omit, for the sake of better readability, the
dependence on the variables $\xv$ and $\vv$.
We consider the following two families of approximate one-step
temporal integrators:

\medskip
\begin{enumerate}
\item\textit{Explicit $\NRK$-stage Runge-Kutta methods (\emph{cf.}
  \cite{HNW93} and references therein).}\\
  The sequence $\{\Yv^{\tau}(\xv,\vv)\}_{\tau}$ of numerical solutions
  of \eqref{eq:PDE} with initial condition $\Yv^0$
  is obtained as
  \begin{equation}\label{eq:RK}
    \begin{aligned}
      \Yv^{\tau+1}=\phiRK(\Yv^{\tau}) &:= \Yv^{\tau} + \Delta\ts^{\tau}\sum_{i=1}^{\NRK}\bs_{i}\,\Fs(\calY_i\,;t^\tau+c_i \Delta\ts^{\tau}),\\
      \calY_{1} &:= \Yv^{\tau},\\
      \calY_{i} &:= \Yv^{\tau} + \Delta\ts^{\tau}\sum_{j=1}^{i-1}a_{ij}\,
      \Fs(\calY_j;t^\tau+c_j \Delta\ts^{\tau}),
      \qquad i=2,\ldots,\NRK.
    \end{aligned}
  \end{equation}
  The explicit $\NRK$-stage Runge-Kutta method is, thus, characterized
  by the Runge-Kutta matrix $(a_{ij})_{1\leq j<i\leq \NRK}$ and the
  set of coefficients ${(b_i)}_{i=1}^\NRK$, ${(c_i)}_{i=1}^\NRK$.
  
  \medskip
\item\textit{Gauss-Legendre methods of order $2\NGL$ (\emph{cf.}
  \cite[Chapter IV.5]{Hairer-Wanner:1996} and references therein).}\\
  Gauss-Legendre temporal integrators are collocation methods based on
  Gaussian quadrature formulas and belong to the family of implicit
  Runge-Kutta methods.
  More specifically, the sequence $\{\Yv^{\tau}(\xv,\vv)\}_{\tau}$ of
  numerical solutions of \eqref{eq:PDE} with initial condition $\Yv^0$
  is derived as
  \begin{equation}\label{eq:GL}
    \begin{aligned}
      \Yv^{\tau+1} = \phiGL(\Yv^{\tau}) &:= \Yv^{\tau} + \Delta\ts^{\tau}\sum_{i=1}^{\NGL}\bs_{i}\,\Fs(\calY_i\,;t^\tau+c_i \Delta\ts^{\tau}),\\
      \calY_{i} &:= \Yv^{\tau} + \Delta\ts^{\tau}\sum_{j=1}^{\NGL}a_{ij}\,
      \Fs(\calY_j;t^\tau+c_j \Delta\ts^{\tau}),
      \qquad i=1,\ldots,\NGL\,,
    \end{aligned}
  \end{equation}
  with
  \begin{equation*}
    a_{ij} := \int_0^{c_i}\dfrac{q_j(\sigma)}{q_j(c_j)}\, d\sigma,\qquad
    b_i := \int_0^1 \dfrac{q_i(\sigma)}{q_i(c_i)}\, d\sigma,\qquad i,j=1,\ldots,\NGL,
  \end{equation*}
  and
  \begin{equation*}
    q_i(\sigma) := \dfrac{\Pi_{j=1}^{\NGL} (\sigma-c_j)}{\sigma-c_i},\qquad i=1,\ldots,\NGL.
  \end{equation*}
  The Gauss-Legendre method of order $2$ is the implicit midpoint
  rule. 
  For higher orders, the values of the coefficients $(a_{ij})_{1\leq
    i,j\leq \NGL}$, ${(b_i)}_{i=1}^\NGL$, and ${(c_i)}_{i=1}^\NGL$
  in~\eqref{eq:GL} can be found in compendia on numerical methods for
  ordinary differential equations (ODE), see, e.g.,
  \cite[Chapter~IV.5]{Hairer-Wanner:1996}.
\end{enumerate}
Throughout, we will adopt the shorthand notation
$\FN{i}:=\FNs(\calY_i;t^{\tau}+c_i\Delta t^{\tau})$.

\medskip
The choice of an explicit or an implicit temporal integrator is
dictated by several factors, and it is usually problem-dependent.
Explicit methods are computationally attractive since their
implementation to solve \eqref{eq:PDE} only requires $\NRK$
evaluations of the function $\Fs$ per time step: this allows accurate
approximations at a competitive computational cost.
However, explicit schemes suffer from time step restrictions due to
stability requirements and, as problems become increasingly stiff,
implicit methods might become more convenient.
Moreover, \rev{traditional} explicit temporal integrators usually fail to preserve at the discrete level the conservation properties of the continuous system.
\rev{Indeed such schemes are not even guaranteed to preserve polynomial invariants of degree strictly larger than 1.}
This might trigger spurious behavior of the approximate problem and
yield unphysical solutions, particularly for long-time integration.

Implicit temporal integrators exhibit superior stability properties
when compared to explicit schemes.
In particular, Gauss-Legendre methods are also symmetric and
symplectic, and, hence, are well-suited for the approximation of
problems over long temporal intervals.
Although the stability properties allow for larger time steps, the
implementation costs of implicit methods are dominated by the
iterative solvers required for the solution of a nonlinear system of
equations at each time step.
Therefore, the efficiency of implicit methods strongly relies on the
availability of fast linear solvers which, for problems featuring a
large number of unknowns, also implies the need for efficient
preconditioning strategies.
The development of efficient solvers is out of the scope of the
present work and might provide an interesting direction for future
investigation.


%% file: sec4_discrete_space.tex

\section{Hermite-DG approximation of the Vlasov-Maxwell equations in space and velocity}
\label{sec:variational:formulation}

Concerning the discretization in space and velocity of
system~\eqref{eq:PDE}, we consider the spectral-DG method proposed
in~\cite{koshkarov2021}, where
Hermite functions provide the spectral approximation of the Vlasov
equation in the velocity space, while a discontinuous Galerkin method
is used for the spatial discretization. \revtwo{For reader's convenience and to introduce the relevant notation, the method described in \cite[Section 3]{koshkarov2021} is briefly summarized in this section. Differently from \cite{koshkarov2021}, we present the semi-discretization in space and velocity from an equivalent variational formulation perspective.}

\medskip
The spatial domain $\Ox$ is
partitioned into $N_c=N_xN_yN_z$ cubic or
regular hexahedral cells, with $N_x$. $N_y$, $N_z$ the number of cells in the $x$, $y$, $z$ direction, respectively.
We assume that the partition is uniform with mesh size $\Delta\xs$,
$\Delta\ys$, and $\Delta\zs$, in each direction, and label the mesh
elements by the indices $i,j,k$ running from $1$ to $N_x$,
$N_y$, and $N_z$, respectively.
For convenience of exposition, the generic mesh cell is labeled by the
letter $I$ and the summation over all mesh cells is indicated by
$\sum_{I}$, without the summation bounds being specified.
With some abuse of notation, we may subindex $\Is$ as $\Is_{i,j,k}$.
Accordingly, triplets with an half-integer index, e.g.
$(i\pm\frac{1}{2},j,k)$, $(i,j\pm\frac{1}{2},k)$ and
$(i,j,k\pm\frac{1}{2})$ label the cell interfaces that are orthogonal
to the $x$-, $y$-, and $z$-direction, respectively, and delimiting
cell $\Is_{i,j,k}$.
With this notation, for example, two consecutive cells in the
$x$-direction are denoted by $\Is_{i,j,k}$ and $\Is_{i+1,j,k}$ and are
separated by the cell interface $\Fip$.
The faces are oriented such that the normal vector to each face always
points outwards.

For the spatial discretization, we consider multivariate polynomial
functions whose restriction to a given mesh cell $\Is$ is a polynomial
of degree at most $N_{\DG{}}$.
The basis for the local polynomial space $\PS{N_{\DG{}}}(\Is)$ on the
mesh element $\Is$ are denoted by $\{\varphi^{\Il}\}$ for
$l=1,\ldots,\Nl$,
so that
$\PS{N_{\DG{}}}(I)=\SPAN\big\{\varphi^{\Il}\}_{l=1}^{\Nl}$. Here $\Nl\in\mathbb{N}_+$ depends on the polynomial degree $N_{\DG{}}$ and on the spatial dimension.
We construct the local polynomials $\varphi^{\Il}$ as tensor product
of the univariate Legendre polynomials
$\{\Ls_\zeta\}_{\zeta=0}^{N_{\DG{}}}$ of degree $\zeta$, which are
defined in the interval $[-1,1]$ and, then, suitably rescaled and
translated on every mesh cell $\Is$.
The polynomials $\Ls_\zeta$ form an orthogonal basis for
$\PS{N_{\DG{}}}([-1,1])$ \cite{Funaro:2008}.

\medskip
Concerning the spectral approximation in velocity of the Vlasov
equation, we consider, as
in~\cite[Section~3.1]{koshkarov2021},
the univariate asymmetrically weighted Hermite functions defined as
\begin{align}
  \psi_{\zeta}(\xi^s_{\beta})
  = \big(\pi\,2^{\zeta}\,\zeta!\big)^{-\frac{1}{2}}\Hs_{\zeta}(\xi^s_{\beta})\exp\big(-(\xi^s_{\beta})^2\big),\qquad
  \psi^{\zeta}(\xi^s_{\beta}) =
  \big(2^{\zeta}\,\zeta!\big)^{-\frac{1}{2}}\Hs_{\zeta}(\xi^s_{\beta}),\label{eq:Psi-updw}
\end{align}
where $\Hs_{\zeta}(\xi^s_{\beta})$ is the $\zeta$-th univariate
Hermite polynomial in the velocity direction $v_\beta$ for
$\beta(\zeta)\in\{\xs,\ys,\zs\}$, $\zeta\in\{n,m,p\}$ and
$n=0,\ldots,\Nn$, $m=0,\ldots,\Nm$, $p=0,\ldots,\Np$. 
Moreover, \revtwo{$\xiv^{s} = \big(\xi^{s}_{x},\xi^{s}_{y},\xi^{s}_{z}\big)^T$, with $\xi^{s}_{\beta}=\frac{\vs_{\beta}-\us_{\beta}}{\as_{\beta}}$}
with the quantities $\us_\beta\in\mathbb{R}$ and
$\as_\beta\in\mathbb{R}$ constant factors that depend on the plasma
species (and that are specified by the user).
The Hermite functions $\Psi_{n,m,p}$ and $\Psi^{n,m,p}$ are given by
the tensor product of the univariate Hermite
functions~\eqref{eq:Psi-updw} as
\begin{align*}
  \Psi_{n,m,p}(\xiv^s)=\vpsi_{n}(\xi^{s}_{x})\vpsi_{m}(\xi^{s}_{y})\vpsi_{p}(\xi^{s}_{z})
  \quad\textrm{and}\quad
  \Psi^{n,m,p}(\xiv^s)=\vpsi^{n}(\xi^{s}_{x})\vpsi^{m}(\xi^{s}_{y})\vpsi^{p}(\xi^{s}_{z}).
\end{align*}

\medskip
The Hermite-\DG{} variational formulation of the Vlasov-Maxwell system
requires the introduction of the following finite dimensional spaces:
\begin{equation}\label{eq:VHN}
\begin{aligned}
  \calHN &:=\textrm{span}\Big\{\Psi_{n,m,p},\,\textrm{for~}(n,m,p)\in\big\{(0,0,0),\ldots,(N_{v_x},N_{v_y},N_{v_z})\big\}\Big\},\\[0.5em]
  \calHpN&:=\textrm{span}\Big\{\Psi^{n,m,p},\,\textrm{for~}(n,m,p)\in\big\{(0,0,0),\ldots,(N_{v_x},N_{v_y},N_{v_z})\big\}\Big\},\\[0.5em]
  \calVN &:=\textrm{span}\Big\{\vphi^{\Il},\,\textrm{for~}\Is\equiv\Is_{i,j,k},(i,j,k)\in\big\{(1,1,1),\ldots(N_x,N_y,N_z)\big\},\,\,l=1,\dots,\Nl\Big\}.
\end{aligned}
\end{equation}


\subsection{Semi-discrete variational formulation}
\label{ssec:semi-disc-VF}

For any time $t\in[0,T]$, we assume that the numerical distribution
function $\fsN(\cdot,\cdot,t)$ belongs to $\calHN\times\calVN$ for any
plasma species $s$ and therefore it can be written as a linear
combination of the Hermite and \DG{} basis functions as
\begin{align}
  \fsN(\xv,\vv,t) =
  \sum_{n,m,p}\sum_{\Il}\C_{n,m,p}^{s,\Il}(t)\Psi_{n,m,p}(\xiv^s)\varphi^{\Il}(\xv),
  \qquad \forall\, t>0,\,\xv\in\Ox,\, \vv\in\Ov\,,
  \label{eq:fsN:def}
\end{align}
where, to ease the notation, we did not specify the summation bounds.
Similarly, we take the numerical electromagnetic fields $\vecEN$ and
$\vecBN$ in the finite-dimensional space $\calVN$.
Hence, they admit the expansions
\begin{align}
  \vecEN(\xv,t)
  &= \sum_{\Il}\vecEIl(t)\,\varphiIl(\xv),\qquad \forall\, t>0,\,\xv\in\Ox\,,
  \label{eq:Ev:expansion}
  \\[0.5em]
  \vecBN(\xv,t)
  &= \sum_{\Il}\vecBIl(t)\,\varphiIl(\xv),\qquad \forall\, t>0,\,\xv\in\Ox\,,
  \label{eq:Bv:expansion}
\end{align}
where $\vecEIl(t)=\big(\ExIl(t),\EyIl(t),\EzIl(t)\big)^T$ and
$\vecBIl(t)=\big(\BxIl(t),\ByIl(t),\BzIl(t)\big)^T$ are the \DG{}
expansion coefficients of the spatial components of the electric and
magnetic fields, $\vecEN=(\ExN,\EyN,\EzN)^T$ and
$\vecBN=(\BxN,\ByN,\BzN)$, respectively.

The semi-discrete variational formulation of the Hermite-\DG{} method
reads as: \emph{For every species $s$, and any time $t\in(0,T]$, find
  $\fsN(\cdot,\cdot,t)\in\calHN\times\calVN$ and
  $\vecEN(\cdot,t),\vecBN(\cdot,t)\in\calVN$ such that}
\begin{subequations}
  \label{eq:semi-discrete:Vlasov-Maxwell}
  \begin{align}
    \As\big( (\fsN,\vecEN,\vecBN),(\Psi,\varphi)\big) &=0
    \phantom{\Ls(\varphi)\fsNz\vecENz\vecBNz}\hspace{-1cm}
    \forall\,(\Psi,\varphi)\in\calHpN\times\calVN,
    \label{eq:semi-discrete:Vlasov}\\[0.5em]
    \Bs\big( (\vecEN,\vecBN), \varphi \big) &= \Ls(\varphi)
    \phantom{0\fsNz\vecENz\vecBNz}\hspace{-1cm}
    \forall\,\varphi\in\calVN,
    \label{eq:semi-discrete:Maxwell}
    \\[0.5em]
    \fsN(\cdot,\cdot,0)&=\fsNz   \phantom{0\Ls(\varphi)\vecENz\vecBNz}\hspace{-1cm}\textrm{in~}\Ox\times\Ov\,,  \\[0.5em]
    \vecEN(\cdot,0)    &=\vecENz \phantom{0\Ls(\varphi)\fsNz\vecBNz}\hspace{-1cm}  \textrm{in~}\Ox\,,\\[0.5em]
    \vecBN(\cdot,0)    &=\vecBNz \phantom{0\Ls(\varphi)\fsNz\vecENz}\hspace{-1cm}  \textrm{in~}\Ox\,,
  \end{align}
  \emph{where $\fsNz$, $\vecENz$ and $\vecBNz$ are the orthogonal
    projections of the initial conditions $\fs(\cdot,\cdot,0)$, $\Ev(\cdot,0)$ and
    $\Bv(\cdot,0)$ onto the spaces $\calHN\times\calVN$, $\calVN$ and $\calVN$,
    respectively, and suitable boundary conditions are prescribed on $\partial\Ox$.}
\end{subequations}

\medskip 
To define the multilinear form $\As$ in~\eqref{eq:semi-discrete:Vlasov},
we first introduce the auxiliary vector function
\begin{align*}
  \gv^{s,N}_{\Psi}(\xv,t) := \int_{\Ov}\vv\,\fsN(\xv,\vv,t)\Psi(\xiv)\,d\xiv,
  \quad\forall\,\Psi\in\calHpN.
\end{align*}
Then, for any $\fsN\in\calHN\times\calVN$, $\vecEN\in\calVN$,
$\vecBN\in\calVN$, and $(\Psi,\varphi)\in\calHpN\times\calVN$, we
define
\begin{align}
  &\As\big( \big(\fsN,\vecEN,\vecBN\big),(\Psi,\varphi)\big)
  := \sum_{I}\left(
    \int_{I}\int_{\Ov}\PT{\fsN}\,\Psi(\xiv)\,\varphi(\xv)\,d\xiv\dxv 
    -\int_{I}\gv^{s,N}_{\Psi}\cdot\GRADX\varphi(\xv)\,\dxv
  \right.
  \nonumber\\[0.5em] 
  &\qquad\qquad
  \left.
    +\int_{\partial\Is}\NFLX{\nv\cdot\gv^{s,N}_{\Psi}}\,\varphi(\xv)\,\dS
    +\csVM
    \int_{I}\int_{\Ov}\big(\vecEN+\vv\times\vecBN)\cdot\GRADV\fsN\,\Psi(\xiv)\,\varphi(\xv)\,d\xiv\dxv
  \right).
  \label{eq:bilA:def}
\end{align}
Concerning the discretization of
Maxwell's equations in \eqref{eq:semi-discrete:Maxwell}, we first consider the fluxes of the conservative
formulation~\eqref{eq:Maxwell:conservative:00}, and observe that, if
we partition the vector flux $\FF(\UvC)$ in a column-wise form, so
that
\begin{align}
  \nabla_{\xv}\cdot\FF(\UvC) 
  = \dfrac{\partial}{\partial x}\Fvx(\UvC)
  +\dfrac{\partial}{\partial y}\Fvy(\UvC)
  +\dfrac{\partial}{\partial z}\Fvz(\UvC),
  \label{eq:FF:def}
\end{align}
we can write $\Fvx(\UvC)=\matF_{x}\UvC$, $\Fvy(\UvC)=\matF_{y}\UvC$,
and $\Fvz(\UvC)=\matF_{z}\UvC$, with flux matrices
$\matF_{x},\matF_{y},\matF_{z}\in\mathbb{R}^{6\times6}$
\revtwo{defined as in \cite[Eqs. (50)-(52)]{koshkarov2021}.}
Let $\nv=(\ns_x,\ns_y,\ns_z)^T$ be a generic vector in $\mathbb{R}^3$.
We use the notation $\FF(\UvC)\nv$ to indicate
$\FF(\UvC)\nv:=\Fvx(\UvC)\ns_x+\Fvy(\UvC)\ns_y+\Fvz(\UvC)\ns_z$,
and use an analogous definition for $\FFE(\UvC)\nv$ and
$\FFB(\UvC)\nv$.
\revtwo{Using \cite[Eq. (53)]{koshkarov2021}, one has
$\FFE(\UvC)\nv\cdot\vecE = \FFB(\UvC)\nv\cdot\vecB =\nv\cdot\big(\vecE\times\vecB\big)$.}

\begin{figure}[!ht]
  \begin{center}
    \begin{overpic}[width=10cm]{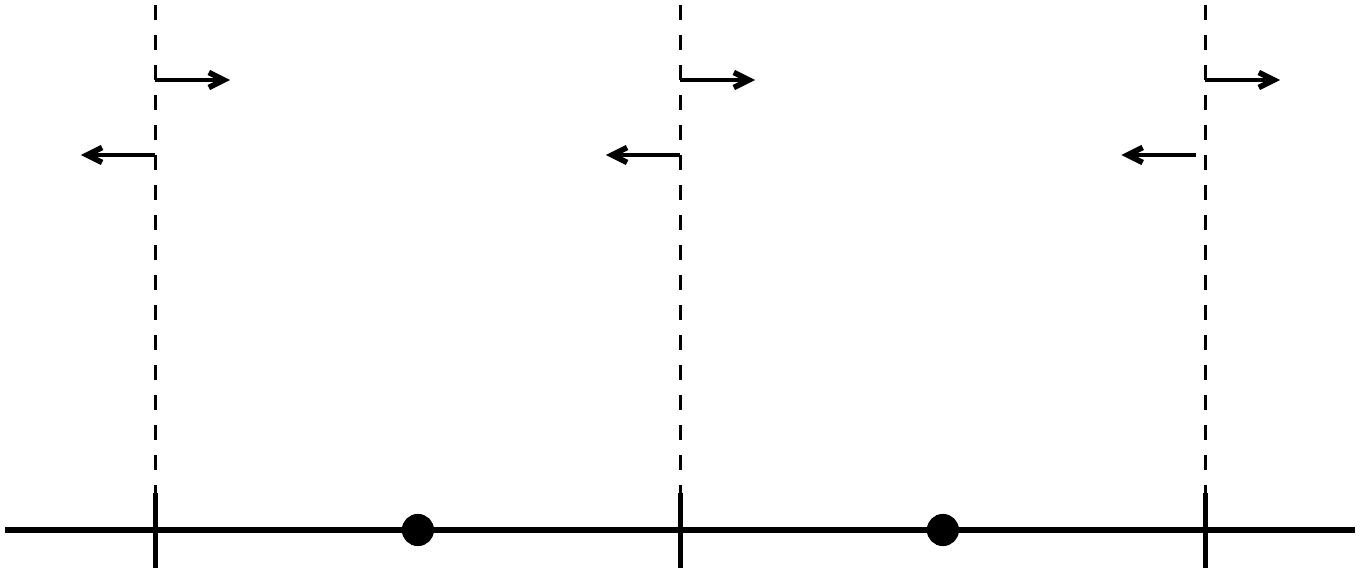}
      \put(85,-5){$i+\frac32$}
      \put(46,-5){$i+\frac12$}
      \put( 8, -5){$i-\frac12$}
      \put(66,-3){$i+1$}
      \put(30,-3){$i$}
      \put(76,36){$\mathbf{n}^{-}_{\F_{i+\frac32}}$}
      \put(92,40) {$\mathbf{n}^{+}_{\F_{i+\frac32}}$}
      \put(38,36){$\mathbf{n}^{-}_{\F_{i+\frac12}}$}
      \put(54,40){$\mathbf{n}^{+}_{\F_{i+\frac12}}$}
      \put(  0,36){$\mathbf{n}^{-}_{\F_{i-\frac12}}$}
      \put(16,40){$\mathbf{n}^{+}_{\F_{i-\frac12}}$}
      \put(42,22){$\matFx^{+}$}
      \put(53,22){$\matFx^{-}$}
      \put(39,10){${\Uv}^{-}_{|\F_{i+\frac12}}$}
      \put(53,10){${\Uv}^{+}_{|\F_{i+\frac12}}$}
    \end{overpic}
  \end{center}
  \caption{Notation of the numerical flux functions along direction
    $x$.
    The normal vectors $\nv^{+}_{\F_{\ell}}$ are oriented along the
    positive real axis (from left to right); vectors
    $\nv^{-}_{\F_{\ell}}$ are in the opposite sense.
    %
    At the interface $\F_{i+\frac12}$, the matrix
    $\matFx=\matFx^{+}+\matFx^{-}$ is decomposed into characteristics
    waves traveling from left to right, $\matFx^{+}$, which transport
    solution ${\Uv}^{-}$ from inside cell $i$ toward the cell
    interface, and from right to left, $\matFx^{-}$, which transport
    the solution ${\Uv}^{+}$ from inside cell $i+1$ towards the cell
    interface.
    Here ${\Uv}^{\pm}_{|\F_{i+\frac12}}$ are shortcut notations for
    $\Uv(x^{\pm}_{i+\frac{1}{2}},y,z,t)$.\protect\footnotemark
    }
  \label{fig:numflux:notation}
  \vspace{0.5cm}
\end{figure}

Let $\Uv(t) := ((\vecEN)^T,(\vecBN)^T)^T$ denote the \DG{} approximation
of the vector-valued function $\uv$
in~\eqref{eq:Maxwell:conservative:05}.
The bilinear form associated with the discretization of Maxwell's
equations is
\begin{align}
  \Bs\big( (\vecEN,\vecBN), \varphi \big)
  &:= \BsE\big( (\vecEN,\vecBN), \varphi \big) + \BsB\big( (\vecEN,\vecBN), \varphi \big),
  \qquad\qquad
  \forall\,
  \vecEN,\,
  \vecBN,\,
  \varphi\in\calVN,
  \label{eq:bilB:def}
\end{align}
where
\begin{align} 
  \BsE\big( (\vecEN,\vecBN), \varphi \big) 
  &=
  \sum_{I}\left(
    \int_{I}\dfrac{\partial\vecEN}{\partial t}\varphi(\xv)\,\dxv 
    -\int_{I}\FFE(\Uv)\GRADX\varphi(\xv)\,\dxv
    +\int_{\partial\Is}\NFLX{\FFE(\Uv)\nv}\,\varphi(\xv)\,\dS
  \right),
  \label{eq:semi-discrete:Maxwell:E}\\[1.em]
  \BsB\big( (\vecEN,\vecBN), \varphi \big) 
  &=
  \sum_{I}\left(
    \int_{I}\dfrac{\partial\vecBN}{\partial t}\varphi(\xv)\,\dxv 
    -\int_{I}\FFB(\Uv)\GRADX\varphi(\xv)\,\dxv
    +\int_{\partial\Is}\NFLX{\FFB(\Uv)\nv}\,\varphi(\xv)\,\dS
  \right).
  \label{eq:semi-discrete:Maxwell:B}
\end{align}
The linear functional $\Ls$ in~\eqref{eq:semi-discrete:Maxwell} is
defined, for all $\varphi\in\calVN$, as
\begin{equation*}
  \Ls(\varphi)
  := -\frac{\ope}{\oce}\sum_{I}\int_{I}\vecJN(\xv,t)\varphi(\xv)\,\dxv,
  \qquad
 \vecJN(\xv,t) = \sum_{s}\qs\int_{\Ov}\vv\fsN(\xv,\vv,t)\,d\vv.
\end{equation*}
The quantity $\NFLX{\nv\cdot\gv^{s,N}_{\Psi}}$ in~\eqref{eq:bilA:def}
and $\NFLX{\FFE(\Uv)\nv}$ and $\NFLX{\FFB(\Uv)\nv}$
in~\eqref{eq:semi-discrete:Maxwell:E}-\eqref{eq:semi-discrete:Maxwell:B}
are the numerical fluxes at the interfaces of the element boundaries.
At the boundaries of the domain $\Ox$ these quantities are defined in
accordance with the prescribed boundary conditions.
We consider two different kinds of numerical fluxes for Maxwell's
equations, the \emph{central numerical flux} and the \emph{upwind
numerical flux}, while the Vlasov equation is treated with the upwind
numerical flux.
\revtwo{We remark that the choice of the numerical flux plays a crucial role in the conservation properties of the resulting semi-discrete problem and, hence, in the construction of the numerical time integrators, as it will be shown in Section~\ref{sec:cons-general}.}
We provide a brief description of the numerical treatment of the
boundary integral in the $x$ direction, and, in particular, at the
face $\Fip$.
We refer to~\cite[Section~4]{koshkarov2021} for a detailed derivation.
Figure \ref{fig:numflux:notation} illustrates the meaning of the main
symbols that we adopt in the paper.
Let $\nv=(n_x,n_y,n_z)^T$ be the unit vector that is orthogonal to the
boundary $\partial I$ of $I$ and $\nv_{\F_{*}}$ the normal vector to
the face $\F_{*}$ indicated by a specific triple of indices, e.g.,
$\nv_{\F_{i+\frac12,j,k}}$.
Since, by construction, $n_x=\pm 1$, $n_y=n_z=0$ on the two faces
$\F_{i\pm\frac{1}{2},j,k}$, we obtain,
\begin{align*}
  \int_{\partial_x I}n_x\Fvx(\Uv)\,\varphiIl(\xv)\dS = 
  \int_{\F_{i+\frac{1}{2},j,k}}\Fvx(\Uv)\,\varphiIl(\xs_{i+\frac12},\ys,\zs)\dy\dz -
  \int_{\F_{i-\frac{1}{2},j,k}}\Fvx(\Uv)\,\varphiIl(\xs_{i-\frac12},\ys,\zs)\dy\dz.
\end{align*}

\footnotetext{Figure reprinted from \cite{koshkarov2021}, with permission from Elsevier.}

Let 
$\Uv(\xv^{\pm}_{i+\frac{1}{2},j,k},t)$ denote
$\lim_{\epsilon\to0}\Uv(\xv_{i+\frac{1}{2},j,k}\pm\epsilon\,\nv_{\F_{i+\frac12,j,k}},t)$.
Then,
\begin{description}
\item[$\bullet$] the \textit{central numerical flux} is given by the formula:
\begin{align}\label{eq:fluxC}
  \NFLX{\matF_x\Uv}(\xv_{i+\frac{1}{2},j,k},t) 
  =\matF_x\,\frac{1}{2}\big(\Uv(\xv^{-}_{i+\frac{1}{2},j,k},t)+\Uv(\xv^{+}_{i+\frac{1}{2},j,k},t)\big);
\end{align}

\item[$\bullet$] the \textit{upwind numerical flux} is based on the
  characteristic decomposition of the flux matrices:
  $\matF_{\beta}=\matF_{\beta}^{+}+\matF_{\beta}^{-}$, for each
  $\beta\in\{x,y,z\}$.
  We perform such decomposition numerically as described in
  \revtwo{\cite[Remark 3.4]{koshkarov2021}.}
  On face $\F_{i+\frac{1}{2},j,k}$, we approximate the flux integral
  using the upwind flux decomposition of the matrix $\matF_x$ as
  follows
  \begin{align}\label{eq:fluxU}
    \NFLX{\matF_x\Uv}(\xv_{i+\frac{1}{2},j,k},t) 
    = \underbrace{\matF_x^{+}\Uv(x^{-}_{i+\frac12},y,z,t)}_{\textrm{from cell $I=I_{i,j,k}$}}
    + \underbrace{\matF_x^{-}\Uv(x^{+}_{i+\frac12},y,z,t)}_{\textrm{from cell $I^{+}=I_{i+1,j,k}$}}.
  \end{align}
\end{description}


\begin{remark}[Consistency]
  Since upwind and central numerical fluxes are consistent whenever
  the associated physical flux is sufficiently regular, the Hermite-DG
  discretization \eqref{eq:semi-discrete:Vlasov-Maxwell} is
  consistent.
  This means that, if the exact solution fields $\fs$, $\vecE$,
  $\vecB$, and $\vecJ$ are sufficiently regular (at least continuous),
  then they satisfy the variational problem
  \begin{align*}
    \As\big( \big(\fs,\Ev,\Bv\big),(\Psi,\varphi)\big) &=0
    \qquad\forall\,(\Psi,\varphi)\in\calHpN\times\calVN,\\[0.5em]
    \Bs\big( \big(\Ev,\Bv\big), \varphi \big) &= \Ls(\varphi)
    \qquad\forall\,\varphi\in\calVN,
  \end{align*}
  defined in the finite dimensional spaces $\calHN$, $\calHpN$, and
  $\calVN$.
\end{remark}


%% file: sec5_invariants.tex
\section{Semi-discrete problem and conservation properties of the Hermite-DG scheme}\label{sec:invariants}

With the Hermite-DG
discretization~\eqref{eq:semi-discrete:Vlasov}-\eqref{eq:semi-discrete:Maxwell}
of the Vlasov-Maxwell equations, we can derive the semi-discrete
formulation of the evolution problem \eqref{eq:PDE}.
Analogously to~\eqref{eq:yvecC}, we introduce the finite-dimensional
vector of unknowns defined, for any $t\in[0,T]$, as
\begin{align}\label{eq:yvec}
  \ys(t) = 
  \left(
    \begin{array}{c}
      \big(\C_{n,m,p}^{s,\Il}(t)\big)_{n,m,p}^{s,\Il}\\[0.5em]
      \big(\vecEIl(t)\big)_{\Il}\\[0.5em]
      \big(\vecBIl(t)\big)_{\Il}
    \end{array}
  \right)\in\mathbb{R}^M,
\end{align}
where the coefficients $\C_{n,m,p}^{s,\Il}(t)$ are defined
in~\eqref{eq:fsN:def}, the coefficients $\vecEIl(t)$
in~\eqref{eq:Ev:expansion}, and the coefficients $\vecBIl(t)$
in~\eqref{eq:Bv:expansion}.
Testing
equations~\eqref{eq:semi-discrete:Vlasov}-\eqref{eq:semi-discrete:Maxwell}
against the basis functions $\{\Psi^{n,m,p}\}$ in $\calHpN$ and
$\{\phi^{\Il}\}$ in $\calVN$ yields an autonomous system of ordinary
differential equations that can formally be written as
\begin{align}\label{eq:IBVP:def}
  \ys'(t) = \FNs(\ys(t);t)
  \qquad\mbox{in}\quad (0,T],
\end{align}
where $\FNs$ corresponds to the Hermite-DG discretization of
\eqref{eq:rhs}, and a suitable initial condition $\ys(0)$ is provided.
We assume
that, for a fixed initial condition $\ys(0)$, the solution $\ys(t)$ at time $t\in(0,T]$ exists and is unique.

The fully discrete approximation of the Vlasov-Maxwell equations is
then obtained by solving \eqref{eq:IBVP:def} with one of the RK
methods introduced in Section \ref{sec:Runge-Kutta:methods}, after
setting $\ys^0=\ys(0)$.

\subsection{Invariants and conservation properties}
%

In order to study the conservation properties of the fully discrete
problem, we consider polynomial conserved quantities of a dynamical
system.
Suppose that $\calI\in C^{\infty}(\mathbb{R}^M)$ is an invariant of
the autonomous ODE \eqref{eq:IBVP:def}.
Then, $\calI$ is called a \emph{polynomial invariant} of degree $p$ if
it is a polynomial of degree at most $p$ in the variable
$y\in\mathbb{R}^M$, namely $\calI\in\mathbb{P}_p(\mathbb{R}^M)$.
For $p = 1$ we call $\calI$ a linear invariant and it is of the form
$\calI(y) = \mu^\top y + c$, where $\mu\in\mathbb{R}^M$ and
$c\in\mathbb{R}$.
For $p = 2$, we have a quadratic invariant that can be written as
$\calI(y) = y^\top S y + \mu^\top y + c$, for $S\in\mathbb{R}^{M\times
  M}$ symmetric, $\mu\in\mathbb{R}^M$ and $c\in\mathbb{R}$.


The conservation laws of general physical systems might not be
expressed as invariants of motion in the classical sense, since they
might need to take into account the exchanges of the system with the
external ambient space.
The flux of a quantity at the boundary of the domain contributes to a
conservation law as an integral over the temporal interval that
expresses the balance between what has entered and what has left the
domain during the system evolution.
In this case, we talk about conservation laws \emph{in weak form}, as
specified in the following definition.

\medskip
\begin{definition}\label{def:weak-inv}
  Let $\calD, \calC:\mathbb{R}^M\rightarrow\mathbb{R}$ be given
  functions of $\ys:[0,T]\rightarrow\mathbb{R}^M$.  A quantity of the
  form
  \begin{equation}\label{eq:I}
    \calI(\ys(t),t):=\calD(\ys(t)) + \int_{0}^t \calC(\ys(s))\,ds,
  \end{equation}
  provides a conservation law in weak form if it satisfies
  \begin{equation}\label{eq:dIdt}
    \dfrac{d\, \calI(\ys(t),t)}{dt} = 0,\qquad\forall\,y(t)\in\mathbb{R}^M.
  \end{equation}
\end{definition}

Notice that if $\calC$ in \eqref{eq:I} is identically zero, then
$\calI(\ys(t),t)=\calD(\ys(t))$ is an invariant of motion of
\eqref{eq:IBVP:def} in the classical sense.

\medskip
The Vlasov-Maxwell equations possess an infinite set of conserved
quantities in weak form, among those, the number of particles, the
total momentum, and the total energy. In the notation of
Definition~\ref{def:weak-inv}, these quantities are encoded in the
function $\calD$.
The term $\calC$, instead, describes the flux of the conserved
quantity at the boundaries of the computational domain and it is
associated with the possible presence of sources.
The term $\calC$ also takes into account the numerical contributions
to the conserved quantity $\calD$ that result from the numerical
approximation of the Vlasov-Maxwell equations in space and velocity.

In the next subsections, we study the conservation laws associated
with the total number of particles and the total energy for the
Vlasov-Maxwell system \eqref{eq:Vlasov:fs}, \eqref{eq:dEdt} and
\eqref{eq:dBdt}.
Note that we do not consider the conservation of total momentum since
this quantity is already not conserved at the level of the
semi-discrete system, as demonstrated in Ref.~\cite{koshkarov2021}.
For each of these, we characterize the corresponding conserved
quantity of the ODE \eqref{eq:IBVP:def}, obtained from the
approximation of the continuous conservation laws.
Finally, in Section~\ref{sec:cons-general}, we analyze the capability
of each RK temporal integrator introduced in Section
\ref{sec:Runge-Kutta:methods} to preserve such quantities.

The conservation properties of the semi-discrete scheme
\eqref{eq:semi-discrete:Vlasov-Maxwell} were already studied in
\cite{koshkarov2021}.
Here we characterize them as conservation properties of
\eqref{eq:IBVP:def}, in the sense of Definition~\ref{def:weak-inv},
and we refer the reader to \cite{koshkarov2021} for further details.
To fix the notation, for the generic quantity $\Phi^N(\xv,\vv,t)$
involved in the discretization of the Vlasov-Maxwell
equations~\eqref{eq:semi-discrete:Vlasov}-\eqref{eq:semi-discrete:Maxwell},
we use the superscript $\tau$ to indicate its approximation at time
$t^\tau$, namely $\Phi^{N,\tau}(\xv,\vv) \approx
\Phi^N(\xv,\vv,t^\tau)$.

\subsection{Number of particles}

\begin{lemma}\label{theo:number-of-particles}
  Let $\fsN(\xv,\vv,t)$ be the numerical solution of the semi-discrete
  Vlasov-Maxwell problem \eqref{eq:semi-discrete:Vlasov-Maxwell}
  and let $y(t)$ be as in \eqref{eq:yvec} for $t\in[0,T]$.
  The approximate number of particles of each plasma species $s$,
  \begin{equation}
    \cN^{s}(y(t)) = \sum_{I}\int_{I}\int_{\Ov}\fsN(\xv,\vv,t)\dxv\dvv,
    \label{eq:number-of-particles}
  \end{equation}
  satisfies
  \begin{equation}\label{eq:particles-balance}
    \dfrac{d \cN^{s}(y(t))}{dt} +\int_{\partial \Ox}\int_{\Ov}\big(\nv\cdot\vv\big)\,\fsN(\xv,\vv,t)\dvv \,\dS=0,
    \qquad\forall\,\ss.
  \end{equation}
\end{lemma}

\medskip
An analogous result holds for the total number of particles in the
plasma $\cN^{tot} := \sum_{s}\cN^{s}$.
Notice that the integral term
in~\eqref{eq:particles-balance} is the flux of particles at the
boundary of $\Ox$ and represents the balance between the number of
particles entering the domain and the particles leaving the domain at
time $t$.
This quantity vanishes when periodic boundary conditions are imposed
or when the plasma is entirely confined to $\Ox$.


\medskip
The approximate number of particles $\cN^s$ of a given plasma species
$s$ is a linear function of the vector $\ys$ whose components are
specified in \eqref{eq:yvec}.
Indeed, the Hermite-DG expansion \eqref{eq:fsN:def} of $\fsN$ yields
\begin{equation*}
  \cN^s(\ys) = \sum_{I}\int_{I}\int_{\Ov}\fsN(\xv,\vv,t)\dxv\dvv
  = \Delta x\Delta y\Delta z\,
    \alpha_x^s\alpha_y^s\alpha_z^s\sum_I C_{0,0,0}^{s,I,0}.
\end{equation*}
\subsection{Total energy}\label{sec:En}
In the continuum, the total energy of the Vlasov-Maxwell problem is defined as
\begin{equation*}
  \calE := \dfrac12 \sum_{s}\ms\int_{\Ox}
  \int_{\Ov}\ABS{\vv}^2\fs(\xv,\vv,t)\,d\vv\dxv
  + \dfrac12\left(\frac{\oce}{\ope}\right)^2 \int_{\Ox}\Big(\ABS{\vecE(\xv,t)}^2+\ABS{\vecB(\xv,t)}^2\Big)\dxv.
\end{equation*}

To derive a conservation law for \eqref{eq:IBVP:def} associated with
the total energy we consider, for $\ys$ defined as in \eqref{eq:yvec},
the following quantities,
\begin{align}
  \calE_{\KIN}(\ys) & := \dfrac12 \sum_{s}\ms\sum_{\Is}\int_{\Is}
  \int_{\Ov}\ABS{\vv}^2\fsN(\xv,\vv,t)\,d\vv\,\dxv,\label{eq:calEkin}\\
  \calE_{\ELE}(\ys) & :=\dfrac12\left(\frac{\oce}{\ope}\right)^2 \sum_{\Is}\int_{\Is}\Big(\ABS{\vecEN(\xv,t)}^2+\ABS{\vecBN(\xv,t)}^2\Big)\dxv,\label{eq:calEBE}\\
  \calE_{\JMP}(\ys) &:= \left(\dfrac{\oce}{\ope}\right)^2 \sum_{\F}\JMPF(\ys(t)),
  \label{eq:calE:jmp:def:cont}\\
  \calE_{\BND}(\ys) &:= \dfrac12  \sum_{s}\ms
  \int_{\partial\Ox}
  \int_{\Ov}\big(\nv\cdot\vv\big)\,\ABS{\vv}^2\fsN(\xv,\vv,t)\,d\vv\,\dS.\label{eq:calE:bnd}
\end{align}
The first two quantities $\calE_{\KIN}$ and $\calE_{\ELE}$ are the
approximate kinetic and electromagnetic energy.
The remaining quantities \eqref{eq:calE:jmp:def:cont} and
\eqref{eq:calE:bnd} are associated with the energy flux at the domain
boundaries and mesh interfaces.
Specifically, term \eqref{eq:calE:bnd} describes the exchange of
kinetic energy at the boundaries of the domain $\Ox$ at time $t$,
while term $\JMPF$ in \eqref{eq:calE:jmp:def:cont} is defined as
\begin{align}\label{eq:JF}
  \JMPF(\ys(t)) :=
  \begin{cases}
    0 
    & \qquad\textrm{for the central numerical flux \eqref{eq:fluxC}},\\[0.5em]
    \displaystyle\frac{1}{2}\int_{\F}\JUMP{\Uv(t)}_{\F}^\top\cdot\ABS{\FF}\cdot\JUMP{\Uv(t)}_{\F}\,\dS
    & \qquad\textrm{for the upwind numerical flux \eqref{eq:fluxU}},
  \end{cases}
\end{align}
where $\ABS{\FF}=\ABS{\matF_x}+\ABS{\matF_y}+\ABS{\matF_z}$ with
$\ABS{\matF_{\beta}}=\matF^+_{\beta}-\matF^-_{\beta}$, for
$\beta\in\{x,y,z\}$, and
\begin{align*}
  \JUMP{\Uv(t)}_{\F}
  = \nv^+\cdot\Uv(\xv^+_{\F},t)+\nv^-\cdot\Uv(\xv^-_{\F},t)
  =: \nv^+\cdot\Uv^+_{\F}(t)+\nv^-\cdot\Uv^-_{\F}(t),
\end{align*}
is the jump of $\Uv$ across the mesh face $\F$, and
$\Uv^{\pm}_{\F}(t)=\Uv(\xv^{\pm}_{\F},t)$ denotes the trace of $\Uv$
on opposite sides of $\F$.
Note that Eq. (\ref{eq:calE:jmp:def:cont}) depends only on the jumps
at the mesh interfaces of the electromagnetic field.
The fluxes at the mesh interfaces associated with the Vlasov equation
do not enter the time derivative of the kinetic energy (since it is a
linear quantity and one can take advantage of properties of telescopic
sums).
Hence, the fact that we use upwind fluxes for the Vlasov equation does
not affect explicitly the conservation of energy.
Note, also, that the energy flux \eqref{eq:calE:bnd} at the domain
boundaries vanishes whenever periodic boundary conditions are
prescribed.

The following result was shown in \cite[Theorem 5.3]{koshkarov2021}
under the assumption of periodic boundary conditions on $\partial\Ox$,
and we report it here, without proof, in the more general setting with
arbitrary boundary conditions.

\medskip 
\begin{lemma}
  \label{theo:energy}
  Let $\ys(t)$ be the numerical solution at time $t\in(0,T]$ of the
  semi-discrete Vlasov-Maxwell problem \eqref{eq:IBVP:def} with
  initial condition $\ys(0)$.
  Let $\calE_{\KIN}$ and $\calE_{\ELE}$ be the discrete kinetic and
  electromagnetic energies introduced in \eqref{eq:calEkin} and
  \eqref{eq:calEBE}, respectively.
  The variation in time of the total discrete energy of the system
  defined as
  \begin{equation}\label{eq:energyN:def}
    \calE_{\TOT}(\ys) = \calE_{\KIN}(\ys) + \calE_{\ELE}(\ys),
  \end{equation}
  satisfies
  \begin{equation}\label{eq:law-energy}
    \frac{d \calE_{\TOT}(\ys(t))}{\dt} 
    + \dfrac12  \sum_{s}\ms
    \int_{\partial\Ox}
    \int_{\Ov}\big(\nv\cdot\vv\big)\,\ABS{\vv}^2\fsN(\xv,\vv,t)\,d\vv\,\dS
    = -\left(\dfrac{\oce}{\ope}\right)^2\sum_{\F}\JMPF(y(t))\leq 0,
  \end{equation}
  where $\JMPF$ is defined in \eqref{eq:JF}.
\end{lemma}
In other words, the quantity
\begin{equation}
  \calE_{\MOD}(\ys(t)) := \calE_{\KIN}(\ys(t)) + \calE_{\ELE}(\ys(t)) +
  \int_{0}^{t}\calE_{\BND}(\ys(\sigma))\,\ds + \int_{0}^{t}\calE_{\JMP}(\ys(\sigma))\,\ds,
  \label{eq:Emodcont}
\end{equation}
is conserved in the sense of \eqref{eq:I}, \eqref{eq:dIdt} for the
evolution problem \eqref{eq:IBVP:def}.
When central numerical fluxes are used in Maxwell's equations, the
last term vanishes and $\calE_{\MOD}$ coincides with the semi-discrete
total energy which is then conserved.
Instead, the term $\calE_{\JMP}$, appearing when upwind numerical
fluxes are used in Maxwell's equations, corresponds to a numerical
dissipation and hence, in this case, the semi-discrete total energy is
not conserved.

\medskip
The total energy \eqref{eq:energyN:def} is a quadratic function of
the vector unknown $\ys$ whose components are specified in \eqref{eq:yvec}.
Indeed, the kinetic energy can be expressed in the components of $\ys$ as
\begin{equation*}
  \begin{aligned}
    \calE_{\KIN}(\ys)=
    &\;\dfrac14\Delta x\Delta y\Delta z\,
    \sum_s\alpha_x^s\alpha_y^s\alpha_z^s\,
    \ms \sum_I
    \bigg[
      \sum_{\beta\in\{x,y,z\}}
      \big((\alpha^s_{\beta})^2+ 2(u^s_{\beta})^2\big)\,
      C^{s,I,0}_{0,0,0}+
      2\sqrt{2}\alpha^s_x u^s_x\,C^{s,I,0}_{1,0,0}\\
      & + (\alpha^s_x)^2\, C^{s,I,0}_{2,0,0} +
      2\sqrt{2}\alpha^s_y u^s_y\,C^{s,I,0}_{0,1,0}+
      (\alpha^s_y)^2\,C^{s,I,0}_{0,2,0}+
      2\sqrt{2}\alpha^s_z u^s_z\,C^{s,I,0}_{0,0,1}+
      (\alpha^s_z)^2\,C^{s,I,0}_{0,0,2}\bigg].
  \end{aligned}
\end{equation*}
The kinetic energy is linear in $\ys$ so that
$\calE_{\KIN}(\ys)=\mu^\top\ys$ where $\mu\in\mathbb{R}^M$ can be
derived from the expression above by using a suitable application that
maps the indices $(s,I,l,n,m,p)$ into a single integer index in
$[1,M]$.
Concerning the electromagnetic energy, the orthogonality properties of
the Legendre polynomials yield
\begin{equation*}
  \calE_{\ELE}(\ys) = \dfrac12\left(\dfrac{\oce}{\ope}\right)^2
  \Delta x\Delta y\Delta z
  \sum_{I,l}\big[(E^{\Il})^2 + (B^{\Il})^2\big].
\end{equation*}
Hence, $\calE_{\ELE}(\ys) =\ys^\top S\ys$ where
$S\in\mathbb{R}^{M\times M}$ is the zero matrix in the part that
multiplies the components $(C_{n,m,p}^{s,\Il})$ of $\ys$ and it is the
identity matrix in the part that multiplies the DG degrees of freedom
of the electromagnetic fields.


%% file: sec6_cons_GL_RK.tex

\section{Conservation properties of the fully discrete Vlasov-Maxwell equations}
\label{sec:cons-general}

This section pertains to the conservation of invariants and balance
laws in weak form by the numerical time integrators introduced
in Section~\ref{sec:Runge-Kutta:methods}, when applied to the
semi-discrete problem \eqref{eq:IBVP:def} ensuing from the Hermite-DG
discretization \eqref{eq:semi-discrete:Vlasov-Maxwell} of the
Vlasov-Maxwell equations.

\medskip
\begin{definition}
  A numerical time integrator preserves polynomial invariants of
  degree $p$ if, for any autonomous ODE of the form
  \eqref{eq:IBVP:def} with an invariant
  $\calI\in\mathbb{P}_p(\mathbb{R}^M)$, the function $\calI$ satisfies
  \begin{equation*}
    \calI(y^\tau) = \calI(y^0),\qquad \forall\, \tau,
  \end{equation*}
  where $y^\tau$ is the approximate solution of $\eqref{eq:IBVP:def}$
  at time $t^\tau$ with initial condition $y^0\in\mathbb{R}^M$.
\end{definition}

\medskip
A complete characterization of Runge-Kutta methods that preserve
polynomial invariants can be found in \cite[Chapter IV]{HLW06}.
In particular, all explicit and implicit RK methods preserve linear
invariants, and all GL methods preserve quadratic invariants,
\cite{Cooper87}.

In the following result, we characterize the conservation properties
of Runge-Kutta methods when dealing with conservation laws in weak
form, in the sense of Definition~\ref{def:weak-inv}.

\medskip
\begin{theorem}\label{thm:conserved-quantities}
  Let $\calD, \calC:\mathbb{R}^M\rightarrow\mathbb{R}$ be given
  functions of $\ys:[0,T]\rightarrow\mathbb{R}^M$.
  Assume that the quantity
  \begin{equation*}
    \calI(\ys(t),t):=\calD(\ys(t)) + \int_{0}^t \calC(\ys(s))\,ds,
  \end{equation*}
  satisfies a conservation law in weak form in the sense of
  Definition~\ref{def:weak-inv}, namely
  \begin{equation}\label{eq:dIdt1}
    \dfrac{d\, \calI(\ys(t),t)}{dt} = 
    \dfrac{d\, \calD(\ys(t))}{dt} + \calC(\ys(t)) = 0,
  \end{equation}
  for every $y(t)\in\mathbb{R}^M$.
  Consider the two cases:
  \begin{enumerate}
  \item $\{\ys^{\tau}\}_{\tau=1,\ldots,T/\Delta\ts}$ is the sequence
    of numerical solutions of \eqref{eq:IBVP:def} with initial
    condition $\ys^0$, obtained with an explicit RK scheme of the
    family \eqref{eq:RK}.  Assume, then, that $\calI$ is linear in
    $\ys$ and let $\Ng=\NRK$.
  \item $\{\ys^{\tau}\}_{\tau=1,\ldots,T/\Delta\ts}$ is the sequence
    of numerical solutions of \eqref{eq:IBVP:def} with initial
    condition $\ys^0$, obtained with a GL scheme of the family
    \eqref{eq:GL}.
    Assume, in this case, that $\calI$ is at most quadratic in $\ys$
    and let $\Ng=\NGL$.
  \end{enumerate}
  Under the assumptions above, it holds that
  \begin{equation}\label{eq:approxInv}
    \calD(\ys^{\tau+1}) = \calD(\ys^{\tau}) -
    \Delta t \sum_{i=1}^{\Ng} b_i \,\calC(\calY_i),\qquad\forall\,\tau,
  \end{equation}
  where $(\calY_i)_{i=1}^{\Ng}$ are defined in \eqref{eq:RK} for case
  (i) and in \eqref{eq:GL} for case (ii). Equation (\ref{eq:approxInv}) is the discrete equivalent of Eq. (\ref{eq:dIdt1}) and expresses the conservation law in weak form for the family of RK methods.
\end{theorem}
\medskip
\BEGINPROOF
Since the quantity $\calI$, and hence $\calD$, is assumed to be at
most quadratic in $\ys$ we can recast it as $\calD(\ys)=\ys^\top S\ys
+ \mu^\top \ys + \eta$, for $\mu\in\mathbb{R}^M$, $\eta\in\mathbb{R}$
and a symmetric matrix $S\in\mathbb{R}^{M\times M}$ that is
identically zero when $\calI$ is linear, as assumed in case (i).
Using the definition of temporal integrator in \eqref{eq:RK} and
\eqref{eq:GL}, we have $\ys^{\tau+1}=\ys^{\tau}+ \delta\ys^{\tau}$
with
\begin{equation}\label{eq:deltay}
\delta\ys^{\tau}:= \Delta t\sum_{i=1}^{\Ng} b_i \FN{i},
\end{equation}
with $\FN{i}:=\FNs(\calY_i;t^{\tau}+c_i\Delta t^{\tau})$,
and $\Ng$ either equal to $\NRK$ or $\NGL$.
Then,
\begin{equation*}
  \calD(\ys^{\tau+1}) =(\ys^{\tau+1})^\top S\ys^{\tau+1} + \mu^\top \ys^{\tau+1} + \eta
  = \calD(\ys^{\tau}) + 2(\ys^{\tau})^\top S\delta\ys^{\tau}
  + (\delta\ys^{\tau})^\top S\delta\ys^{\tau} + \mu^\top\delta\ys^{\tau},
\end{equation*}
where, due to the symmetry of S, we have used $(\ys^{\tau})^\top S\delta\ys^{\tau}=(\delta\ys^{\tau})^\top S\ys^{\tau}$.
The definition of $\delta\ys^{\tau}$ yields
  \begin{equation}\label{eq:invD1b}
\calD(\ys^{\tau+1}) = \calD(\ys^{\tau}) 
  + 2\Delta t\sum_{i=1}^{\Ng} b_i (\ys^{\tau})^\top S\FN{i}
  + \Delta t^2 \sum_{i,j=1}^{\Ng} b_i b_j (\FN{i})^\top S\FN{j}
  + \Delta t \sum_{i=1}^{\Ng} b_i \mu^\top \FN{i}.
\end{equation}
Using the formula $\calY_i=\ys^{\tau}+\Delta t\sum_{j=1}^{\Ng}
a_{ij}\FN{j}$ to substitute $\ys^{\tau}$ in \eqref{eq:invD1b},
results in
\begin{equation}\label{eq:invD2}
  \calD(\ys^{\tau+1}) = \calD(\ys^{\tau})
  + \Delta t\sum_{i=1}^{\Ng} b_i (\mu^\top + 2\calY_i^\top S)\FN{i}
  + \Delta t^2 \sum_{i,j=1}^{\Ng} (b_i b_j-a_{ij}b_i-a_{ji}b_j) (\FN{i})^\top S\FN{j}.
\end{equation}
Here, when $\Ng=\NRK$, $a_{ij}=0$ for any $j\geq i$.
Assumption \eqref{eq:dIdt1}, together with \eqref{eq:IBVP:def}, yield
\begin{equation}\label{eq:dDdt}
  -\calC(\ys(t))= \dfrac{d \calD(\ys(t))}{dt} = \dfrac{d}{dt} (\ys^\top S \ys + \mu^\top \ys) = (\mu^\top+2\ys^\top S) \dfrac{d\ys}{dt} = (\mu^\top + 2\ys^\top S) \FNs(\ys), \qquad\forall\, \ys.
\end{equation}
Taking $\ys=\calY_i$ in \eqref{eq:dDdt} and substituting in \eqref{eq:invD2}
results in
\begin{equation*}
  \calD(\ys^{\tau+1}) = \calD(\ys^{\tau})
  - \Delta t\sum_{i=1}^{\Ng} b_i \calC(\calY_i)
  + \Delta t^2 \sum_{i,j=1}^{\Ng} (b_i b_j-a_{ij}b_i-a_{ji}b_j) (\FN{i})^\top S\FN{j}.
\end{equation*}
For case (i), conclusion \eqref{eq:approxInv} follows from the fact that the matrix $S$ is identically zero.
For case (ii),
since any Gauss-Legendre method satisfies the algebraic constraint
\begin{equation}\label{eq:ab}
b_ib_j=a_{ij}b_i-a_{ji}b_j,\qquad \mbox{for}\quad i,j=1,\ldots, \NGL,
\end{equation}
the result \eqref{eq:approxInv} is proven.
\ENDPROOF

\medskip
\begin{remark}
  In Lemmas~\ref{theo:number-of-particles} and \ref{theo:energy} we
  have characterized the conservation of the total number of particles
  and of the total energy by isolating the contribution associated
  with the flux at the boundaries of the spatial domain $\Ox$.
  This term is associated with the presence of a source of particles
  or energy.
  Knowing the source allows one to quantify exactly the flux of
  particles or energy at the domain boundaries.
  In the case of known boundary fluxes, the last term of \eqref{eq:I}
  can be computed analytically and, hence, the conservation law
  \eqref{eq:dIdt} is satisfied exactly by the fully-discrete problem
  under the assumptions of Theorem~\ref{thm:conserved-quantities}.
\end{remark}

\medskip
In the forthcoming subsections, we assess the capability of each
family of temporal integrators of
Section~\ref{sec:Runge-Kutta:methods} to preserve the conservation
law~\eqref{eq:IBVP:def} stated in Section~\ref{sec:invariants}.

\subsection{Conservation properties of the GL-Hermite-DG scheme}
\label{sec:cons-GL}
Theorem~\ref{thm:conserved-quantities} implies that the discretization
of the Vlasov-Maxwell problem with the Hermite-DG method coupled with
any Gauss-Legendre method \eqref{eq:GL} yields the conservation of the
number of particles and the total energy, including the contribution
of the particles and energy flux at the domain boundaries, when
central numerical fluxes are used in the discretization of Maxwell's
equations.
This is because these two quantities are linear and quadratic in $y$,
respectively.
The numerical tests in Section \ref{sec:numerical:results} corroborate
this finding, when the implicit midpoint rule is used as temporal
integrator.

The conservation of total energy in the case of upwind numerical
fluxes needs to be treated separately.
Indeed the total energy $\calE_{\TOT}$ does not satisfy a conservation
law and can dissipate over time (depending on whether numerical
dissipation dominates over the sources).
As shown in Lemma~\ref{theo:energy}, if $\ys$ is the solution
of \eqref{eq:IBVP:def}, then
\begin{equation}\label{eq:diss}
  \dfrac{d\calE_{\TOT}(\ys(t))}{dt}
  + \calE_{\BND}(y(t)) \leq 0,\qquad\forall\, y(t)\in\mathbb{R}^M,
\end{equation}
which is a consequence of the DG spatial discretization.
In the following result we show that the Gauss-Legendre methods have
two properties: they preserve, at the temporal discrete level, the
dissipative relation \eqref{eq:diss}, and, secondly,
they do not introduce any additional dissipation, in the sense that
the amount of energy dissipated by the fully-discrete scheme equals
the temporal approximation of the dissipative
term \eqref{eq:calE:jmp:def:cont} ensuing from the Hermite-DG
semi-discretization.

\medskip
\begin{proposition}\label{prop:dissipation}
  Let us consider the semi-discrete problem \eqref{eq:IBVP:def}
  obtained with upwind numerical fluxes in the discretization of
  Maxwell's equations and with the initial condition $\ys(0)$.
  Let $\{\ys^{\tau}\}_{\tau=1,\ldots,T/\Delta\ts}$ be the sequence of
  approximate solutions provided by the Gauss-Legendre method
  \eqref{eq:GL}, with $\ys^0=\ys(0)$.
  Then, the total energy \eqref{eq:energyN:def} satisfies the inequality
  \begin{equation}\label{eq:diss-discrete}
    \calE_{\TOT}(\ys^{\tau+1}) - \calE_{\TOT}(\ys^{\tau}) +
    \Delta t \sum_{i=1}^{\NGL} b_i \,\calE_{\BND}(\calY_i)\leq 0,
  \end{equation}                 
  where $(\calY_i)_{i=1}^{\NGL}$ are defined in \eqref{eq:GL} and
  $(b_i)_{i=1}^{\NGL}$ are the weights of the Gauss-Legendre temporal
  scheme.
\end{proposition}
\medskip
\noindent
\BEGINPROOF
Since the total energy \eqref{eq:energyN:def} is a quadratic function
of $\ys$, as highlighted in Section~\ref{sec:En}, we can recast it as
$\calE_{\TOT}(\ys)=\ys^\top S\ys+\mu^\top\ys$, for
$\mu\in\mathbb{R}^M$ and $S\in\mathbb{R}^{M\times M}$ symmetric.
We can then repeat the derivation in the proof of
Theorem~\ref{thm:conserved-quantities}.
From \eqref{eq:invD2}, since any Gauss-Legendre method satisfies $b_i
b_j=a_{ij}b_i-a_{ji}b_j$, we have
\begin{equation*}
  \calE_{\TOT}(\ys^{\tau+1}) = \calE_{\TOT}(\ys^{\tau})
  + \Delta t\sum_{i=1}^{\NGL} b_i(\mu^\top +2 \calY_i^\top S)\FN{i},
\end{equation*}
where $(\calY_i)_{i=1}^{\NGL}$ are defined in \eqref{eq:GL}.
The conclusion \eqref{eq:diss-discrete} follows by showing that
\begin{equation}\label{eq:ineq1}
  \sum_{i=1}^{\NGL} b_i (\mu^\top +2 \calY_i^\top S)\FN{i} \leq 
  -\sum_{i=1}^{\NGL} b_i \,\calE_{\BND}(\calY_i).
\end{equation}
To prove this, we use \eqref{eq:diss} together with
\eqref{eq:IBVP:def}, so that
\begin{equation}\label{eq:ineq2}
  \begin{aligned}
    0\geq \dfrac{d \calE_{\TOT}(y(t))}{dt}+ \calE_{\BND}(\ys) & = \dfrac{d}{dt} (\ys^\top S \ys+\mu^\top\ys)
    + \calE_{\BND}(\ys) \\
    & = (\mu^\top+ 2\ys^\top S) \dfrac{d\ys}{dt} + \calE_{\BND}(\ys)=
    (\mu^\top+ 2\ys^\top S) \FNs(\ys)+ \calE_{\BND}(\ys), \qquad\forall\, \ys.
  \end{aligned}
\end{equation}
Taking $\ys=\calY_i$ in \eqref{eq:ineq2} and considering the fact that
each weight $b_i$ is non-negative, for all $i=1,\ldots,\NGL$ (see
e.g. \cite{IserlesBook}), yields \eqref{eq:ineq1} and the proof is
concluded.
\ENDPROOF

To quantify the amount of local numerical dissipation associated with
the GL-Hermite-DG discretization of \eqref{eq:PDE}, observe that, by
Lemma~\ref{theo:energy}, the quantity $\calE_{\MOD}$ in
\eqref{eq:Emodcont} defined as
\begin{equation*}
  \calE_{\MOD}(\ys(t)) = \calE_{\TOT}(\ys(t)) +
  \int_{0}^{t}\big(\calE_{\BND}(\ys(\sigma)) + \calE_{\JMP}(\ys(\sigma))\big)\,\ds,
\end{equation*}
is quadratic in $\ys$ and satisfies \eqref{eq:dIdt}.
Applying Theorem~\ref{thm:conserved-quantities} 
we can conclude that,
if $\ys^{\tau}$ is the solution of the GL-Hermite-DG discretization of
\eqref{eq:PDE}, then
\begin{equation}\label{eq:Etot-upwind}
  \calE_{\TOT}(\ys^{\tau+1}) = \calE_{\TOT}(\ys^{\tau}) -
  \Delta t \sum_{i=1}^{\NGL} b_i \,\calE_{\BND}(\calY_i) -
  \Delta t \sum_{i=1}^{\NGL} b_i \,\calE_{\JMP}(\calY_i),
\end{equation}
where $(\calY_i)_{i=1}^{\NGL}$ are defined in \eqref{eq:GL}.
Since the jump contribution in \eqref{eq:Etot-upwind} tends to zero as
the spatial and temporal meshes are refined, equation
\eqref{eq:Etot-upwind} entails that the total energy is asymptotically
preserved and the numerical dissipation is locally bounded.

\medskip
\begin{remark}
  An immediate consequence of \eqref{eq:Etot-upwind} is that, for the
  lowest order GL method, namely the implicit midpoint rule, it holds
  \begin{equation*}
    \begin{aligned}
      \dfrac{\calE_{\TOT}(\ys^{\tau+1}) - \calE_{\TOT}(\ys^{\tau})}{\Delta t} &\, +
      \dfrac12  \sum_{s}\ms
      \int_{\partial\Ox}
      \int_{\Ov}\big(\nv\cdot\vv\big)\,\ABS{\vv}^2
      \bigg(\dfrac{\fsNtp(\xv,\vv)+\fsNt(\xv,\vv)}{2}\bigg)\,d\vv\,\dS\\
      & = -\left(\dfrac{\oce}{\ope}\right)^2 \sum_{\F}\JMPF(\ys^{\tau+1/2})\leq 0.
    \end{aligned}
  \end{equation*}
\end{remark}

\subsection{Conservation properties of the RK-Hermite-\DG{} scheme}
\label{sec:fully-discrete:RK}

Theorem~\ref{thm:conserved-quantities} implies that explicit
Runge-Kutta schemes \eqref{eq:RK} applied to the semi-discrete
Vlasov-Maxwell problem \eqref{eq:IBVP:def} conserve the number of
particles.
Concerning quadratic invariants, however, general explicit Runge-Kutta
methods are known not to satisfy the algebraic constraint \eqref{eq:ab}.
This implies, in our case, that the conservation of the total energy
\eqref{eq:Emodcont} is not guaranteed by the scheme \eqref{eq:RK}, not
even when central numerical fluxes are used in the discretization of
Maxwell's equations.
This provides the motivation for the study of a modified version
of the explicit Runge-Kutta methods \eqref{eq:RK} with an additional
projection step that preserves the energy invariant for central
numerical fluxes and maintains the semi-discrete energy dissipation
(i.e., it does not add dissipation due to the time integration) for
the upwind numerical fluxes.
Taking the cue from the incremental direction technique \cite{Calvo-HernandezAbreu-Montijano-Randez:2006} and relaxation RK methods \cite{Ketcheson19},
this task is discussed in the next subsection.

\subsection{Energy-conserving explicit RK schemes}\label{ssec:modifiedRK}

Let us consider the $\NRK$-stage explicit RK method \eqref{eq:RK} defined by the
coefficients $(a_{ij})_{1\leq j<i \leq \NRK}$, ${(c_i)}_{i=1}^\NRK$ and weights
$(b_i)_{i=1}^{\NRK}$.  We use the term
\emph{modified} Runge-Kutta
scheme, to refer to the $\NRK$-stage explicit RK method having coefficients
$(a_{ij})_{1\leq j<i \leq \NRK}$, ${(c_i)}_{i=1}^\NRK$ and weights $(\bhati)_{i=1}^{\NRK}$,
where $\bhati:=\gamma^{\tau}b_i$ for any $1\leq i\leq\NRK$.
In the interval $(t^\tau,t^{\tau+1}]$, this scheme reads
\begin{equation}\label{eq:phiRKh:def}
\begin{aligned}
  \ys^{\tau+1} 
  = \phiRKh(\ys^{\tau}) 
  & := \ys^{\tau} + \gamma^{\tau}\Delta\ts\sum_{i=1}^{\NRK}b_{i}\FNs(\calY_i;t^{\tau}+c_i\Delta t),\\
      \calY_{1} &:= \ys^{\tau},\\
      \calY_{i} &:= \ys^{\tau} + \Delta\ts\sum_{j=1}^{i-1}a_{ij}\FNs(\calY_j;t^{\tau}+c_j\Delta t),
      \qquad i=2,\ldots,\NRK.
\end{aligned}
\end{equation}
The real scalar factors $\{\gamma^{\tau}\}_{\tau}$
depend on the approximate solution of the evolution problem, and are determined by imposing that a given conservation law is satisfied.
Observe that, with the definition of $\delta\ys^{\tau}$ given in
\eqref{eq:deltay}, we have that the standard RK scheme can be written
as $\ysb^{\tau+1}=\ys^{\tau}+\delta\ys^{\tau}$, while the modified RK
scheme reads $\ys^{\tau+1}=\ys^{\tau}+\gamma^{\tau}\delta\ys^{\tau}$
and, hence,
$\ys^{\tau+1}=\ys^{\tau}+\gamma^{\tau}(\ysb^{\tau+1}-\ys^{\tau})$.
This implies that the modified RK scheme can be interpreted as a
projection method \cite[Section IV.4]{HLW06} that consists in
projecting the current RK update onto the set of functions that
satisfy the invariant constraint.
Its implementation can be performed according to
Algorithm~\ref{alg:modified:RK}.
\medskip
\begin{algorithm}[!ht]
  \caption{Modified RK algorithm}\label{alg:modified:RK}
  \begin{algorithmic}[1]
    \Procedure{\texttt{Modified\_RK}}{$\ys^{\tau}$}                                                                                                 
    \smallskip\State $\ysb^{\tau+1}                  \gets \phiRK(\ys^{\tau})$                                                                            \label{algo:RK:update}
    \smallskip\State $\gamma^{\tau}                  \gets \texttt{\textsc{Compute\_gamma}}(\ys^{\tau},\ysb^{\tau+1})$                                             \label{algo:compute:gamma}              
    \smallskip\State $\ys^{\tau+1}                   \gets \ys^{\tau} + \gamma^{\tau}( \bar{\ys}^{\tau+1} - \ys^{\tau} )$ 
    \label{algo:projection:step}
    \smallskip\State \textbf{return} $\ys^{\tau+1}$
    \EndProcedure
  \end{algorithmic}
\end{algorithm}

In line~\ref{algo:RK:update} of Algorithm~\ref{alg:modified:RK}, we
update the solution $\ys^{\tau}$ to the intermediate solution
$\ysb^{\tau+1}$ by using the standard explicit RK scheme
\eqref{eq:RK}.
In line~\ref{algo:compute:gamma} we compute the scalar factor
$\gamma^{\tau}$ that is used in the projection step.  In
line~\ref{algo:projection:step} we perform the projection step
obtaining the final updated solution $\ys^{\tau+1}$ according to
\eqref{eq:phiRKh:def}.

In our setting, we want to ensure that the numerical solution
$\ys^{\tau+1}$ of the modified Runge-Kutta method at time $t^{\tau+1}$
satisfies the discrete relation corresponding to \eqref{eq:diss},
namely
\begin{equation}\label{eq:en-discrete}
  \calE_{\TOT}(\ys^{\tau+1})-\calE_{\TOT}(\ys^{\tau})
  + \Delta t \sum_{i=1}^{\NRK} \bhati \calE_{\BND}(\calY_i) = 
  - \varepsilon\, \Delta t \sum_{i=1}^{\NRK} \bhati\calE_{\JMP}(\calY_i),
\end{equation}
with the quantities defined in \eqref{eq:calEkin}, \eqref{eq:calEBE},
\eqref{eq:calE:jmp:def:cont} and \eqref{eq:calE:bnd}, and where
$\varepsilon\in\{0,1\}$ is a user-defined coefficient that acts as
follows.
The choice $\varepsilon=0$ entails that the modified RK scheme
conserves the discrete energy; while the choice $\varepsilon=1$ implies
that the solution of the modified RK scheme is constrained to maintain
the energy dissipation property of the semi-discrete formulation.
Observe that, when central numerical fluxes are used in the
discretization of Maxwell's equations \eqref{eq:bilB:def}, the term
$\calE_{\JMP}$ vanishes, and thus $\varepsilon$ does not play any role.

To determine the value of the scalar factor
$\gamma^{\tau}\in\mathbb{R}$ we impose the equality constraint
associated with the energy conservation \eqref{eq:en-discrete}.
For the sake of better readability, we introduce the notations
\begin{equation}\label{eq:deltaE}
  \delta\calE_{\JMP}^{\tau,\tau+1}:=\Delta t
  \sum_{i=1}^{\NRK} b_i \,\calE_{\JMP}(\calY_i),
  \qquad
  \delta\calE_{\BND}^{\tau,\tau+1}:=\Delta t
  \sum_{i=1}^{\NRK} b_i \,\calE_{\BND}(\calY_i).
\end{equation}
The following result extends to our setting the derivation of the incremental direction \cite[Theorem 2.1 (i)]{Calvo-HernandezAbreu-Montijano-Randez:2006} that guarantees energy conservation.
\begin{theorem}
  \label{theorem:modified:RK:discrete:energy:conservation}
  Let $\{\ys^{\tau}\}_{\tau=1,\ldots,T/\Delta\ts}$ be the sequence of
  numerical solutions provided by the modified RK method
  in~\eqref{eq:phiRKh:def}, with $\NRK\geq 2$ and initial condition
  $\ys^0$.
  Assume that the discrete total energy $\calE_{\TOT}$ and the
  electromagnetic energy $\calE_{\ELE}$, defined in
  \eqref{eq:energyN:def} and \eqref{eq:calEBE}, respectively, satisfy,
  for any $\tau$,
  \begin{equation}\label{eq:cond-inv}
    \Big(\calE_{\TOT}(\phiRK(\ys^{\tau})) - \calE_{\TOT}(\ys^{\tau})
    - \calE_{\ELE}(\phiRK(\ys^{\tau})-\ys^{\tau})
    + \delta\calE_{\BND}^{\tau,\tau+1} + \varepsilon\,\delta\calE_{\JMP}^{\tau,\tau+1}\Big)
    \calE_{\ELE}(\phiRK(\ys^{\tau})-\ys^{\tau})\neq 0,
  \end{equation}
  with $\delta\calE_{\BND}^{\tau,\tau+1}$ and
  $\delta\calE_{\JMP}^{\tau,\tau+1}$ as in \eqref{eq:deltaE}.
  Then, the conservation law \eqref{eq:en-discrete} associated with
  the total energy \eqref{eq:energyN:def} is satisfied with the choice
  of $\gamma^{\tau}\in\REAL$ given by
  \begin{align}\label{eq:gamma:def}
    \gamma^{\tau} = 1 - 
    \dfrac
        { \calE_{\TOT}( \phiRK(\ys^{\tau}) ) - \calE_{\TOT}(\ys^{\tau}) + \delta\calE_{\BND}^{\tau,\tau+1} + \varepsilon\,\delta\calE_{\JMP}^{\tau,\tau+1}}
        { \calE_{\ELE}(\phiRK(\ys^{\tau})-\ys^{\tau})},\qquad\forall\,\tau.
  \end{align}
\end{theorem}

\medskip
\noindent
\BEGINPROOF
Analogously to the proof of Proposition~\ref{prop:dissipation}, we
write the total energy \eqref{eq:energyN:def} as
$\calE_{\TOT}(\ys)=\ys^\top S\ys+\mu^\top\ys$, for
$\mu\in\mathbb{R}^M$ and $S\in\mathbb{R}^{M\times M}$ symmetric.
Separating the kinetic and electromagnetic energy contributions, we
can formally write $\calE_{\KIN}(\ys)=\mu^T\ys$ and
$\calE_{\ELE}(\ys)=\ys^\top S\ys$.
In addition, we introduce the notation $\big<a,b\big>:=a^\top S\,b$
for any $a, b\in\REAL^M$.
The goal is to find $\gamma^{\tau}$ in \eqref{eq:phiRKh:def} such that
the energy conservation~\eqref{eq:en-discrete} is satisfied.  With the
notations introduced above, we impose
\begin{align*}
  \mu^T\ys^{\tau+1} + \big<\ys^{\tau+1},\ys^{\tau+1}\big> - 
  \mu^T\ys^{\tau} - \big<\ys^{\tau},\ys^{\tau}\big> + \,\gamma^{\tau}
  \big(\delta\calE_{\BND}^{\tau,\tau+1} + \varepsilon\,\delta\calE_{\JMP}^{\tau,\tau+1}\big) = 0.
\end{align*}
Substituting expression \eqref{eq:phiRKh:def} for $\ys^{\tau+1}$,
and using the bilinearity and symmetry of the form
$\big<\cdot,\cdot\big>$, yield
\begin{equation}\label{eq:eq-gamma}
  \gamma^{\tau}\big(\mu^T\delta\ys^{\tau}+2\big<\ys^{\tau},\delta\ys^{\tau}\big>
  +\delta\calE_{\BND}^{\tau,\tau+1} + \varepsilon\,\delta\calE_{\JMP}^{\tau,\tau+1}\big) +
  (\gamma^{\tau})^2\big<\delta\ys^{\tau},\delta\ys^{\tau}\big> = 0.
\end{equation}
Such relation is trivially satisfied by $\gamma^{\tau}=0$
corresponding to the steady state solution $\ys^{\tau+1}=\ys^{\tau}$,
which clearly implies that the energy is the same at $\ts^{\tau+1}$
and $\ts^{\tau}$.
To obtain nontrivial solutions of \eqref{eq:eq-gamma}, we need to
assume that
\begin{equation}\label{eq:cond-gamma}
    \big(\mu^T\delta\ys^{\tau}+2\big<\ys^{\tau},\delta\ys^{\tau}\big>
  +\delta\calE_{\BND}^{\tau,\tau+1} + \varepsilon\,\delta\calE_{\JMP}^{\tau,\tau+1}\big) \big<\delta\ys^{\tau},\delta\ys^{\tau}\big> \neq 0.
\end{equation}
Under this condition, we obtain the nontrivial solution
\begin{align}
  \gamma^{\tau} 
  = -\dfrac{ \mu^T\delta\ys^{\tau}+2\big<\ys^{\tau},\delta\ys^{\tau}\big> +
    \delta\calE_{\BND}^{\tau,\tau+1} + \varepsilon\,\delta\calE_{\JMP}^{\tau,\tau+1}}{ \big<\delta\ys^{\tau},\delta\ys^{\tau}\big> }.
  \label{eq:gamma:00}
\end{align}
The proof is completed by reformulating $\gamma^{\tau}$ in terms of
the discrete plasma energies.
First, we note that, for $\overline{\ys}^{\tau+1}=\phiRK(\ys^{\tau})$ as in Algorithm~\ref{alg:modified:RK}, it holds that
\begin{align*}
  \big<\overline{\ys}^{\tau+1},\overline{\ys}^{\tau+1}\big> 
  = \big<\ys^{\tau},\ys^{\tau}\big> + 2\big<\ys^{\tau},\delta\ys^{\tau}\big>
  + \big<\delta\ys^{\tau},\delta\ys^{\tau}\big>.
\end{align*}
This allows us to rewrite the second term of the numerator
of~\eqref{eq:gamma:00} as follows,
\begin{align*}
  2\big<\ys^{\tau},\delta\ys^{\tau}\big> =
  \big<\overline{\ys}^{\tau+1},\overline{\ys}^{\tau+1}\big> -
  \big<\ys^{\tau},\ys^{\tau}\big> -
  \big<\delta\ys^{\tau},\delta\ys^{\tau}\big>.
\end{align*}
Hence, since $\delta\ys^{\tau}=\overline{\ys}^{\tau+1}-\ys^{\tau}$, we have
\begin{align*}
  \mu^T\delta\ys^{\tau} + 2\big<\ys^{\tau},\delta\ys^{\tau}\big> & =
  \big( \mu^T\overline{\ys}^{\tau+1} + \big<\overline{\ys}^{\tau+1},\overline{\ys}^{\tau+1}\big> \big) -
  \big( \mu^T\ys^{\tau}         + \big<\ys^{\tau},\ys^{\tau}\big>  \big) - 
  \big<\delta\ys^{\tau},\delta\ys^{\tau}\big>\\
  & = \calE_{\TOT}(\overline{\ys}^{\tau+1}) - \calE_{\TOT}(\ys^{\tau})
  - \calE_{\ELE}(\delta\ys^{\tau}).
\end{align*}
Therefore, equation \eqref{eq:gamma:00} becomes,
\begin{align*}
  \gamma^{\tau} = 1 -
  \dfrac{ 
    \calE_{\TOT}(\overline{\ys}^{\tau+1}) - \calE_{\TOT}(\ys^{\tau})
    + \delta\calE_{\BND}^{\tau,\tau+1} + \varepsilon\,\delta\calE_{\JMP}^{\tau,\tau+1} }{ \calE_{\ELE}(\delta\ys^{\tau})},
\end{align*}
and condition \eqref{eq:cond-gamma} equals \eqref{eq:cond-inv}.
\ENDPROOF

\medskip
The modified method \eqref{eq:phiRKh:def}, \eqref{eq:gamma:def} and
Theorem~\ref{theorem:modified:RK:discrete:energy:conservation} offer
the discrete counterpart of the result of Lemma~\ref{theo:energy}, for
explicit temporal integrators of the Runge-Kutta family.
Therefore, with the modified RK scheme, the discrete total energy is
exactly conserved when central numerical fluxes are used in the \DG{}
approximation of Maxwell's equations and, hence, $\gamma^\tau=1$.
When upwind fluxes are considered, two options are possible: the
discrete total energy is exactly preserved by the choice of
$\gamma^{\tau}$ in \eqref{eq:gamma:def} with $\varepsilon=0$; or a
dissipative term appears in the energy conservation equation -- for
$\varepsilon=1$ -- and this depends on a quadratic form of the cell
interface jumps of the electromagnetic fields and is expected to go to
zero with the accuracy of the DG scheme.

\medskip
In the numerical implementation, Algorithm~\ref{alg:modified:RK} is
combined with Algorithm~\ref{alg:compute:gamma} that computes the
scalar factor $\gamma^{\tau}\in\mathbb{R}$.
The computation of the correction factor $\gamma^{\tau}$ in each
temporal interval requires $\NRK$ evaluations of the terms
$\calE_{\BND}$ and $\calE_{\JMP}$, two evaluations of the total energy
at the RK update $\overline{\ys}^{\tau+1}$, and one evaluation of the
electromagnetic energy $\calE_{\ELE}$ at $\delta \ys^{\tau}$.
\begin{algorithm*}[!ht]
  \caption{Compute $\gamma$ in the Maxwell solver}\label{alg:compute:gamma}
  \begin{algorithmic}[1]
    \Procedure{\texttt{Compute\_gamma}}{$\ys^{\tau},\ysb^{\tau+1}$}
    \smallskip\State $\delta\calE_{\JMP}^{\tau,\tau+1}\gets\Delta t
    \sum_{i=1}^{\NRK} b_i \,\calE_{\JMP}(\calY_i)$
    \smallskip\State $\delta\calE_{\BND}^{\tau,\tau+1}\gets\Delta t
    \sum_{i=1}^{\NRK} b_i \,\calE_{\BND}(\calY_i)$
    \smallskip\State $\gamma^{\tau}                  \gets 1 - \dfrac{ \calE_{\TOT}(\overline{\ys}^{\tau+1})-\calE_{\TOT}(\ys^{\tau})
      +\delta\calE_{\BND}^{\tau,\tau+1} + \varepsilon\,\delta\calE_{\JMP}^{\tau,\tau+1}}{ \calE_{\ELE}(\overline{\ys}^{\tau+1}-\ys^{\tau}) }$
    \EndProcedure
  \end{algorithmic}
\end{algorithm*}

\medskip
\begin{remark}[Explicit Euler scheme]\label{rmk:explicitE}
  For the lowest order RK scheme \eqref{eq:RK}, namely explicit Euler
  method, the result of
  Theorem~\ref{theorem:modified:RK:discrete:energy:conservation} does
  not apply.
  This is due to the fact that condition \eqref{eq:cond-inv} is not
  satisfied.
  Indeed, in the explicit Euler scheme, $\delta\ys^{\tau}=\Delta
  t\FNs(\ys^{\tau})$ in \eqref{eq:deltay} and
  $\delta\calE_{\BND}^{\tau,\tau+1}+\delta\calE_{\JMP}^{\tau,\tau+1}=\Delta
  t \,(\calE_{\BND}(\ys^{\tau})+\calE_{\JMP}(\ys^{\tau}))$ in
  \eqref{eq:deltaE}. Hence, the first factor in condition
  \eqref{eq:cond-gamma} reads
  \begin{equation}\label{eq:cond-EE}
    \Delta t\big(\mu^T\FNs(\ys^{\tau}) +
    2\big<\ys^{\tau},\FNs(\ys^{\tau})\big> +
    \calE_{\BND}(\ys^{\tau})+\varepsilon\,\calE_{\JMP}(\ys^{\tau})\big).
  \end{equation}
  In view of Lemma~\ref{theo:energy}, the total energy
  $\calE_{\TOT}(\ys) = \mu^\top y+\big<\ys,\ys\big>$ satisfies, for
  any $\ys\in\mathbb{R}^M$, equation \eqref{eq:dDdt} with
  $\calD=\calE_{\TOT}$ and
  $\calC=\calE_{\BND}+\varepsilon\,\calE_{\JMP}$.  This implies that
  \begin{equation*}
    -\calE_{\BND}(\ys^{\tau})-\varepsilon\,\calE_{\JMP}(\ys^{\tau}) = 
    \mu^\top \FNs(\ys^{\tau}) +
    2\big<\ys^{\tau},\FNs(\ys^{\tau})\big>.
  \end{equation*}
  Hence, the term in \eqref{eq:cond-EE} vanishes and condition
  \eqref{eq:cond-gamma} cannot be satisfied.
\end{remark}

\medskip
\begin{remark}
  The modified Runge-Kutta method \eqref{eq:phiRKh:def} preserves
  linear conservation laws according to
  Theorem~\ref{thm:conserved-quantities}.
  Indeed, the modified scheme can still be expressed as linear
  combinations of explicit Euler steps, since it is an explicit RK
  method of the family \eqref{eq:RK} with weights
  $(\gamma^{\tau}b_i)_{i=1}^{\NRK}$ for any $\tau$.
\end{remark}

As proven in \cite[Theorem 2.1 (ii)]{Calvo-HernandezAbreu-Montijano-Randez:2006} the projection step introduced by the modified RK method causes the loss of one order of convergence if compared with the standard RK scheme.
However, if the numerical solution of the modified RK method at step $\tau+1$ is interpreted as an approximation of the exact solution at time $t^{\tau}+\gamma^{\tau}\Delta t$ the accuracy of the base scheme is recovered \cite[Theorem 2.7 and Corollary 2.10]{Ketcheson19}.
\rev{In the following result we show that our specific $\gamma^{\tau}$, defined in \eqref{eq:gamma:def}, satisfies the condition
$\gamma^\tau=1+\mathcal{O}(\Delta t^{p-1})$ for $\Delta t\rightarrow
0$.}

\begin{proposition}\label{prop:accuracy}
  Consider an explicit RK scheme \eqref{eq:RK} of order $p\geq 2$ for
  the numerical approximation of \eqref{eq:IBVP:def}.
  Let $\{\ys^{\tau}\}_{\tau=1,\ldots,T/\Delta t}$ be the sequence of
  numerical solutions of \eqref{eq:IBVP:def} obtained with the
  modified RK method \eqref{eq:phiRKh:def} with initial condition
  $\ys^0=\ys(0)$ and $\gamma^{\tau}$ defined as in \eqref{eq:gamma:def}.
  Under the assumptions of
  Theorem~\ref{theorem:modified:RK:discrete:energy:conservation}, if $\calE_{\ELE}(\ys^{\tau})\neq 0$ for all $\tau$, the
  modified RK scheme has order of accuracy
  \begin{itemize}
      \item[(i)]  $p-1$, if the numerical solution $y^{\tau+1}$ is interpreted as an approximation of $y(t^{\tau}+ \Delta t)$;
      \item[(ii)] $p$, if the numerical solution $y^{\tau+1}$ is interpreted as an approximation of $y(t^{\tau}+\gamma^{\tau} \Delta t)$.
  \end{itemize}
\end{proposition}
\medskip
\BEGINPROOF
For the sake of better readability, in the following derivation, we
use the notation $\calE:=\calE_{\BND}+\varepsilon\,\calE_{\JMP}$ and
$\delta\calE^{\tau,\tau+1}:=\delta\calE_{\BND}^{\tau,\tau+1}+\varepsilon\,\delta\calE_{\JMP}^{\tau,\tau+1}$.
Since the RK method $\phiRK$ is of order $p$, it can be written as
\begin{equation}\label{eq:oytp}
  \overline{\ys}^{\tau+1} = \phiRK(\ys^{\tau})=\ys(t^{\tau+1}) + e^{\tau+1},
\end{equation}
where $e^{\tau+1}:=C\Delta t^{p+1} + O(\Delta t^{p+2})$,
$C\in\mathbb{R}^M$, and $\ys(t^{\tau+1})$ is the exact solution, at
time $t^{\tau+1}$, of
\begin{align*}
  \begin{cases}
    \ys'(t) = \FNs(\ys(t);t),
    & \quad\textrm{in}\;(t^{\tau},t^{\tau+1}],\\
      \ys(t^{\tau})=\ys^{\tau}.
    & 
  \end{cases}
\end{align*}
The exact solution satisfies the conservation law
\eqref{eq:law-energy} so that
\begin{equation*}
  \calE_{\TOT}(\ys(t^{\tau+1})) - \calE_{\TOT}(\ys^{\tau})
  + \int_{t^{\tau}}^{t^{\tau+1}}\calE(\ys(t))\,dt = 0.
\end{equation*}
Exploiting this conservation law and equation \eqref{eq:oytp}, the
numerator in the expression of $\gamma^{\tau}$ in \eqref{eq:gamma:def}
can be written as
\begin{equation}\label{eq:gamma-num}
    \begin{aligned}
  \calE_{\TOT}( \phiRK(\ys^{\tau}) ) - \calE_{\TOT}(\ys^{\tau}) + \delta\calE^{\tau,\tau+1}
  = &\, \mu^\top e^{\tau+1} + \big<e^{\tau+1},e^{\tau+1}\big>
  + 2\big<\ys(t^{\tau+1}),e^{\tau+1}\big>\\
 &  + \Delta t\sum_{i=1}^{\NRK} b_i\calE(\calY_i)
 -\int_{t^{\tau}}^{t^{\tau+1}}\calE(\ys(t))\,dt\\
 = &\,\Delta t^{p+1}\big(\mu^\top C+\big<\ys^{\tau},C\big>\big) + e_q^{\tau+1} + O(\Delta t^{p+2}),
    \end{aligned}
\end{equation}
where $e_q^{\tau+1}$ is the quadrature error introduced by the
approximation $\delta\calE^{\tau,\tau+1}$ in \eqref{eq:deltaE}.
By means of Taylor expansions, it can be shown that the quadrature
error has the order of the corresponding RK method, namely, for some
$c=c(\calE)\in\mathbb{R}$,
\begin{equation*}
  \left|\int_{t^{\tau}}^{t^{\tau+1}}\calE(\ys(t))\,dt -
  \Delta t\sum_{i=1}^{\NRK} b_i \calE(\calY_i)\right|\approx c \Delta t^{p+1} + O(\Delta t^{p+2}).   
\end{equation*}
Moreover, the denominator of $\gamma^{\tau}$ in \eqref{eq:gamma:def}
can be expanded using Taylor expansions and gives
\begin{equation}\label{eq:gamma-den}
  \begin{aligned}
    \calE_{\ELE}(\phiRK(\ys^{\tau})-\ys^{\tau}) & =
    \big<\ys(t^{\tau+1})-\ys^{\tau} + e^{\tau+1},
    \ys(t^{\tau+1})-\ys^{\tau} + e^{\tau+1}\big>\\
    &= \Delta t^2\big<\dfrac{d\ys}{dt}(y^{\tau}),\dfrac{d\ys}{dt}(y^{\tau})\big> + O(\Delta t^3)
    = \Delta t^2\,\calE_{\ELE}(\ys^{\tau}) + O(\Delta t^3).
  \end{aligned}
\end{equation}
Combining \eqref{eq:gamma-num} with \eqref{eq:gamma-den}, it follows
that $\gamma^{\tau}=1+\Delta t^{p-1}$ for $\Delta t\rightarrow
0$.

The conclusion follows from the fact that
\begin{equation*}
\begin{aligned}
 \ys^{\tau+1} & = \phiRK(\ys^{\tau}) + (\gamma^{\tau}-1)(\phiRK(\ys^{\tau})-\ys^{\tau})\\
 & = \ys(t^{\tau+1}) + O(\Delta t^{p+1})
 + (\gamma^{\tau}-1)(\ys(t^{\tau+1}) -\ys^{\tau}+ O(\Delta t^{p+1})).
 \end{aligned}
\end{equation*}
Taylor expansion of $y(t^{\tau})$ around $t^{\tau+1}$ yields
\begin{equation}\label{eq:ymRK}
 \ys^{\tau+1}=
 \ys(t^{\tau+1}) + O(\Delta t^{p+1})
 + (\gamma^{\tau}-1)\Delta t\, \ys'(t^{\tau+1}) + O((\gamma^{\tau}-1)\Delta t^2)).
\end{equation}
We consider the two cases (i) and (ii) separately and use the fact that $\gamma^\tau=1+O(\Delta t^{p-1})$.
\begin{itemize}
    \item[(i)] Eq.~\eqref{eq:ymRK} yields
    $$\ys^{\tau+1}-\ys(t^{\tau+1}) = O(\Delta t^{p}).$$
    \item[(ii)] Taylor expansion of $ \ys(t^{\tau}+\gamma^{\tau}\Delta t)=
    \ys(t^{\tau+1}+(1-\gamma^{\tau})\Delta t)$
    around $t^{\tau+1}$, together with Eq.~\eqref{eq:ymRK}, yields
    $$\ys^{\tau+1}-\ys(t^{\tau}+\gamma^{\tau}\Delta t) = O(\Delta t^{p+1}).$$
\end{itemize}
\ENDPROOF

%% file: sec7_results.tex

\section{Numerical experiments}
\label{sec:numerical:results}

In this section, we investigate and assess the conservation properties
of the numerical methods discussed in the previous sections.
We test seven variants, which differ for the temporal integration algorithm
and the numerical flux in the discretization of Maxwell's equations.
We denote them with the abbreviations reported in
Table~\ref{tab:methods}.

\begin{table}[ht!]
  \centering
\begin{tabular}{|lcl|}
  \hline
  IMC   & ---   & Implicit Midpoint rule with Central flux \\ \hline
  IMU   & ---   & Implicit Midpoint rule with Upwind flux \\ \hline
  RKC   & ---   & Explicit Runge-Kutta with Central flux \\ \hline
  RKU   & ---   & Explicit Runge-Kutta with Upwind flux \\ \hline
  MRKC  & ---   & Modified explicit Runge-Kutta with Central flux \\ \hline
  MRKU0 & ---   & Modified explicit Runge-Kutta with Upwind flux with $\varepsilon = 0$ \\ \hline
  MRKU1 & ---   & Modified explicit Runge-Kutta with Upwind flux with $\varepsilon = 1$ \\ \hline
\end{tabular}
\caption{Abbreviations for the numerical methods used in the numerical
  tests.  The term `explicit Runge-Kutta' refers to the non-adaptive
  third order Runge-Kutta method of Bogacki and Shampine
  \cite{Bogacki-Shampine:1989}.}
\label{tab:methods}
\end{table}

\medskip
Some of the numerical methods tested in this section are implicit schemes;
therefore, they require the solution of large nonlinear systems of
equations at each time step.
A practical implementation requires efficient nonlinear solvers
coupled with preconditioned iterative linear solvers.
We leave the development of efficient implementations of implicit
temporal integrators to future work.
In all tests that use implicit time integration, nonlinear systems
were solved with an unpreconditioned Jacobian-Free-Newton-Krylov (JFNK)
method with relative tolerance of $10^{-8}$ and absolute tolerance of
$10^{-50}$.
The Krylov linear solver that is embeddded in the JFNK method is the
generalized minimum residual (GMRES) method with relative tolerance of
$10^{-5}$ and absolute tolerance of $10^{-50}$.

\medskip
The development of finer and finer scales in velocity space is a
common phenomenon in collisionless plasmas, known as filamentation.
In a spectral discretization, this implies that higher order spectral modes are
progressively excited until recurrence effects develop when
the velocity-space structures reach scales that are no longer resolved
by the truncated Hermite expansion \cite{joyce71,cheng76}.
An artificial collisional operator can then be employed to prevent
recurrence by damping high order modes.
This collisional operator must always be used in a convergence sense, without significantly affecting the collisionless physics of interest.
We adopt here the same collisional operator introduced in
Ref.~\cite{Delzanno:2015}, \rev{by adding to the operator $\As$ in \eqref{eq:semi-discrete:Vlasov} the term
$\mathcal{C}(\fsN)=\nu \nabla_v\cdot(\vv\fsN+\frac12\nabla_v \fsN)$, where $\nu>0$ is the collision rate, tested against any pair of functions $(\Psi,\varphi)\in\calHpN\times\calVN$. When $\Psi=\Psi^{n,m,p}(\xiv^s)$ and $\varphi=\varphi^{\Il}(\xv)$, as defined in Section~\ref{sec:variational:formulation}, then the collisional term reads
\begin{equation}\label{collop}
  \nu 
  \left[\kappa_x^{-1} n(n-1)(n-2)
    +\kappa_y^{-1} m(m-1)(m-2)
    +\kappa_z^{-1} p(p-1)(p-2)\right]C_{n,m,p}^{s,\Il}(t),
\end{equation}
where $\kappa_{\beta}=N_{v_{\beta}}(N_{v_{\beta}}-1)(N_{v_{\beta}}-2)$ for
all $\beta\in\{x,y,z\}$.}
Since the proposed collisional operator does not directly act on the
first three Hermite modes, it does not affect the conservation of
total mass, momentum and energy~\citep{Delzanno:2015}.

\input{results_WI.tex}

\input{results_X.tex}

\input{results_OT.tex}

\input{results_Mom}


%% file: results_WI.tex
    \subsection{Whistler instability}
\label{sec:whistler}
We start with the whistler instability test, where an electromagnetic
whistler wave grows from an unstable particle distribution function with
different temperatures along and perpendicular to the background magnetic field\rev{, see e.g. \cite{Gary:2005}.} \rev{In collisionless plasmas, the whistler instability grows due to a gyroresonance between a small part of the electron population and the wave. As such, this problem represents a rather sensitive test of the model's ability to capture kinetic effects.}
In this test, we consider two different time steps $\Delta t = 0.01$
and $0.005$, and DG polynomial spaces of degree $N_{DG} = 1$ (linear)
and $2$ (quadratic).
All other numerical parameters and the initial condition are fixed and
are outlined below.
The computational domain of the physical space is reduced to one
dimension ($N_y = N_z = 1$), with $\Ox=[0,2\pi]$ partitioned into a
uniform grid with $N_x=72$ elements. Periodic boundary conditions are
considered.
For the numerical approximation in velocity, we use $10$ Hermite basis
function in each direction (i.e. $\Nn=\Nm=\Np=9$).
Further, for the initial configuration we consider the magnetic field
$B_x=1$, while the distribution functions for electrons and protons
(ions) are Maxwellian distributions with uniform density.
Moreover, we set $\alpha_\beta^s=\sqrt{2}v_{T_\beta^s}$, and
$u_\beta^s=0$ for $\beta\in\{x,y,z\}$ and $s\in\{e,i\}$.
We use a realistic ion-to-electron mass ratio $m_i/m_e = 1836$,
thermal velocities $v_{T_\beta^e}=v_{T_\beta^i}\sqrt{m_i/m_e}=0.125$
for $\beta\in\{y,z\}$ with the exception of the reduced electron
thermal velocity along the $x$-axis, $v_{T_x^e}=v_{T_x^i}\sqrt{m_i/m_e}=0.056$, to create an
anisotropic distribution function which is the source of the
instability.
The electron plasma/gyrofrequency ratio is
$\omega_{pe}/\omega_{ce}=4$, and the collisional rate in
\eqref{collop} is $\nu=1$.
In order to seed the whistler instability, we initialize a small
electron current perturbation in the $x$-direction, $j^e_x (\xv) =
10^{-3} \cos (x)$, by exciting the Hermite mode $C^e_{1,0,0}$
in~\eqref{eq:fsN:def}. 

\par
\begin{figure}[H]
  \centering
  \includegraphics[width=0.7\linewidth]{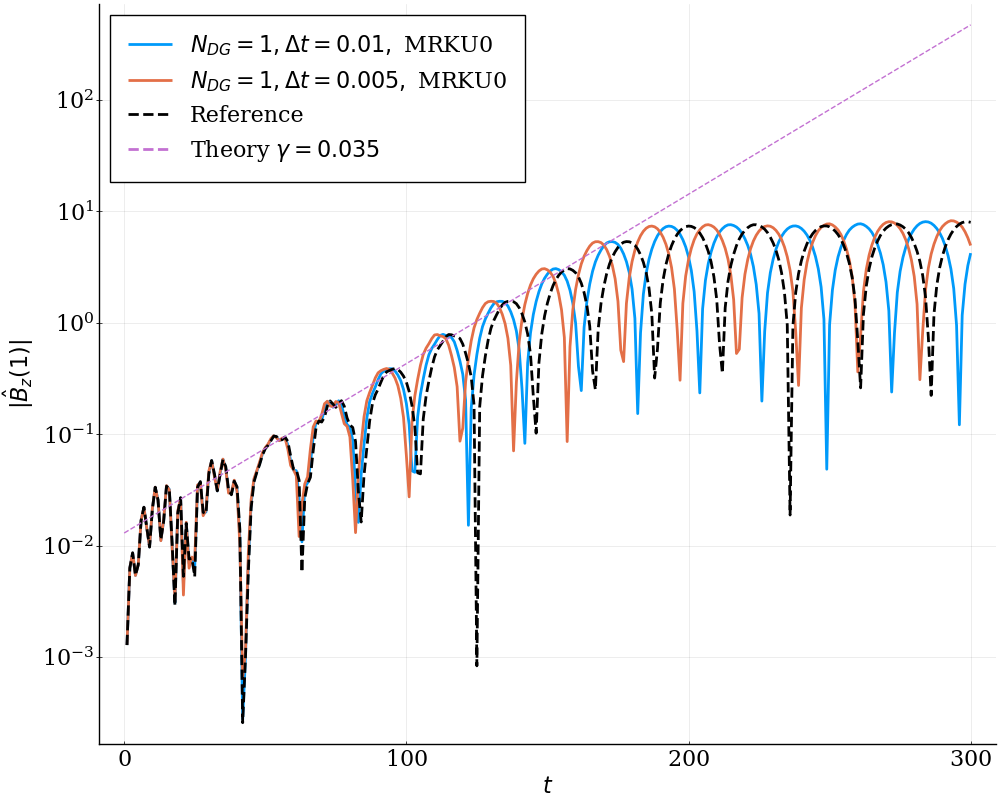}
  \caption{Whistler instability benchmark: time evolution of
    $|\hat{B}_z(1)|$.
    The numerical results show spatial accuracy deterioration for the
    MRKU0 method with $N_{DG}=1$, unlike all other tested methods
    (reference).}
  \label{fig:whistler_x}
\end{figure}

Before proceeding to the numerical study of the conservation
properties, we assess the capability of all considered methods to
reproduce the correct physical results, in particular, the exponential
growth of the electromagnetic whistler wave with the theoretically predicted growth rate from an initial small
perturbation.
It is sufficient to monitor the time evolution of the first Fourier
mode of the magnetic field $\hat B_z (k)$, with $k=1$, and
\begin{equation}
  \hat{B}_z(\ks) = \sum_I B_z(\xs_c^I) e^{i\ks\xs_c^I},
  \label{Bzfft}
\end{equation}
where $\xs_c^I$ is the center of the $I$-th \DG{} cell.
Most of the considered methods produce a visually indistinguishable
evolution of $\hat B_z(1)$, 
matching the theoretically predicted growth rate $\gamma = 0.035$. This has been labeled as 'Reference' solution in Fig.~\ref{fig:whistler_x}
(the only tests which visually deviate from reference curve are discussed below).
However, the method MRKU0 with $N_{DG} = 1$ exhibits a noticeable
deviation from the reference solution, independent of the time step,
as shown in Figure~\ref{fig:whistler_x}.
Increasing the spatial accuracy of MRKU0 by means of higher order DG
polynomials, $N_{DG}=2$, makes it consistent with the other methods and
produces a curve that is visually indistinguishable from the reference curve in Fig.~\ref{fig:whistler_x}.
It is interesting to note that, on this example, MRKU0 (where energy conservation is imposed exactly) is actually not as accurate as MRKU1 (which respects the energy dissipation introduced by upwind fluxes in Maxwell's equations).
\rev{We will further comment on this behavior at the end of Section~\ref{sec:Xmode}.}

\begin{figure}[H]
  \centering
  \includegraphics[width=0.7\linewidth]{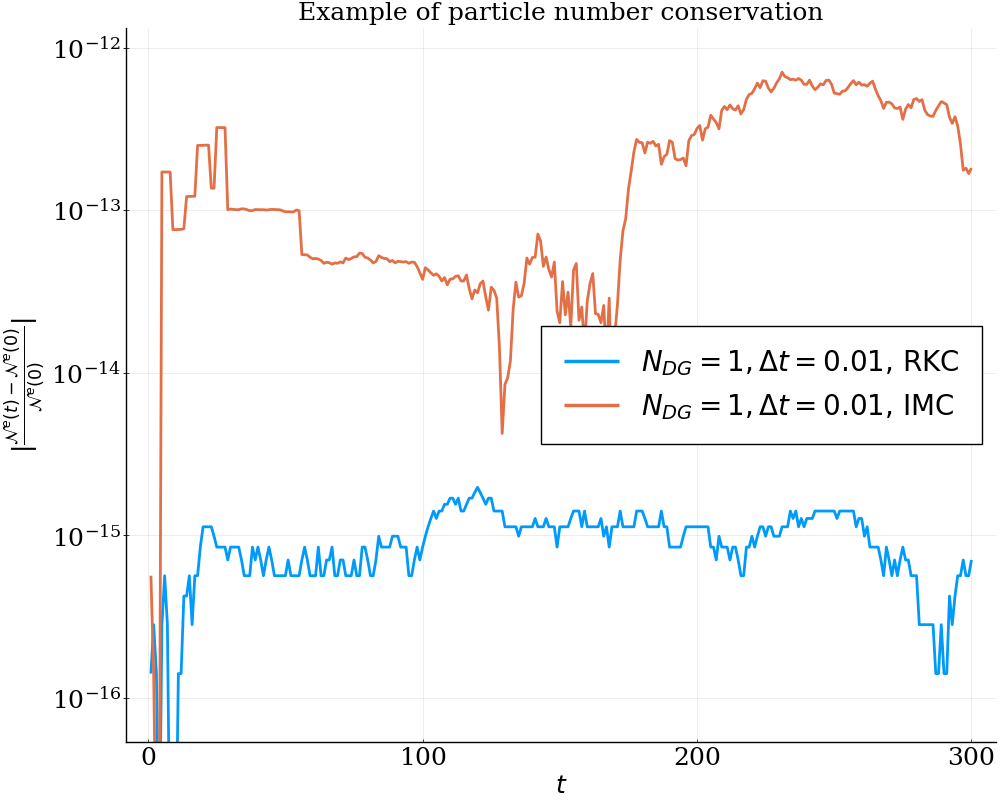}
  \caption{Whistler instability benchmark: time evolution of the error in the number of electrons for RKC and IMC with
    $N_{DG}=1$ and $\Delta t = 0.01$.}
  \label{fig:whistler_x:mass}
\end{figure}

\begin{figure}[H]
  \centering
  \includegraphics[width=0.31\linewidth]{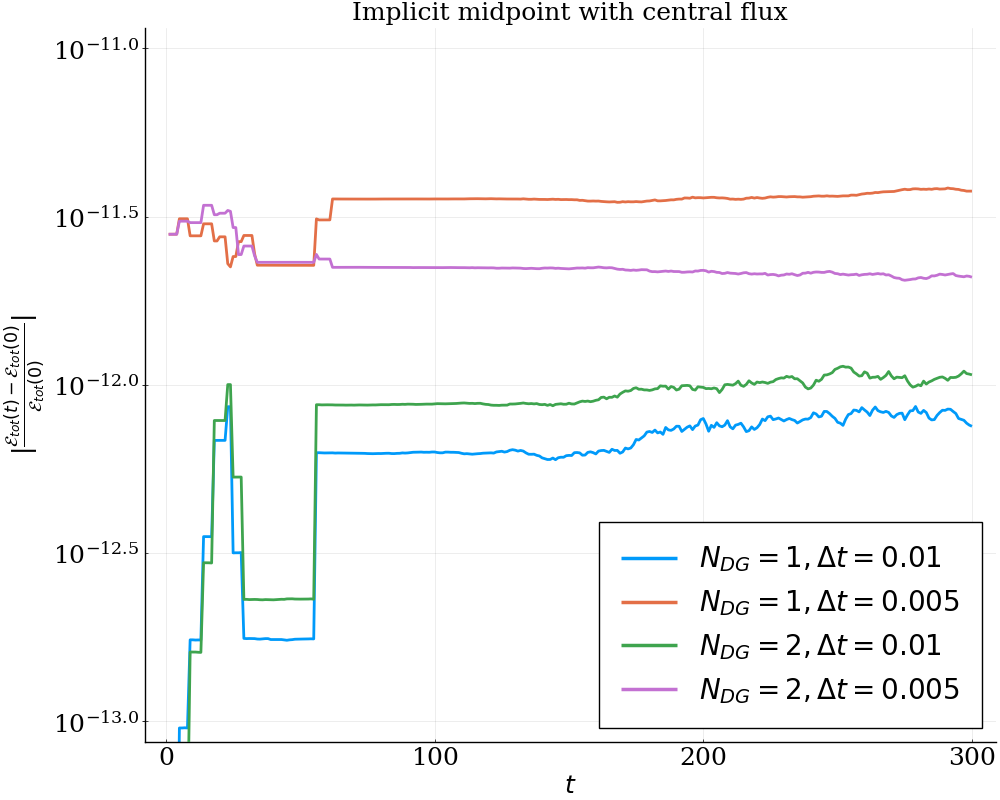}
  \includegraphics[width=0.31\linewidth]{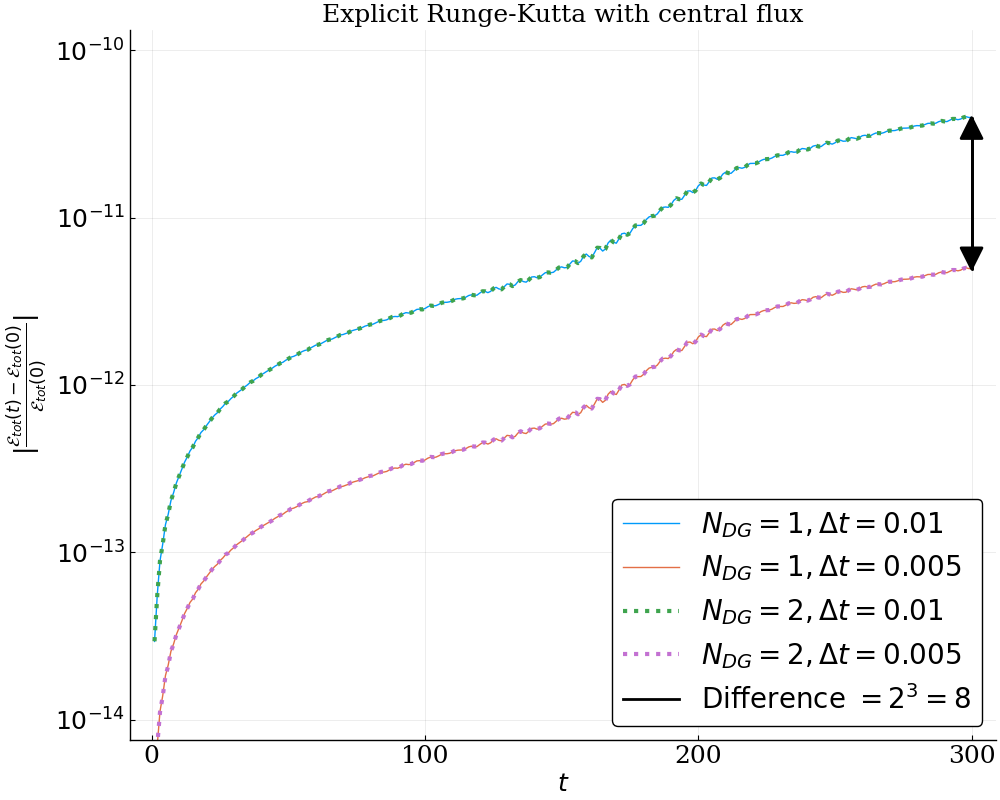}
  \includegraphics[width=0.31\linewidth]{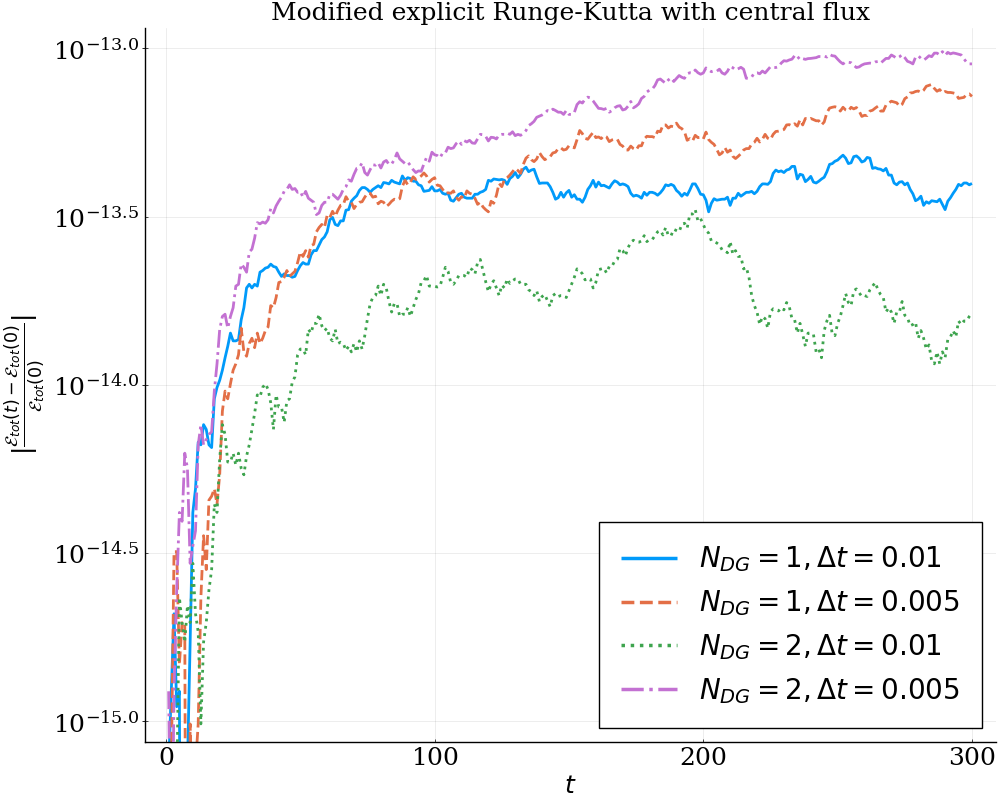}
  \caption{Whistler instability benchmark: time evolution of the relative error 
  in the conservation of the total energy for IMC (left), RKC (middle), and MRKC (right) methods, and
  for different spatial and temporal resolutions.}
  \label{fig:whistler_x:central:energy}
\end{figure}

All tested methods conserve the number of particles: Explicit methods
(RKC, RKU, MRKC, MRKU0, MRKU1) up to machine precision, while implicit
methods (IMC, IMU) with error bound by the tolerance of the nonlinear solver
and independent of the space-time resolution (i.e., $\Delta t$, $N_{DG}$).
In Fig.~\ref{fig:whistler_x:mass} we report the evolution of the
relative error in the number of electrons with respect to the initial
value for two methods, RKC and IMC with $N_{DG}=1$ and $\Delta t =
0.01$.
In the numerical study of the energy conservation properties of the
proposed methods, we first consider the family of methods where
central fluxes are used in the spatial approximation of Maxwell's
equations, i.e., IMC, RKC, and MRKC.
Note that, under this assumption, the energy flux at the mesh
interfaces $\calE_{\JMP}$ in \eqref{eq:calE:jmp:def:cont} vanishes.
In Fig.~\ref{fig:whistler_x:central:energy} we report the evolution
of the error in the total energy, defined in
Eq. \eqref{eq:energyN:def}, for different temporal integrators.
In the left panel, the IMC method is studied for different temporal
and spatial resolutions: as expected, all runs show bounded errors with
magnitude of the order of the nonlinear solver tolerance.
Next, the energy error of the RKC method (middle panel) depends on the
temporal approximation, and is completely independent of the spatial
resolution (curves with different $N_{DG}$, but with identical time
step, coincide).
Moreover, curves with different time step ($\Delta t = 0.01$ and
$0.005$) differ by a factor of $\sim 2^3 = 8$ 
(marked as a black two sided arrow), which is consistent with the third order of the Runge-Kutta method in use.
The right panel in Fig.~\ref{fig:whistler_x:central:energy} shows
the error in the conservation of total energy associated with the MRKC
scheme.
In this case, the error in energy conservation is maintained around
machine precision\footnote{Here and throughout this section, ``machine precision'' is used to describe anything from $10^{-15}$ to $10^{-13}$.}, with a small dependence on temporal and
spatial resolutions.

\begin{figure}[H]
  \centering
  \includegraphics[width=0.31\linewidth]{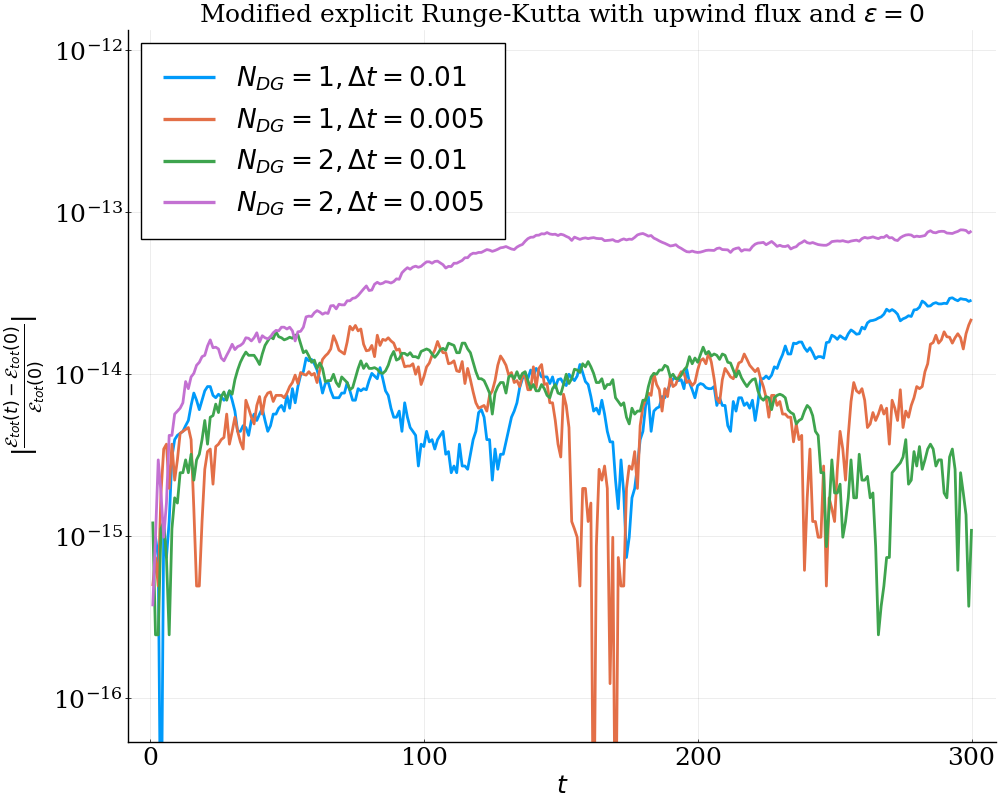}
  \includegraphics[width=0.31\linewidth]{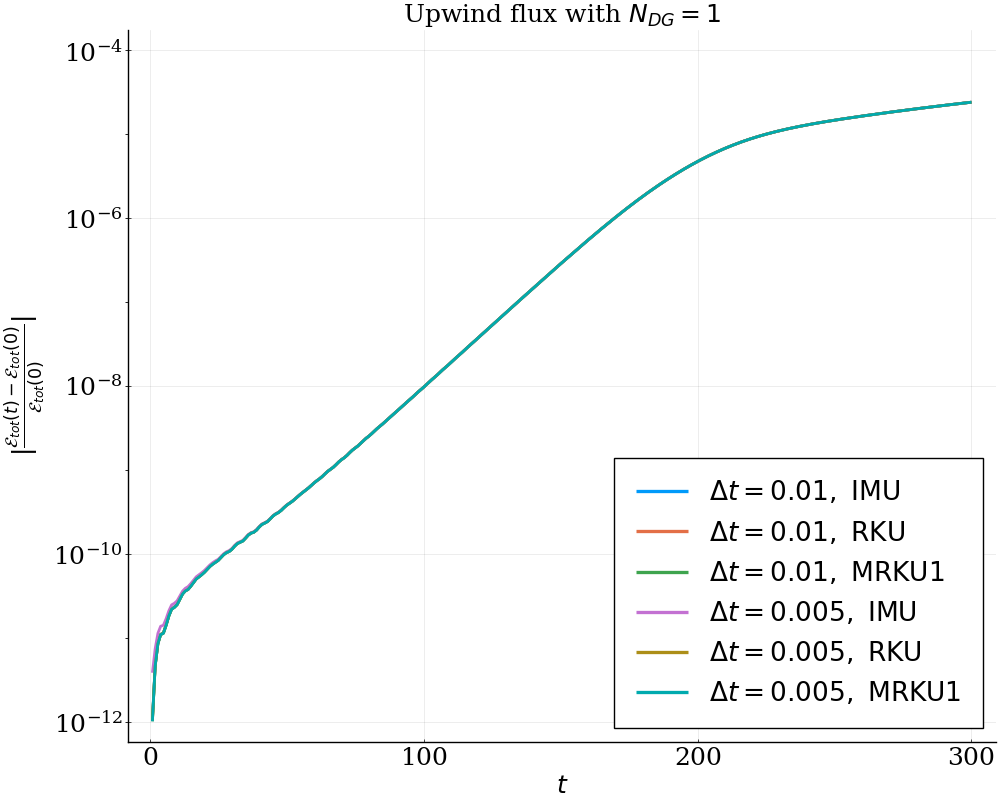}
  \includegraphics[width=0.31\linewidth]{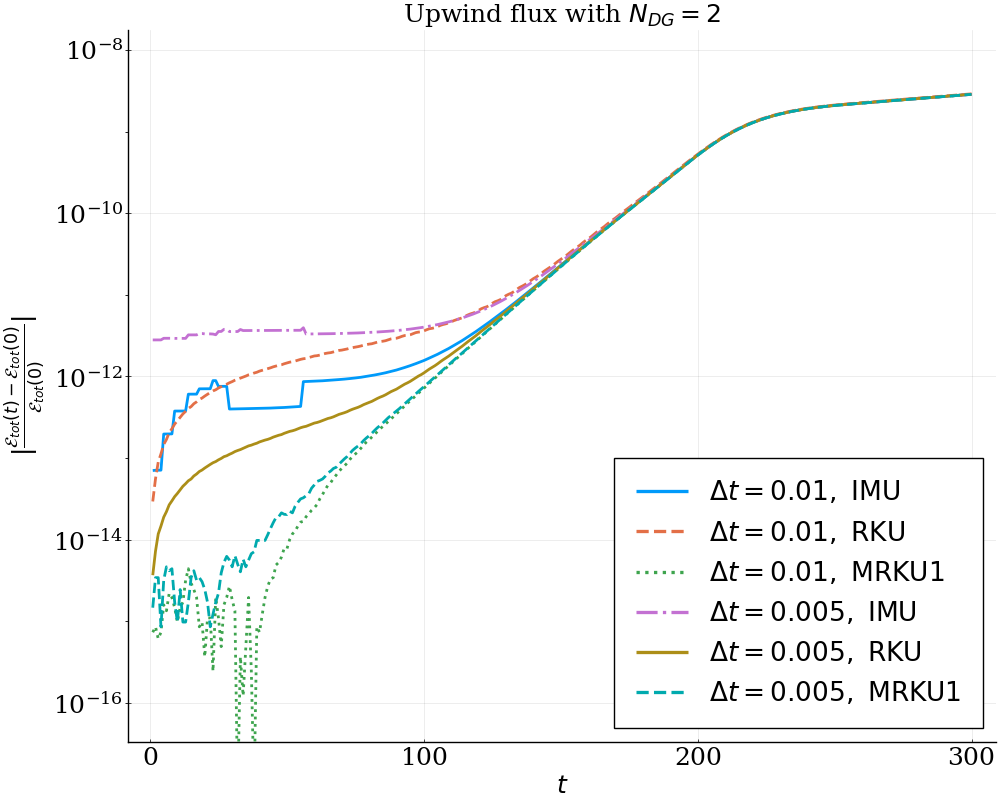}
  \caption{Whistler instability benchmark: time evolution of the relative error 
  in the conservation of the total energy for
  MRKU0 (left), 
      IMU, RKU, MRKU1 with $N_{DG}=1$ (middle),
  and IMU, RKU, MRKU1 with $N_{DG}=2$ (right),
  for different time steps.}
  \label{fig:whistler_x:upwind:energy}
\end{figure}

Finally, we consider the error in the conservation of the total energy
when upwind fluxes are employed in the spatial approximation of
Maxwell's equations.
The numerical results are presented in
Fig.~\ref{fig:whistler_x:upwind:energy} and corroborate the findings of
Theorem~\ref{theorem:modified:RK:discrete:energy:conservation}.
The MRKU0 method ensures conservation of the total energy to machine
precision and this is fairly independent of the space-time resolution, as
shown in the left panel of Fig.~\ref{fig:whistler_x:upwind:energy}.
All the other upwind methods, i.e., IMU, RKU, MRKU1, dissipate energy: for
completeness of exposition, we report, in the middle and right panels
of Fig.~\ref{fig:whistler_x:upwind:energy}, the absolute value of
the error in the conservation of the total energy.
The spatial polynomial degree is $N_{DG}=1$ for the results shown in
the middle panel, while $N_{DG}=2$ for the results in the right panel.
Note that the quantity shown in the middle and right panels of
Fig.~\ref{fig:whistler_x:upwind:energy} is related to the discrete energy flux
at the mesh interfaces, namely, the right hand side of
Eq.~\eqref{eq:law-energy}. 
The value of this quantity depends on both the temporal discretization
and the spatial approximation.
For $N_{DG}=1$ and sufficiently small time step, the error is
independent of the temporal integrator of choice and is dominated by the spatial discretization error (middle panel).
For $N_{DG}=2$ (right panel), the energy error has a very small magnitude (less than $\sim 10^{-12}$) until around $t=100$.
Before this time, we record different behaviors for different temporal
integrators. 
The error for MRKU1 is around machine precision with a rather negligible dependence on $\Delta t$. The error for IMU is generally higher and is determined by the tolerance of the nonlinear solver (as determined by different runs with the same parameters but a lower tolerance of the nonlinear solver, not shown). Finally, we can see that the traditional Runge-Kutta scheme RKU
performs rather well in this test case, with the error on the total energy remaining comparable to that for IMU. RKU also shows the expected third order scaling of the algorithm with respect to $\Delta t$.
For times $t>150$, all curves converge to the same behavior, signaling that the spatial discretization error is now dominant. This is consistent with Eq. (\ref{eq:energyN:def}), which shows that the variation of the total energy accumulates in time by summing the negative-definite contributions coming from the spatial discretization. 


\revtwo{We note that nonlinear evolution of the whistler instability may in principle result in development of strong non-Maxwellian features in the distribution function. Depending on the spectrum and the saturation amplitude of the excited modes under a given set of conditions, such departures may pose challenges to a method that uses a small number of Hermite modes and applies strong damping to high-order moments.
The growth rate and the saturation amplitude of the solution presented in this section were checked for convergence with respect to the number of Hermite polynomials by performing additional simulations with $20$ Hermite basis functions in each direction and varying collisionality. The solution was found to be fully converged with respect to the number of Hermite polynomials, while reducing the collisionality by almost an order of magnitude resulted in saturation amplitude changing by 7.5\% and a phase shift of the oscillations in saturated state, which is inconsequential in most practical applications.}

%% file: results_X.tex
\subsection{High frequency electromagnetic waves}
\label{sec:Xmode}
From the results presented above, it may be tempting to conclude that using the central fluxes in the approximation of Maxwell's equations is preferable in practice, given their energy conservation properties and the fact that enforcing energy conservation in upwind-based methods leads to a degradation of accuracy. However, methods based on upwind fluxes may be preferable in some situations, in particular because they tend to reproduce the physical dispersion properties of plasma waves in a better way compared to the central-flux counterparts \cite{hesthaven2004high,sarmany2007dispersion}. 
\par 
In this section, we illustrate the difference between central and upwind fluxes in resolving the dispersion properties of X mode waves in the cold-plasma limit.
In the simplest case, an X mode is an electromagnetic
electron wave (the ions remain essentially static)
with wave vector and wave electric field perpendicular to the background magnetic field.
\par
In this test, we use a one dimensional ($N_y = N_z = 1$)
spatial domain $\Ox = [0,100\pi]$ with $N_x = 500$, periodic boundary conditions,
and linear DG polynomials $N_{DG}=1$.
We consider the temporal interval $[0,100\pi]$ with
time step $\Delta t = \pi/25$.
Since treating the plasma as a fluid is sufficient for this test, we use only $4$ Hermite polynomials per velocity direction.
Further, the background magnetic field is initialized to $B_z=1$, and the distribution function for the electrons is Maxwellian.
Ions constitute a static background.
Moreover, we set $\alpha_\beta^e=0.002$, and $u_\beta^e=0$ for $\beta\in\{x,y,z\}$.
The electron plasma/gyrofrequency ratio is $\omega_{pe}/\omega_{ce}=1$
and no artificial collisionality is used ($\nu=0$).
We perturb the initial electron current so that
\begin{align}
    j_y^e (\xv) = 10^{-4} \exp \left ( - \left ( \frac{x - L_x/2}{0.15 L_x} \right )^2 \right ) \sum_{i=1}^4 A_i \sin \left ( k_i x + \Phi_i \right ),
\end{align}
with $k_i=50i k_{\min}$, for $i\in\{1,2,3,4\}$ and
$k_{\min} = 2\pi/L_x$,
$A_i   \in  \{ 0.598, 0.517, 0.193, 0.218\}$,
$\Phi_i \in \{0.305, 0.586, 0.050, 0.089\}$.

\begin{figure}[H]
  \centering
  \includegraphics[width=0.45\linewidth]{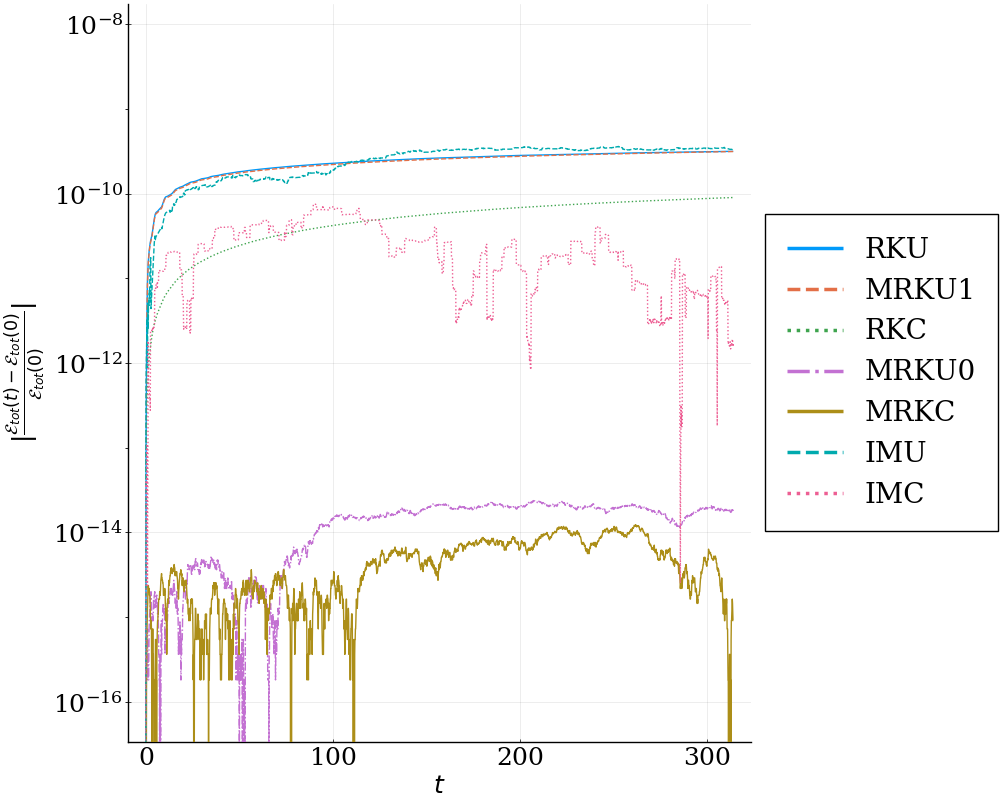}
  \includegraphics[width=0.45\linewidth]{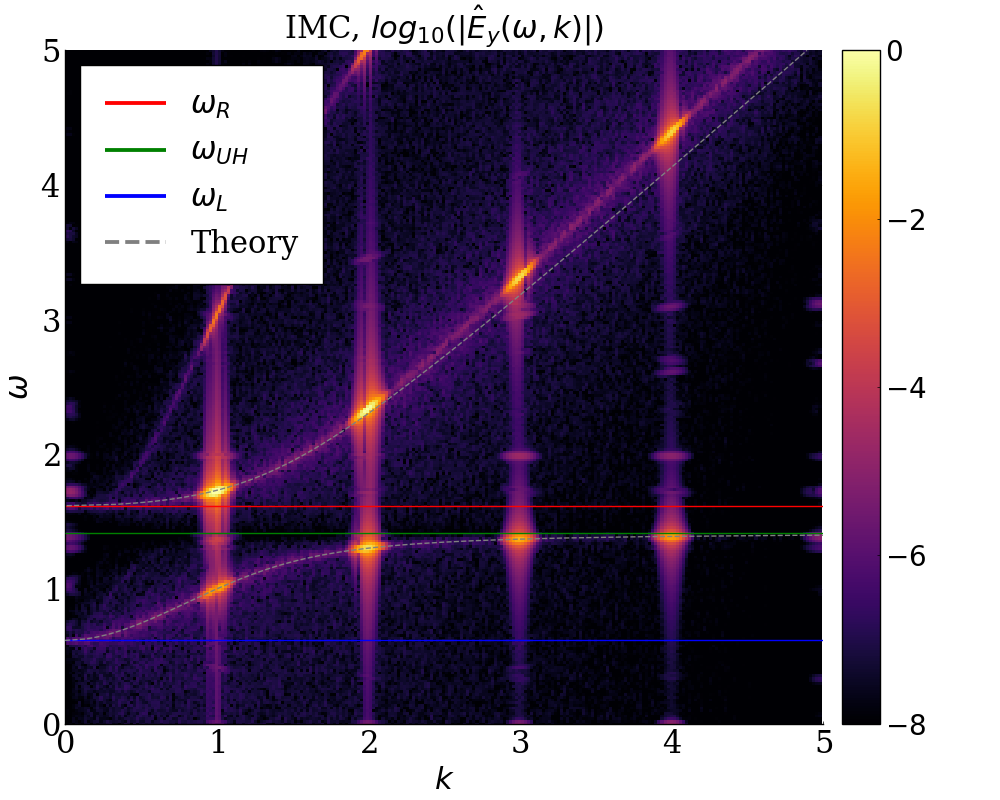} \\
  \includegraphics[width=0.45\linewidth]{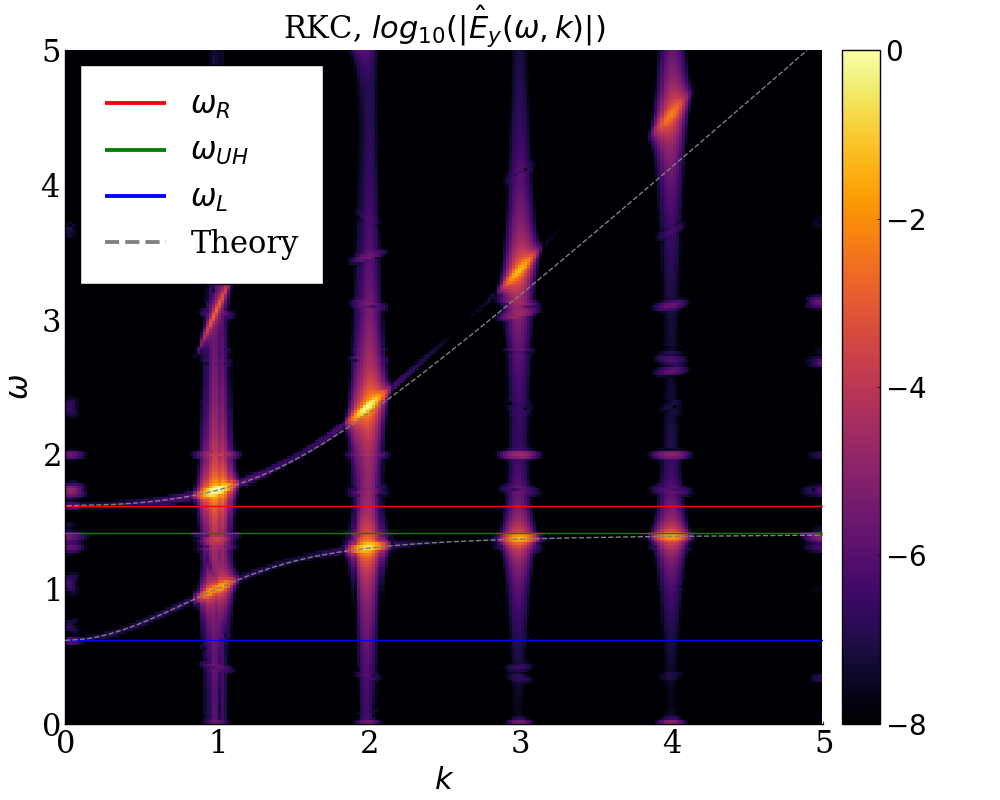}
  \includegraphics[width=0.45\linewidth]{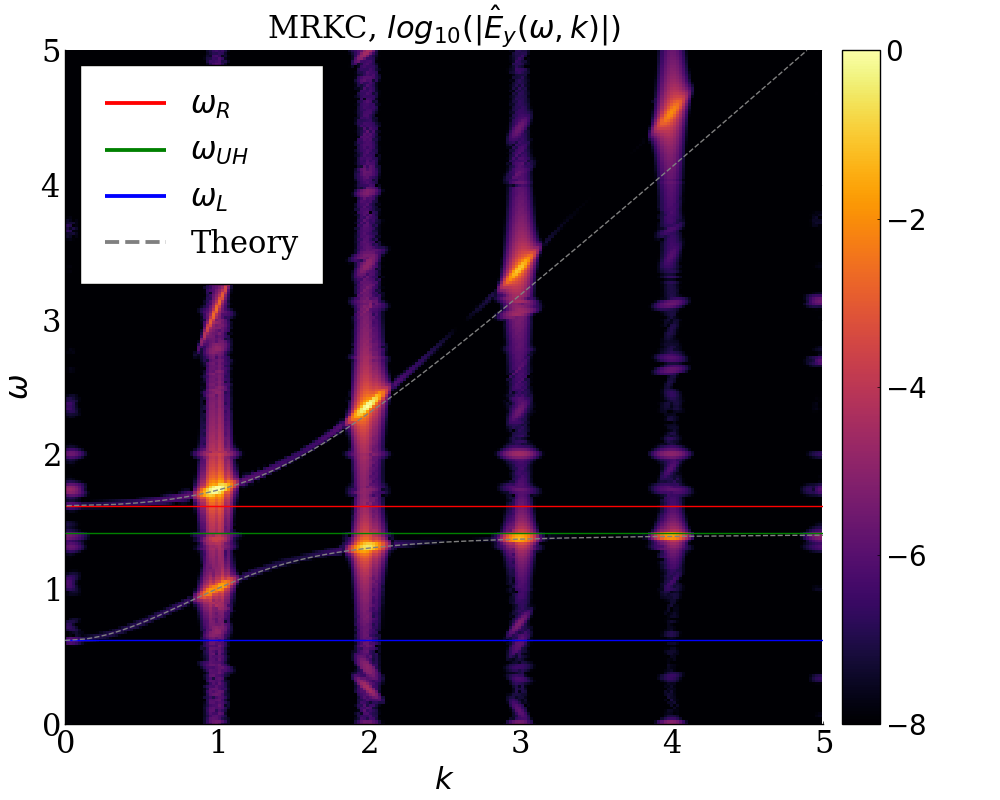} 
  \caption{X mode test: relative error in the conservation of the total energy (top left);
           wave spectrum computed from $E_y(\xv,t)$ 
           for methods with central fluxes, i.e.,
           IMC (top right), RKC (bottom left), and MRKC (bottom right). }
  \label{fig:xmode:c}
\end{figure}

The top left panel of Fig.~\ref{fig:xmode:c} shows the relative error in the conservation of the total energy
(the normalization factor $\mathcal{E}_{tot} (0)$ includes the energy of the background magnetic field)
for all methods summarized in Table~\ref{tab:methods}.
Similar considerations as for the whistler instability can be made,
with the exception that this test does not have any instability, so the excited waves
remain of small amplitude. The relative energy error is thus small for all methods. 
In summary, the
MRKC and MRKU0 methods conserve energy up to machine precision,
IMC conserves energy up to a level determined by the tolerance of the nonlinear solver,
the energy error in RKC is quite small and depends on the accuracy of the temporal discretization 
(as demonstrated by additional runs with different values of $\Delta t$ but all other parameters unchanged, not shown)
and IMU, RKU, and MRKU1 dissipate energy as expected.
Note that the energy errors for RKU and MRKU1 are very similar, suggesting that the RK correction in the upwind-based scheme is not so important for this example.

\begin{figure}[H]
  \centering
  \includegraphics[width=0.45\linewidth]{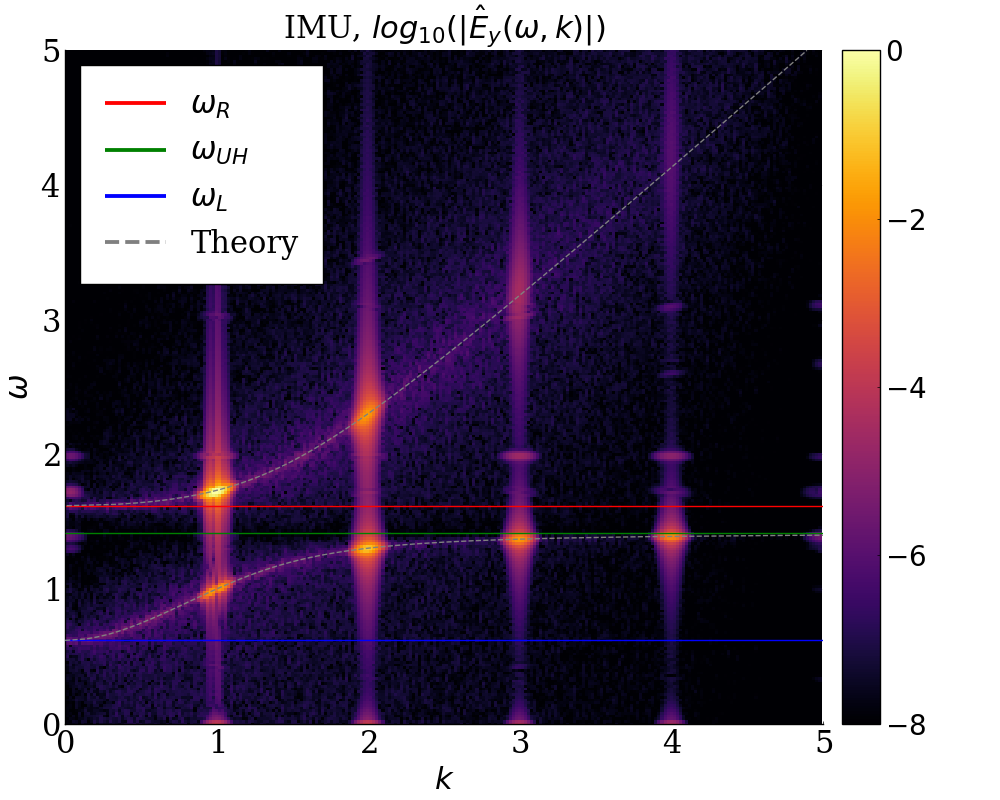}
  \includegraphics[width=0.45\linewidth]{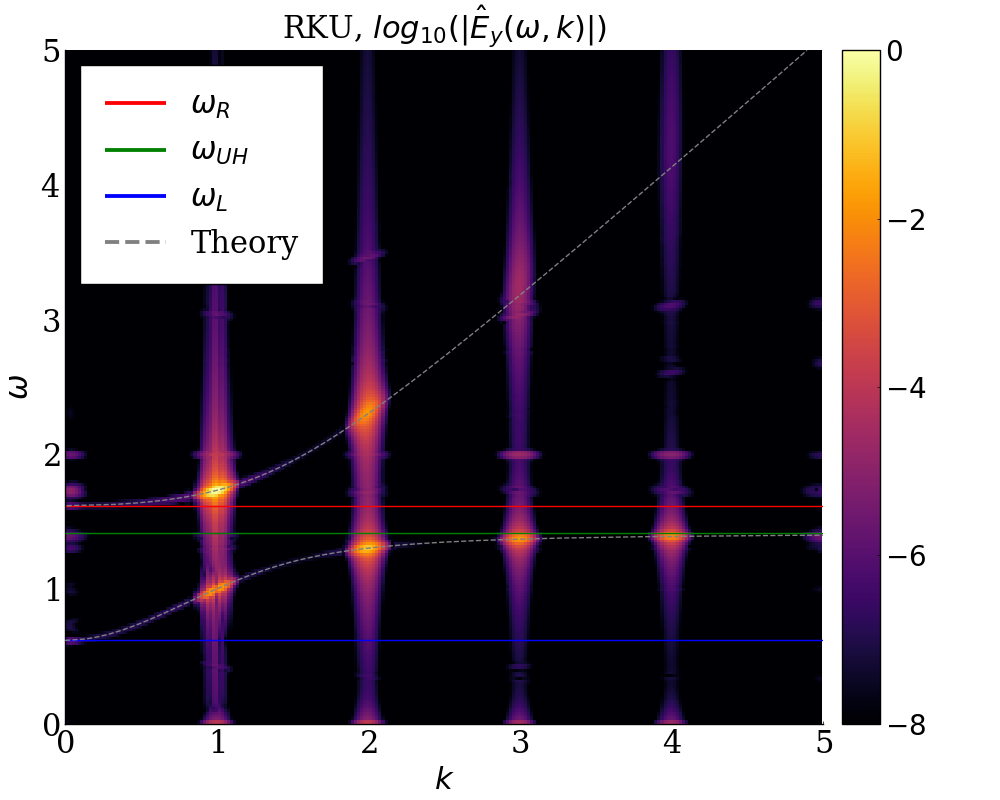} \\
  \includegraphics[width=0.45\linewidth]{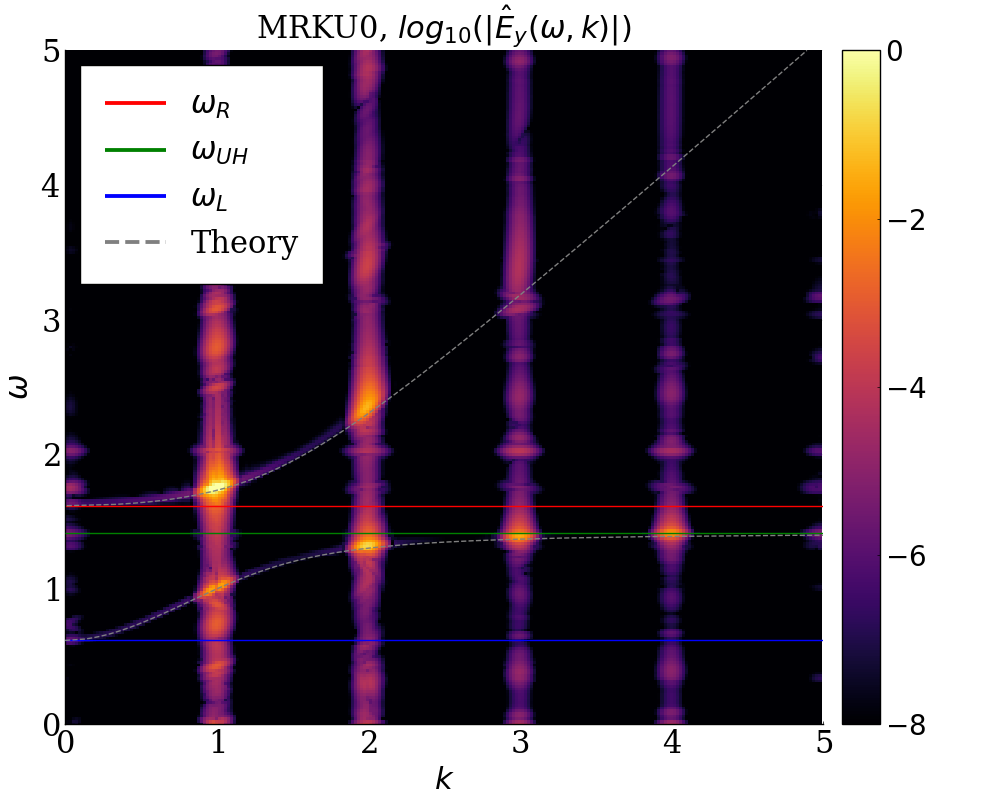}
  \includegraphics[width=0.45\linewidth]{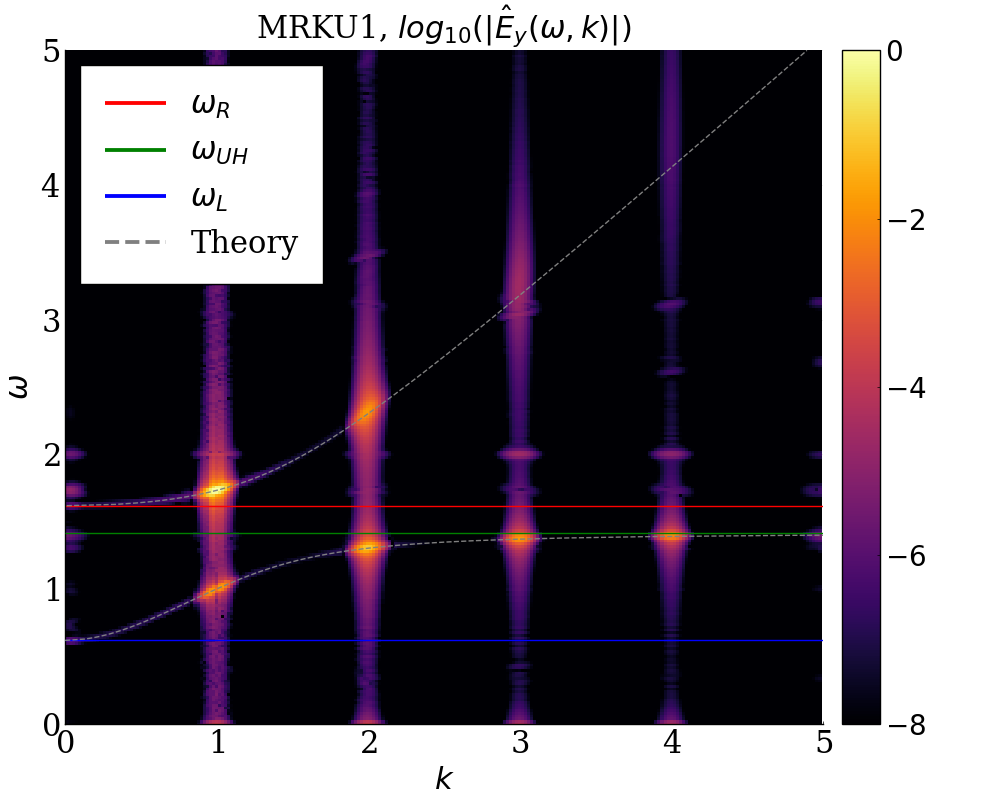} 
  \caption{X mode test: wave spectrum computed from $E_y(\xv,t)$ 
            for methods with upwind fluxes, i.e.,
            IMU (top left), 
            RKU (top right), 
            MRKU0 (bottom left), 
        and MRKU1 (bottom right).}
  \label{fig:xmode:u}
\end{figure}

In order to diagnose the dispersion properties of the X mode waves, we plot the spectrum of the
electric field $E_y$ computed with the discrete Fourier transform, namely
\begin{align}
    \hat E_y (\omega, k) = \sum_{n,I} E_y (x_c^I,t_n) e^{i\omega t_n} e^{ikx_c^I} H_n, \label{Eyhat}
\end{align}
where $t_n = n \Delta t$, $x_c^I$ is the center of $I$-th DG cell,
and $H_n = \sin^2 (\pi n/N_t)$ is a Hann time window applied to the data to reduce potential
artifacts of the Fourier transform originating from non-periodic time signals, 
$N_t$ is a number of time steps.
The logarithm of the absolute value of $\hat{E}_y$ in \eqref{Eyhat}
is shown in Fig.~\ref{fig:xmode:c} and in Fig.~\ref{fig:xmode:u},
for the numerical methods using central and upwind fluxes, respectively.
The theoretical dispersion relation for the X mode (lower and upper branches) 
is plotted with a gray dashed line, 
the resonance frequency $\omega_{UH}$ (the upper hybrid frequency) with a green dashed line,
and the cut-off frequencies $\omega_R$ and $\omega_L$ with red and blue dashed lines, respectively.
All methods resolve the lower frequency branch ($\omega_L < \omega < \omega_{UH}$) satisfactorily.
By contrast, differences can be observed in the resolution of the high frequency branch ($\omega > \omega_R$).
First, high frequency waves which are spatially well-resolved,
i.e., with $k<2.5$ (more than 4 DG cells per wavelength)
are well captured in both upwind- and central-based methods.
However, poorly resolved modes $k>2.5$ (less than 4 DG cells per wavelength)
are mostly suppressed in upwind-based methods,
while they are still present in central-based methods.
This may be attributed to dissipation in upwind methods suppressing the
unresolved modes and making the methods generally more stable.
We also note that the implicit methods (IMC and IMU) are more noisy for this example:
all wavenumbers are excited to amplitudes that are small but visible on the logarithmic-scale plot by the fact that the nonlinear system converges only up to the relatively high tolerance of the nonlinear solver considered here.
However, the IMU method suppresses most of this noise
by intrinsic dissipation, as one can see by comparing with the spectrum of IMC.
Additionally, the MRKU0 method exhibits
numerous unphysical modes in the spectrum, particularly for higher frequencies $\omega$. 

Overall, this example shows a case where upwind fluxes in Maxwell's equations might be preferable to central fluxes to avoid high-frequency modes that are not resolved by the mesh. It also confirms that forcing the exact conservation of the total energy in the explicit upwind scheme (MRKU0) might lead to a lack of accuracy of the overall simulation.
\rev{This behavior, also observed in the whistler instability test, can be interpreted as caused by an erroneous energy redistribution in space or, more precisely, spectrally.
We believe that the projection step, which is only aware of the temporal discretization, can properly correct energy errors which come from the Runge-Kutta step, but not due to the spatial discretization.
Methods based on upwind fluxes in Maxwell's equations damp poorly resolved modes with $k\sim \Delta x^{-1}$ featuring additional stability guarantees compared to methods based on central fluxes (cf. Figure~\ref{fig:xmode:c} vs. Figure~\ref{fig:xmode:u}).
This damping of the energy has magnitude that depends on the spatial resolution.
If one tries to compensate this damping by a projection during the Runge-Kutta time step (as in MRKU0), 
the excessive energy is deposited to other modes and this leads to accuracy loss.
On the other hand, accuracy can be retained in case of upwind fluxes when the upwind damping is retained, i.e., in MRKU1, or when the damping is absent, i.e., in methods based on central fluxes.}

%% file: results_OT.tex
\subsection{Orszag-Tang vortex}
\label{sec:OT}
In the last test, we investigate the performance of the proposed methods in conserving the total energy when solving the
Orszag-Tang vortex problem.
In the Orszag-Tang vortex problem \cite{Orszag:1979}, two large-scale
vortices are initialized and evolve by forming smaller and smaller
current sheets and other filamentary structures.
 \rev{In collisionless plasma, the  Orszag-Tang initial conditions lead to development of turbulence via breaking (reconnection) of the current sheet formed by the initial evolution~\citep[see e.g.][]{Parashar2009}. The overall energy dissipation in such a case is thought to be dominated by kinetic effects that become significant at small scales~\citep[see e.g.][and references therein]{Wan2012}. Therefore, this test is an example of the complex, multi-scale problem involving transition between the large-scale, fluid-like behavior of the plasma and the small-scale, dissipative processes involving kinetic physics.}

\rev{We consider} temporal and spatial discretizations with different resolutions:
in particular, the time steps $\Delta t \in \{ 0.05,0.025\}$,
and DG polynomial spaces $N_{DG} = 1$ (linear) and $N_{DG}=2$ (quadratic).
All other numerical parameters and initial conditions  
are fixed and are outlined below.
The computational domain in physical space is two dimensional ($N_z= 1$)
with $\Ox=[0,L_x]\times[0,L_y]$, $L_x=L_y=50$, and periodic boundary conditions.
It is discretized with a uniform grid with $N_x = N_y = 108$ elements.
We set $\Nn=\Nm=\Np=3$ for the velocity Hermite expansion.
The initial magnetic field is set to
\begin{align*}
  B_x(\xv) &=-\delta B \sin(k_y y + 4.1),\\
  B_y(\xv) &= \delta B \sin(2k_x x + 2.3),\\
  B_z(\xv) &= 1,
\end{align*}
with $\delta B = 0.2$, $k_x = 2 \pi /L_x$, $k_y = 2 \pi/L_y$.
The values $4.1$ and $2.3$ are arbitrary phases that remove any
artificial symmetry in the initial setup.
The distribution functions for electrons and ions are initialized to be
shifted Maxwellian distributions with spatially uniform density
(species superscripts are omitted for clarity),
\begin{equation*}
  f(\xv,\vv,t) = 
  \prod_{\beta\in\{x,y,z\}} \frac{1}{v_{T_{\beta}} \sqrt{2\pi}} 
  \exp \left [ -\frac{ (v_\beta - U_\beta (\xv)) ^2}{2 v_{T_\beta}^2} \right ],
\end{equation*}
with electron and ion velocities
\begin{align*}
  U^e_x(\xv) &= U^i_x(\xv) = -\delta B v_a \sin(k_y y + 0.5),\\
  U^e_y(\xv) &= U^i_y(\xv) = \delta B v_a \sin(k_x x + 1.4),\\
  U^e_z(\xv) &= -\frac{ \delta B \omega_{ce}}{\omega_{pe}} \left (2 k_x \cos(2k_x x + 2.3) + k_y \cos (k_y y + 4.1 )\right ) ,\\
  U^i_z(\xv) &= 0,
\end{align*}
where $v_a = 0.1$ and $\omega_{pe}/\omega_{ce} = 2$.
The values $0.5$ and $1.4$ are random phases and $U_z^e$ was chosen to
satisfy Amp\`ere's law at time $t=0$.
Other parameters include $m_i/m_e = 25$ and the artificial collision rate $\nu=1$.
The species-dependent parameters used in the definition of the Hermite
functions are set as follows:
$\alpha_\beta^e = \sqrt{2} v_{T_\beta}^e = 0.25 $, $u_\beta^e=0$ and
$\alpha_\beta^i = \sqrt{2} v_{T_\beta}^e/\sqrt{m_i/m_e} = 0.05 $,
$u_\beta^i=0$ with $\beta\in\{ x,y,z\}$.

\begin{figure}[H]
  \centering
   \includegraphics[width=0.24\linewidth]{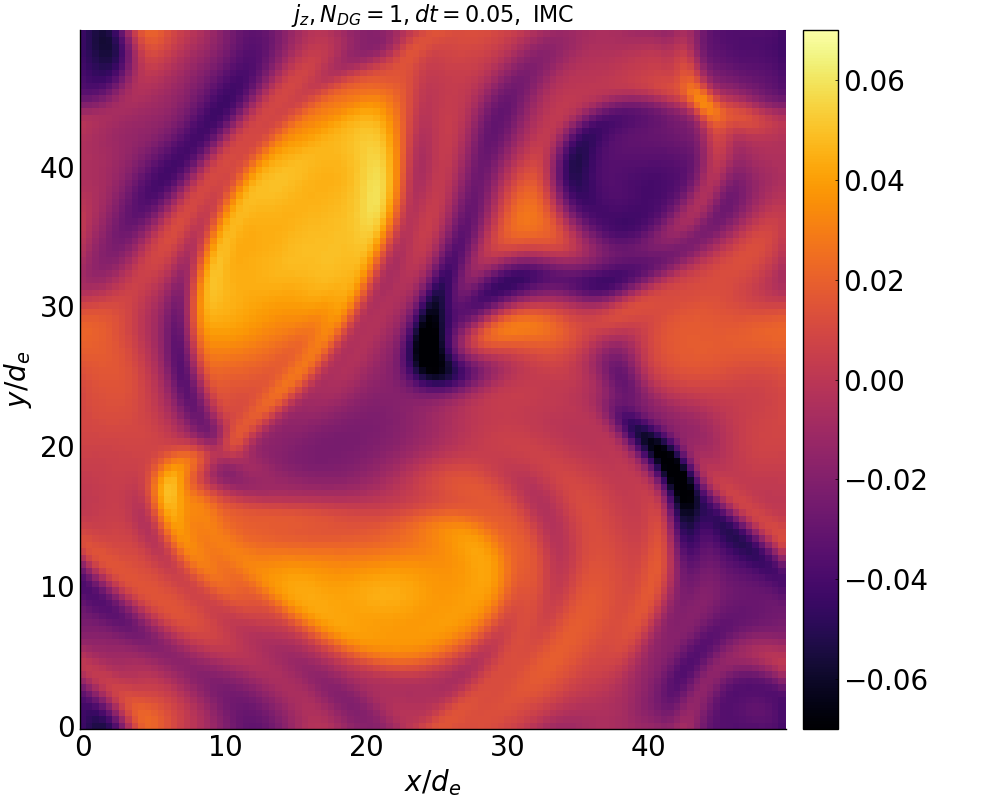}
  \includegraphics[width=0.24\linewidth]{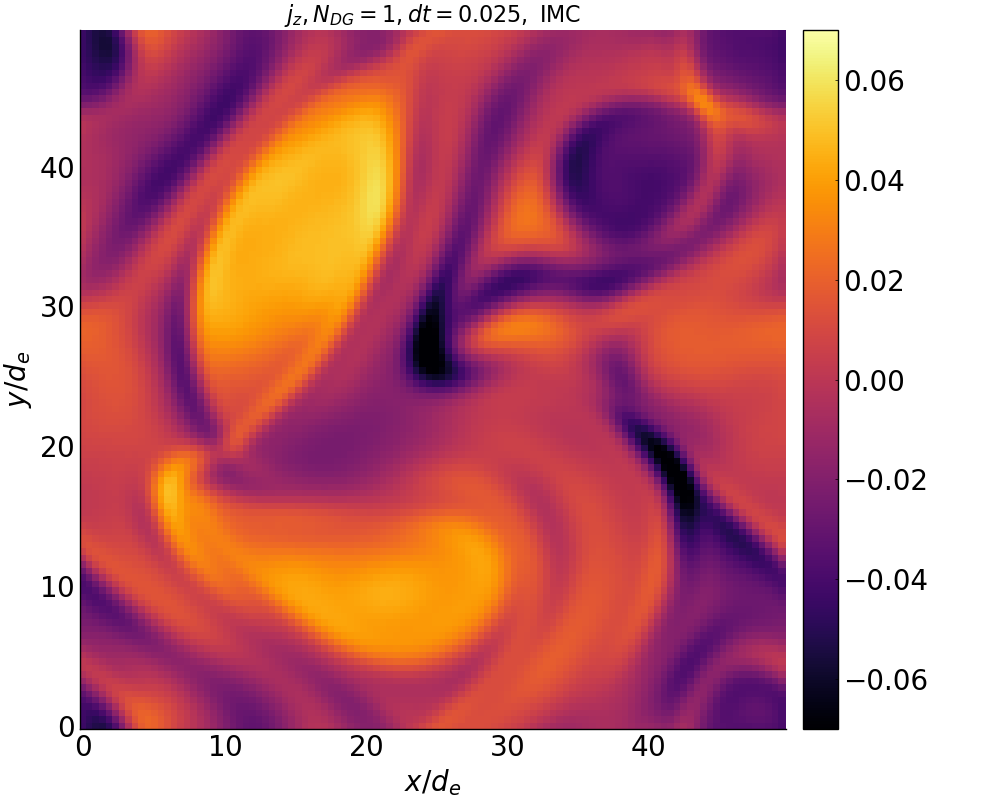}
   \includegraphics[width=0.24\linewidth]{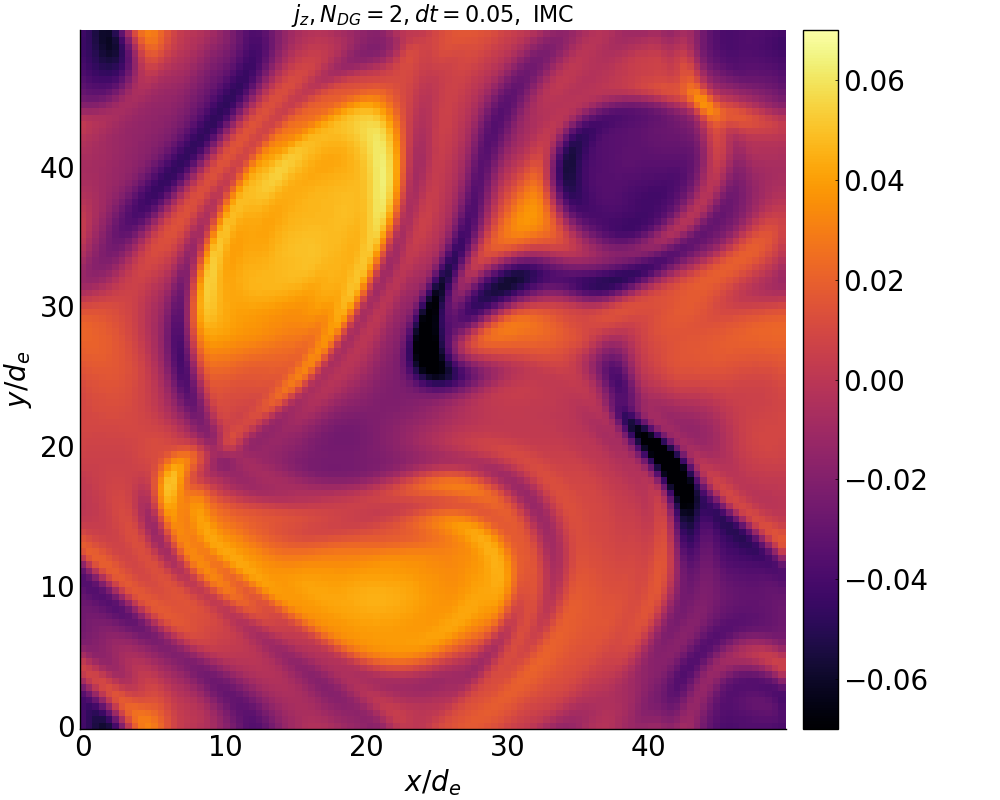}
  \includegraphics[width=0.24\linewidth]{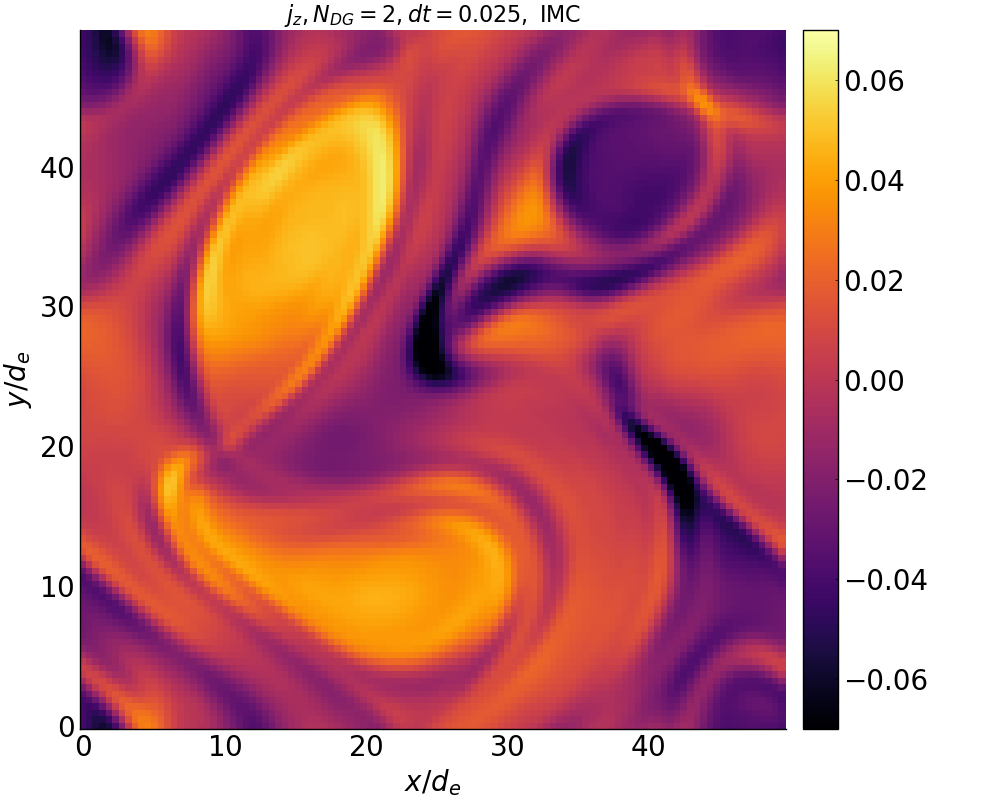} \\
   \includegraphics[width=0.24\linewidth]{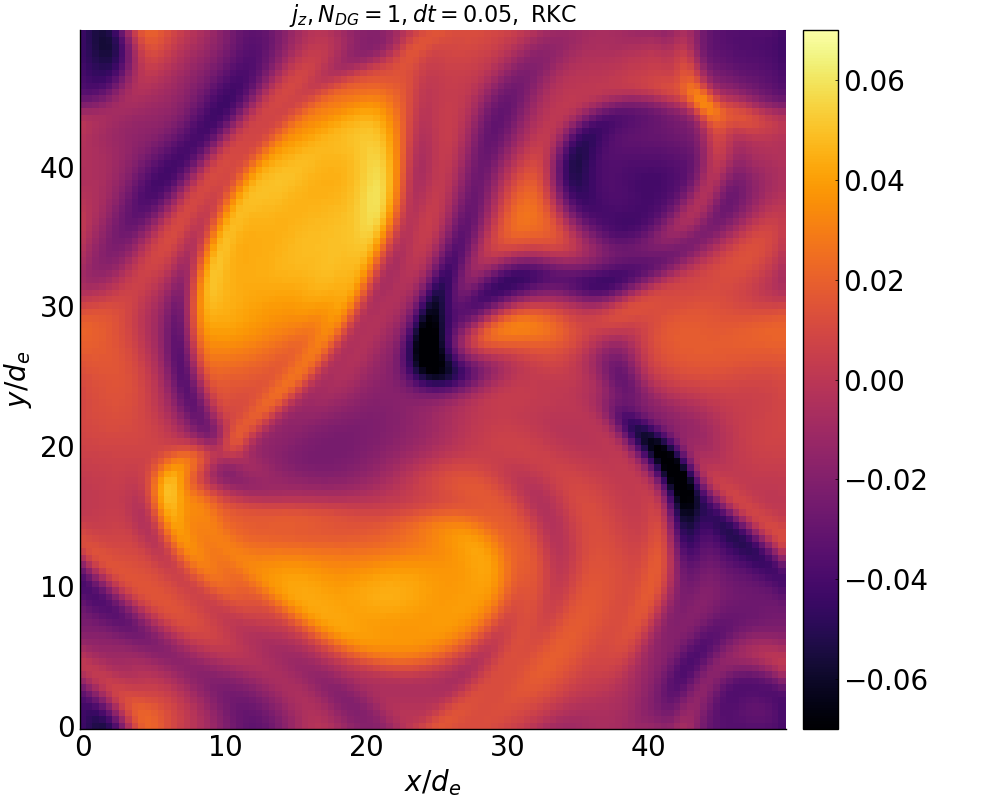}
  \includegraphics[width=0.24\linewidth]{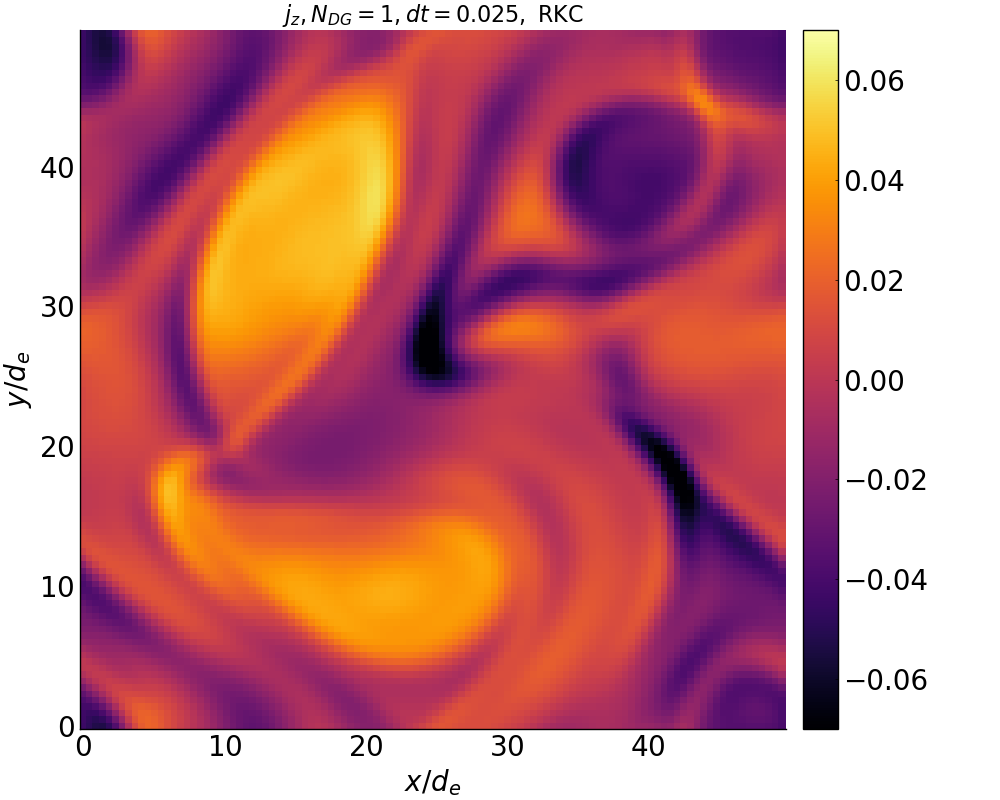}
   \includegraphics[width=0.24\linewidth]{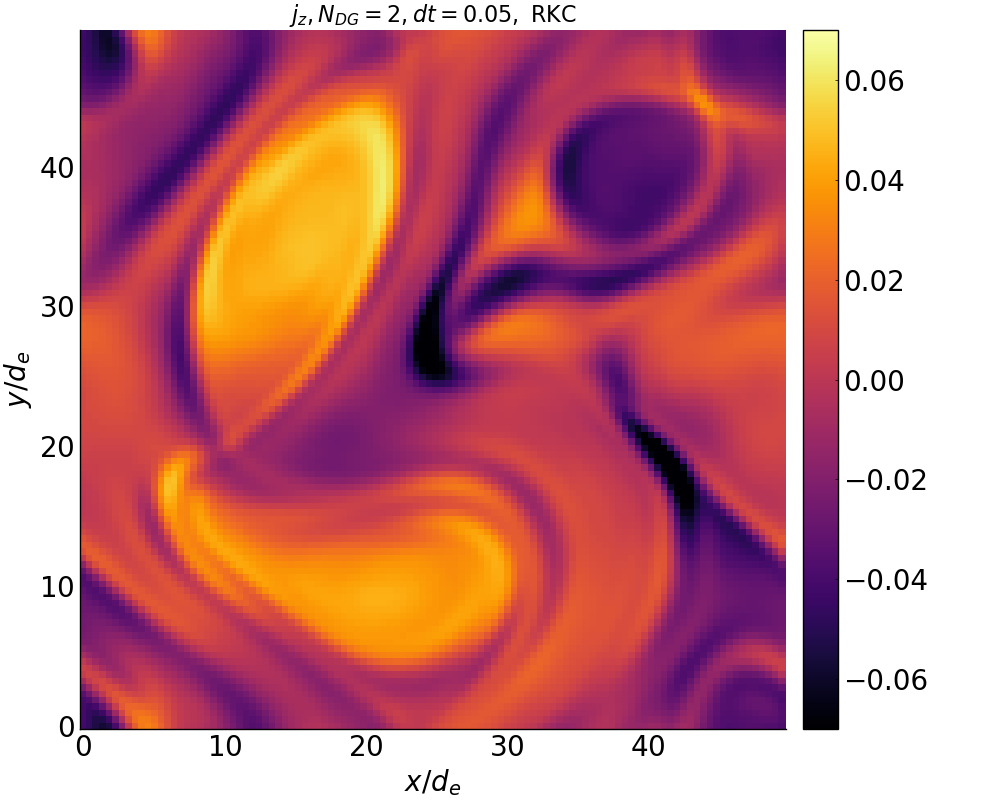}
  \includegraphics[width=0.24\linewidth]{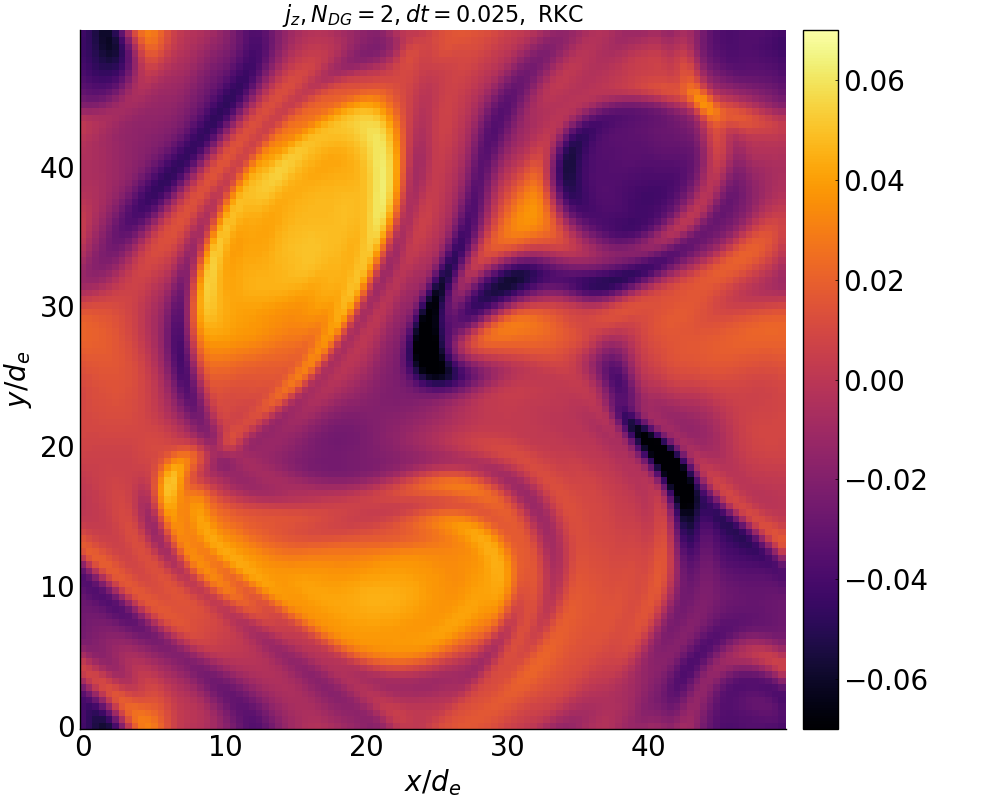} \\
   \includegraphics[width=0.24\linewidth]{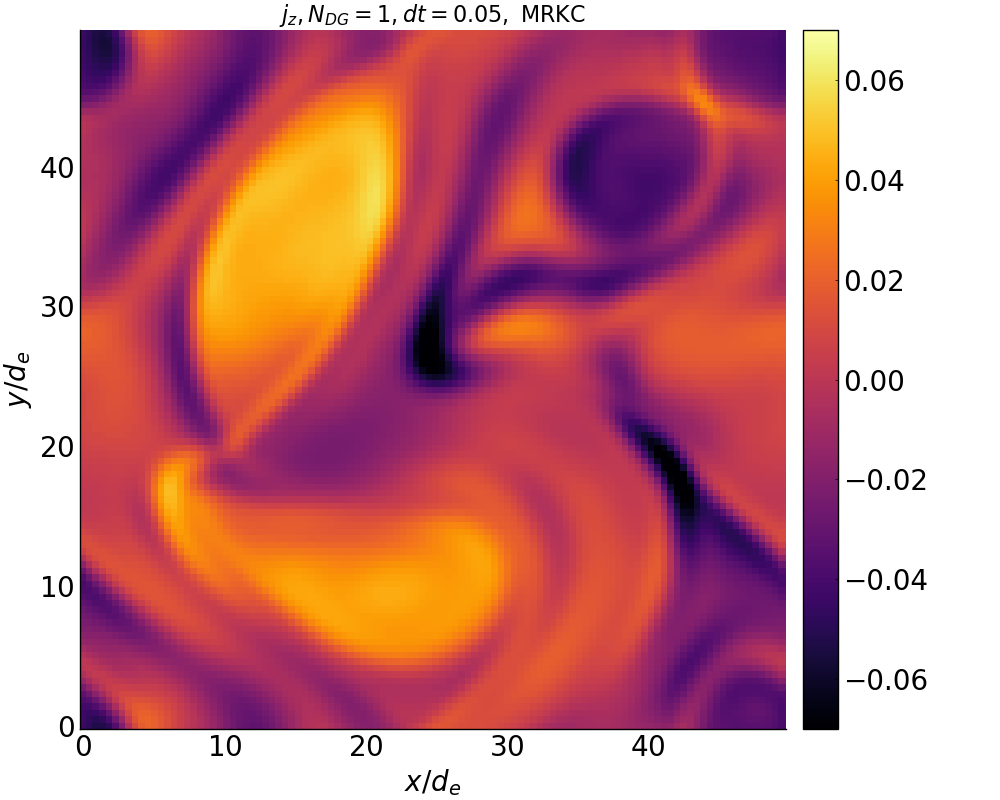}
  \includegraphics[width=0.24\linewidth]{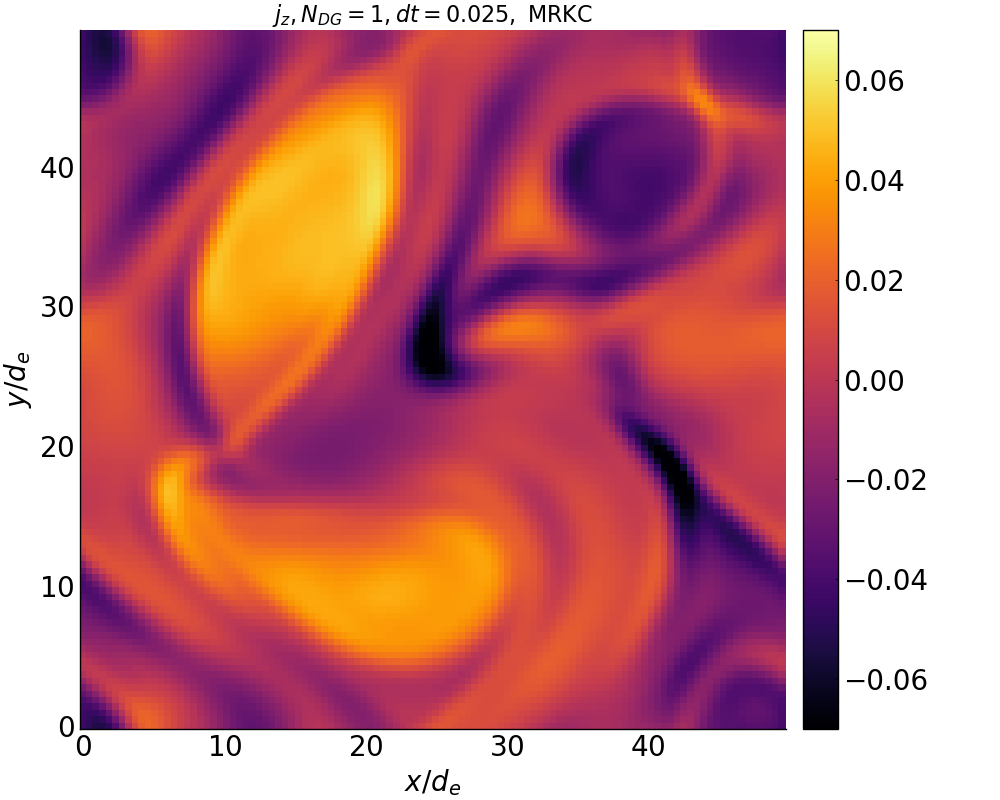}
   \includegraphics[width=0.24\linewidth]{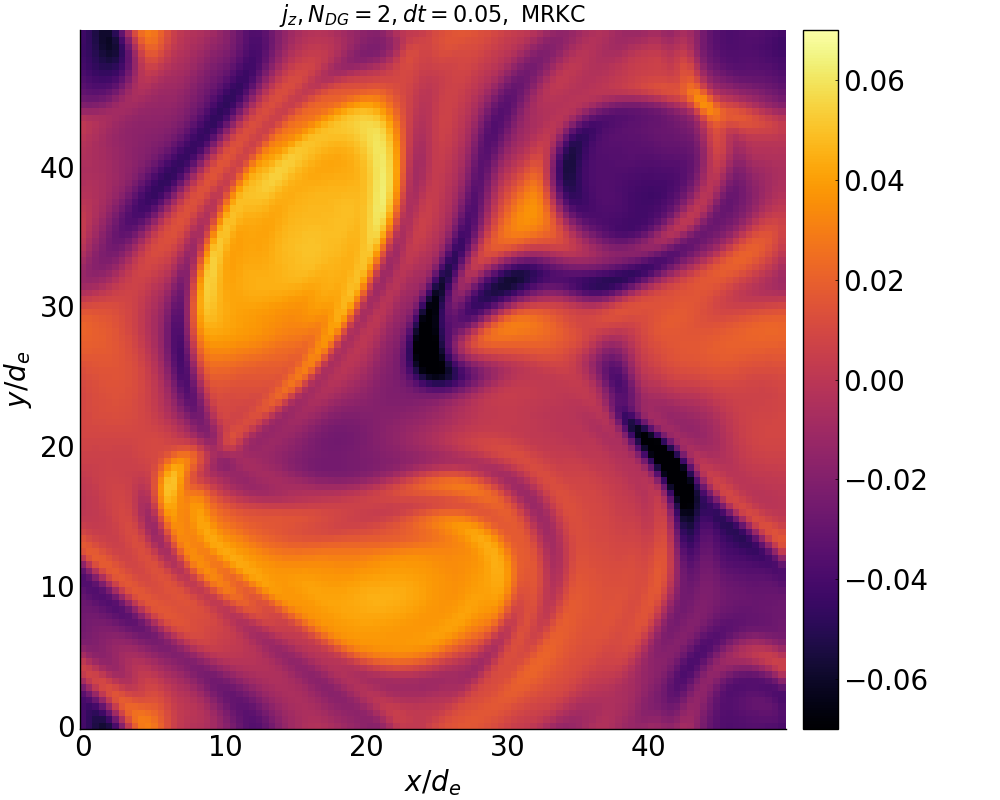}
  \includegraphics[width=0.24\linewidth]{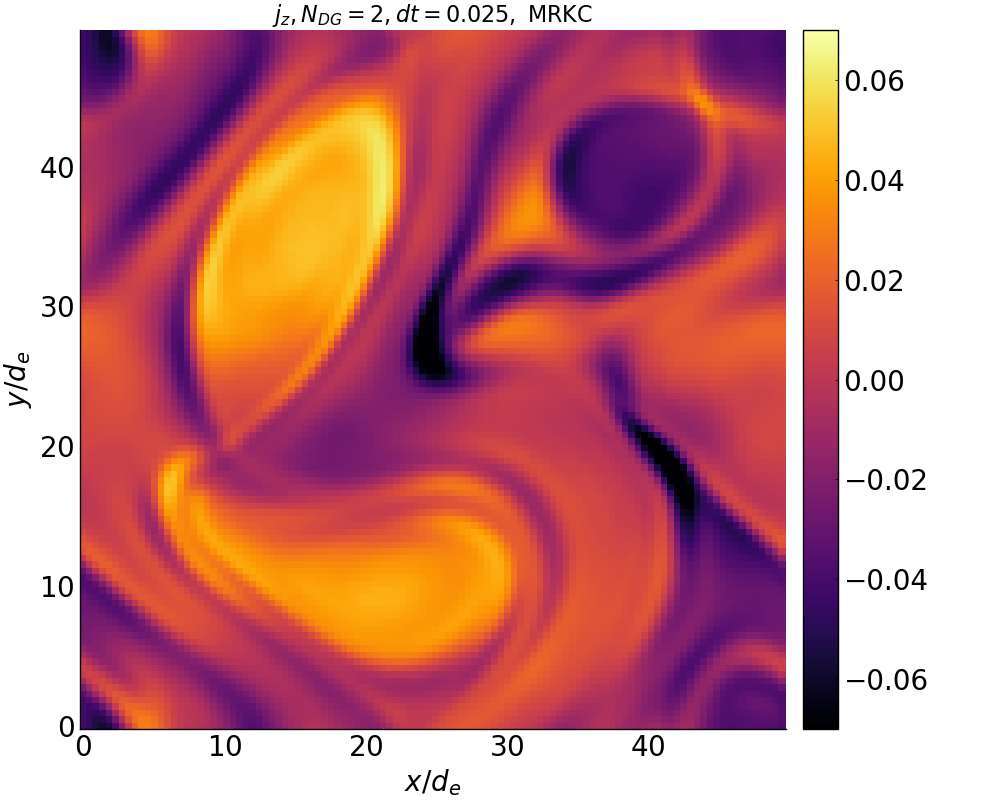} \\
   \includegraphics[width=0.24\linewidth]{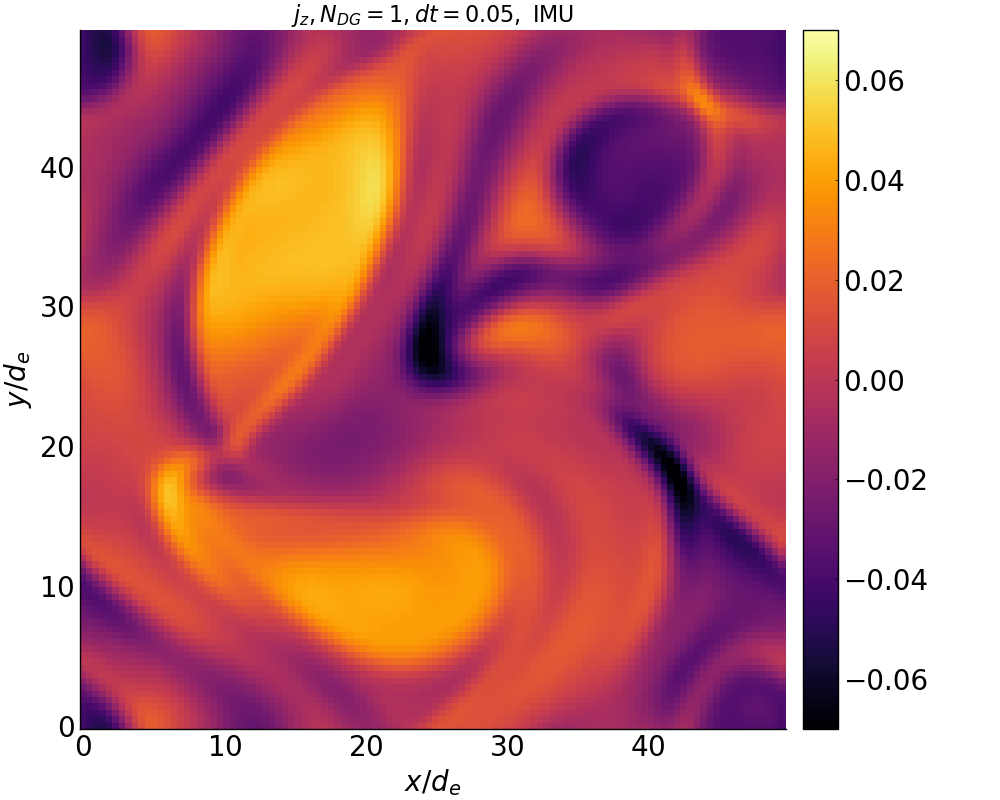}
  \includegraphics[width=0.24\linewidth]{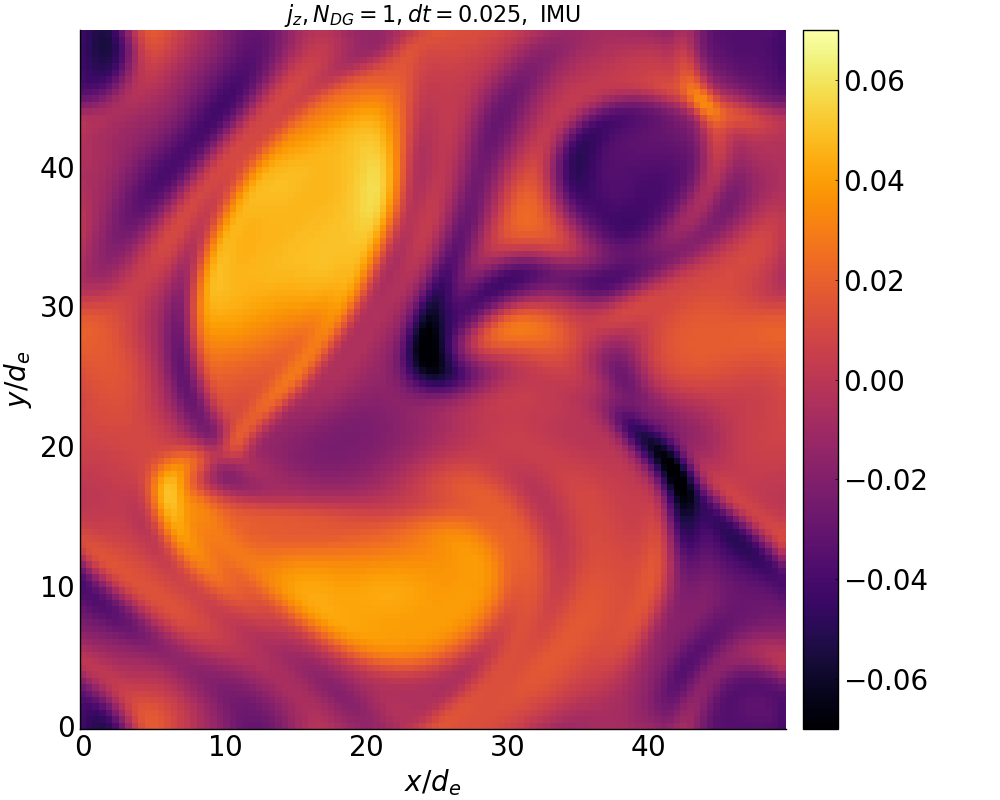}
   \includegraphics[width=0.24\linewidth]{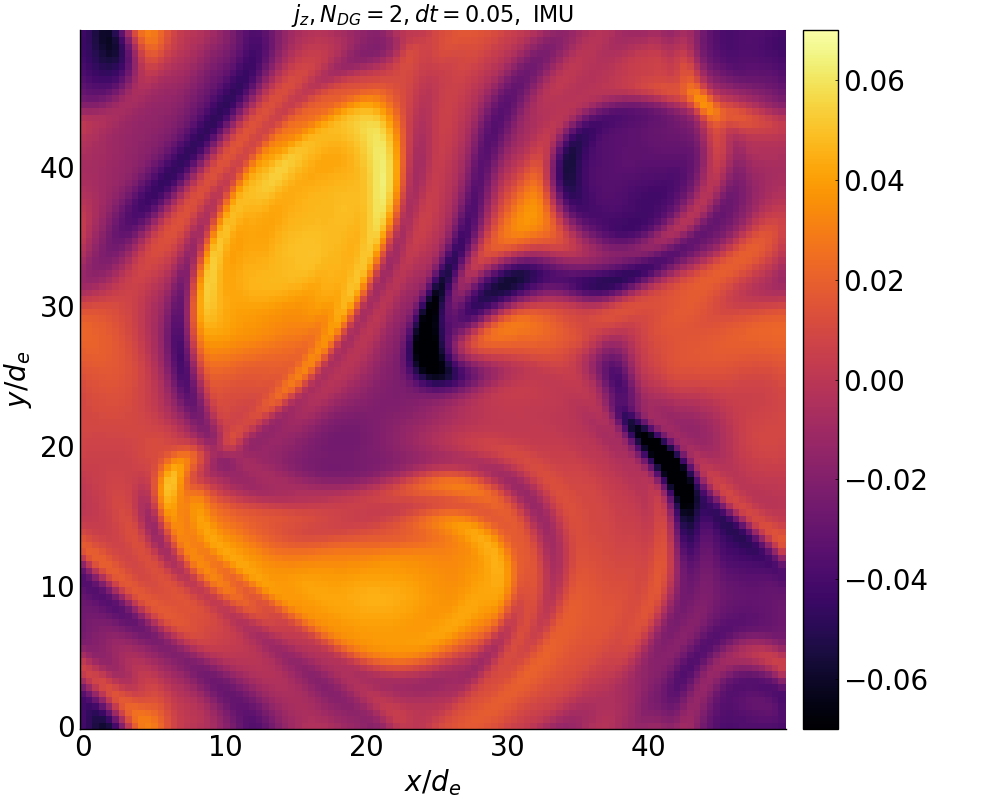}
  \includegraphics[width=0.24\linewidth]{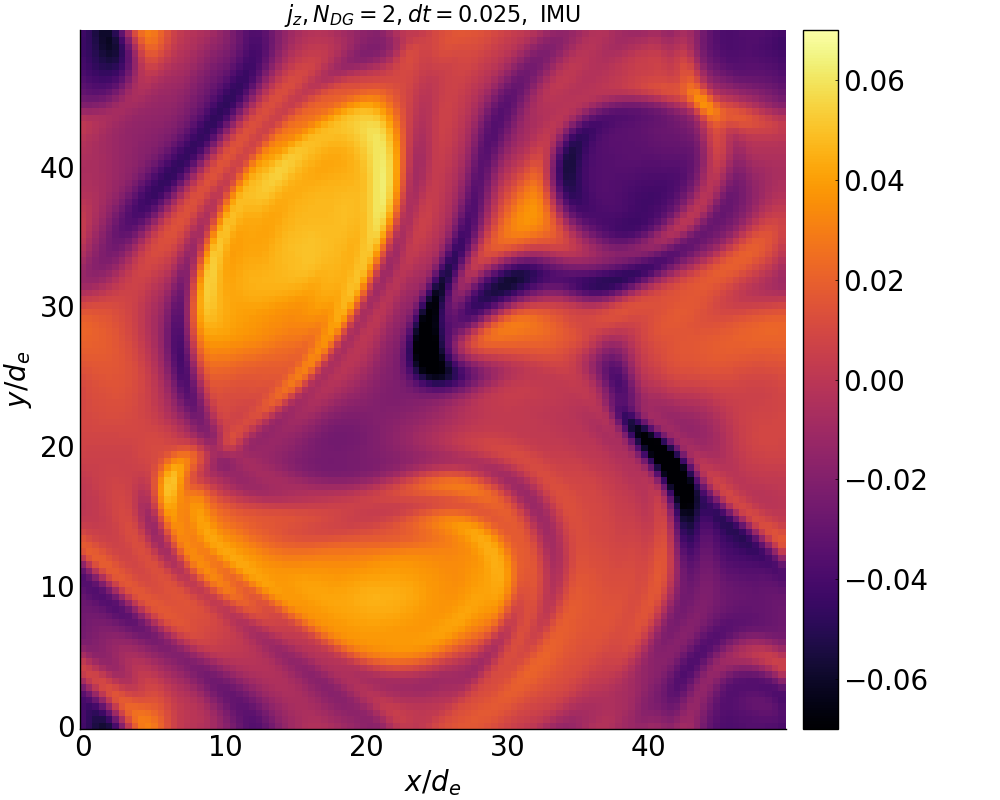} \\ 
   \includegraphics[width=0.24\linewidth]{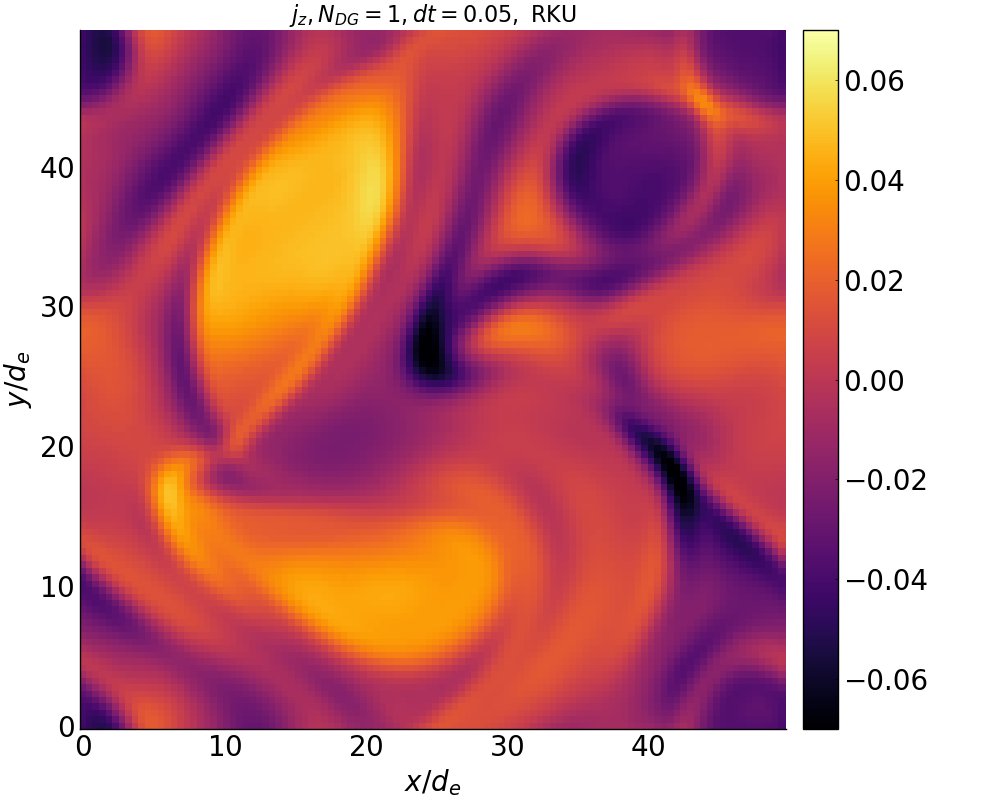}
  \includegraphics[width=0.24\linewidth]{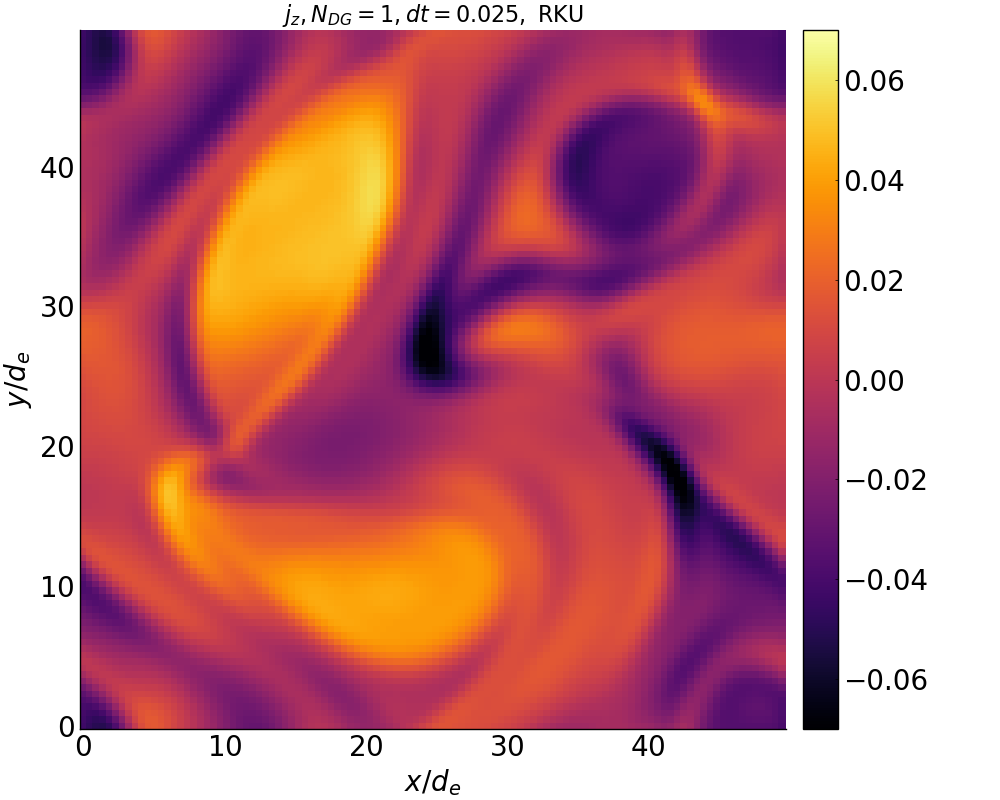}
   \includegraphics[width=0.24\linewidth]{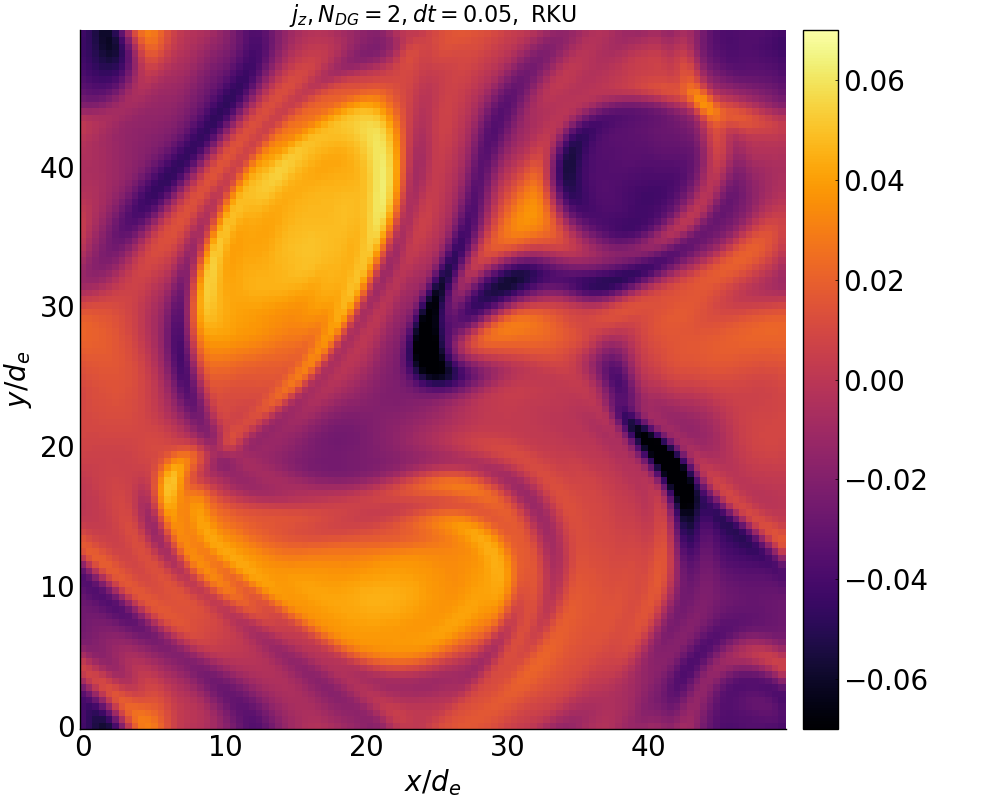}
  \includegraphics[width=0.24\linewidth]{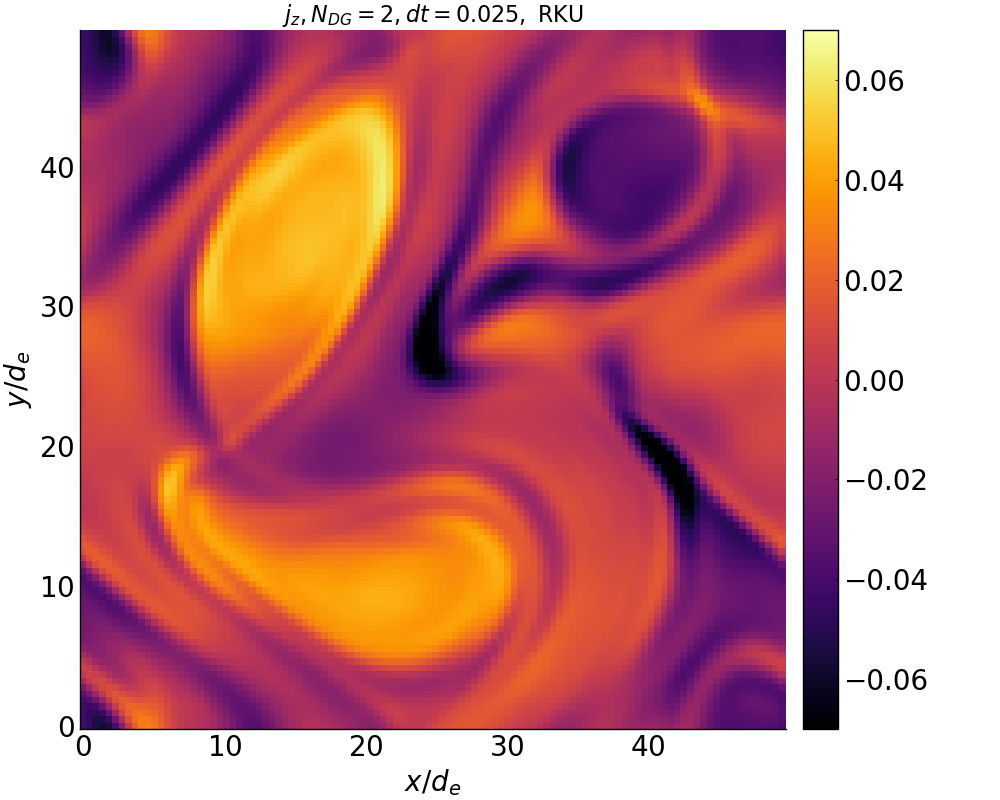} \\ 
   \includegraphics[width=0.24\linewidth]{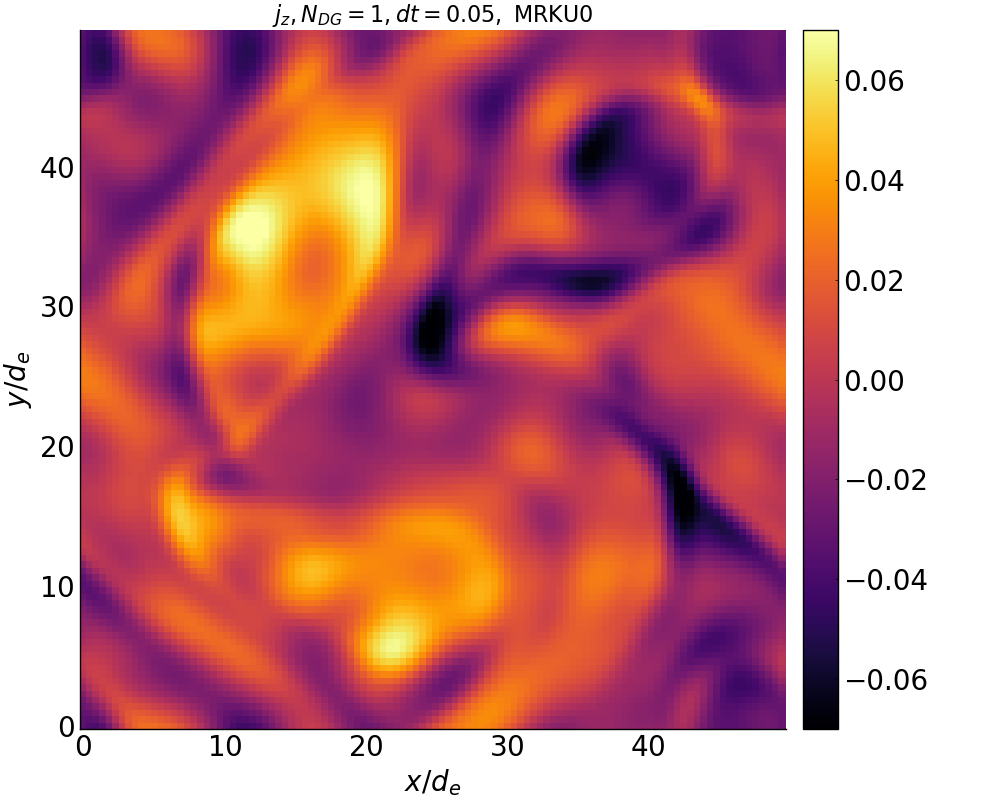}
  \includegraphics[width=0.24\linewidth]{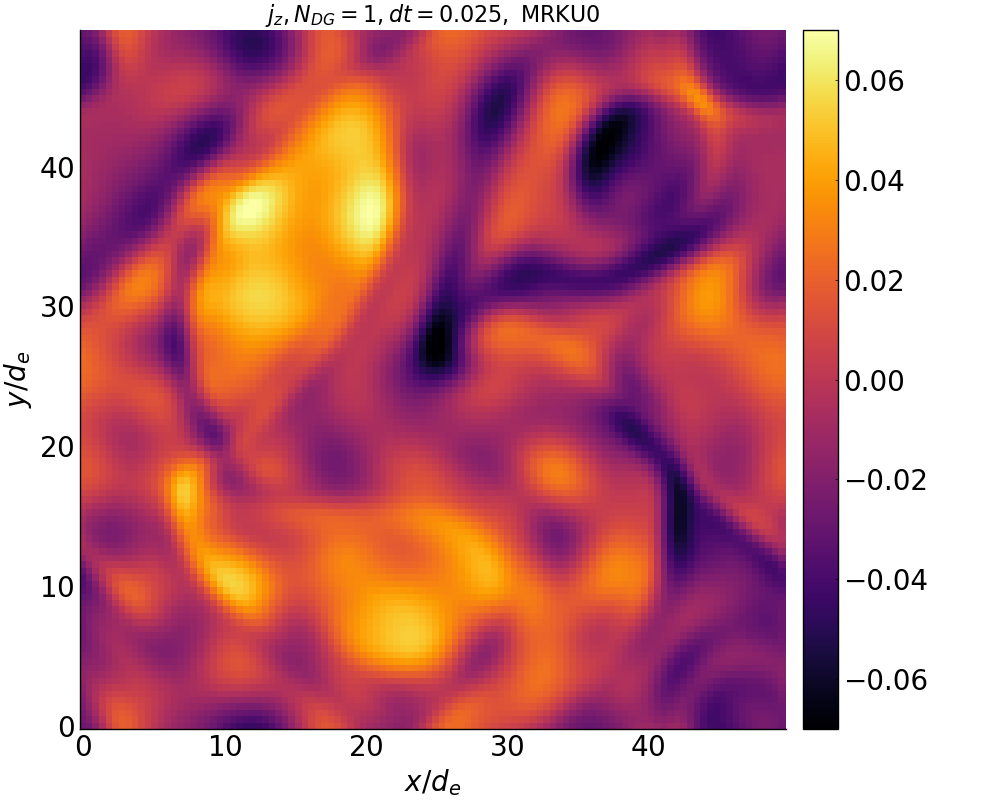}
   \includegraphics[width=0.24\linewidth]{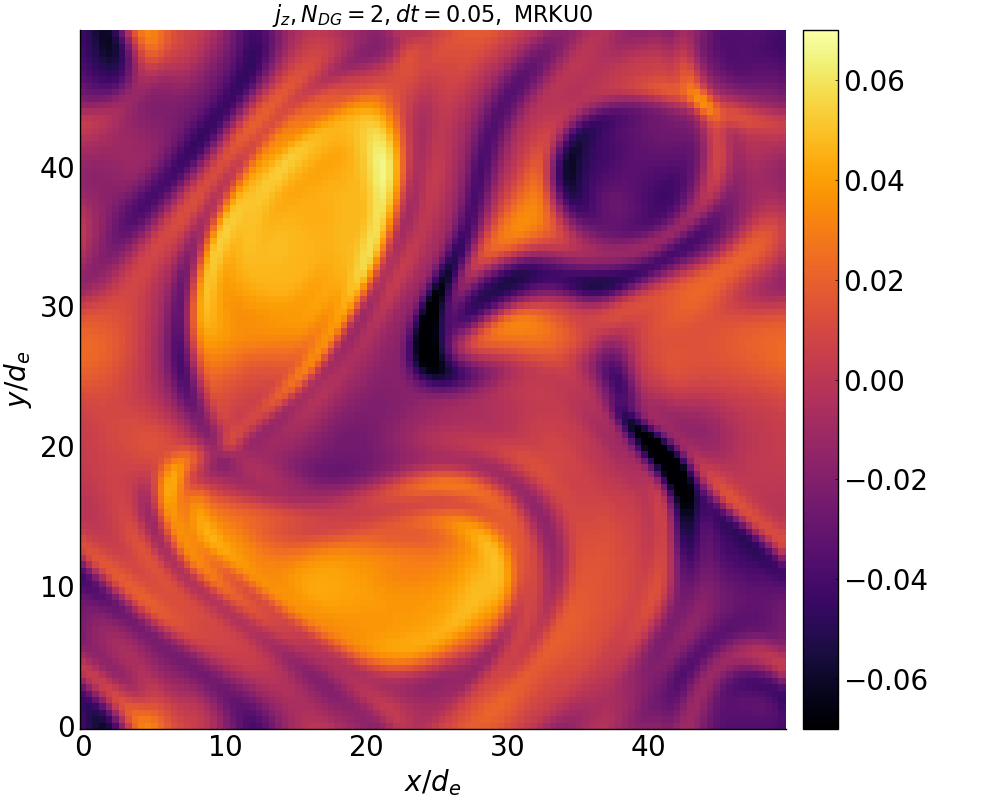}
  \includegraphics[width=0.24\linewidth]{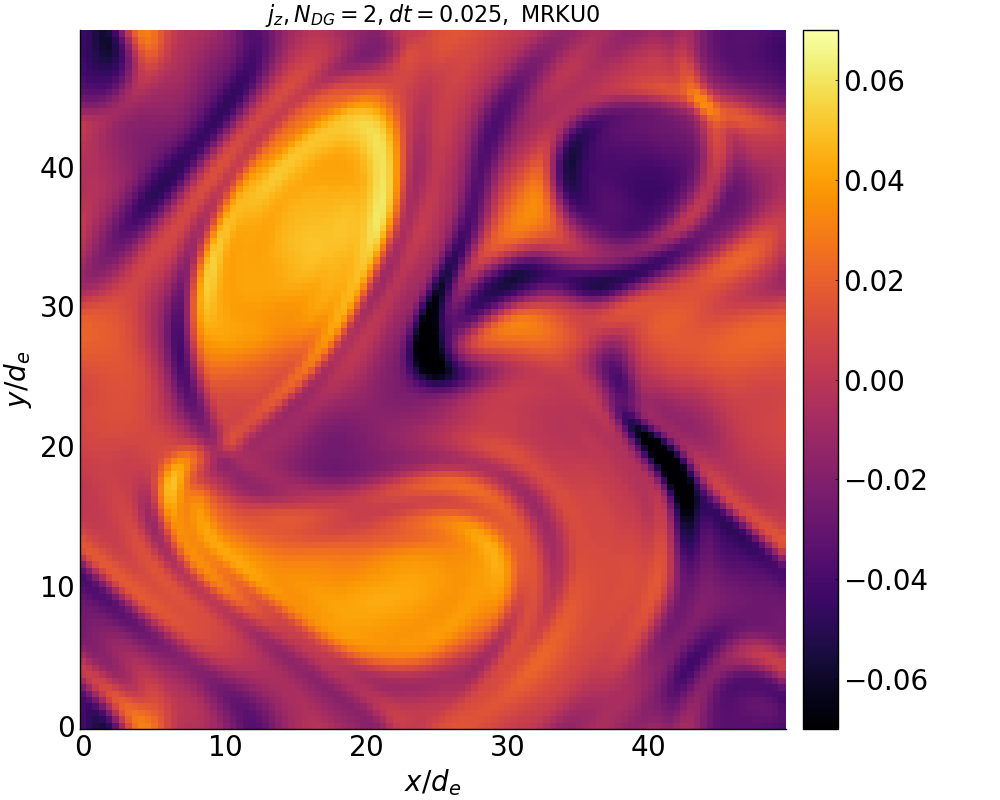} \\ 
   \includegraphics[width=0.24\linewidth]{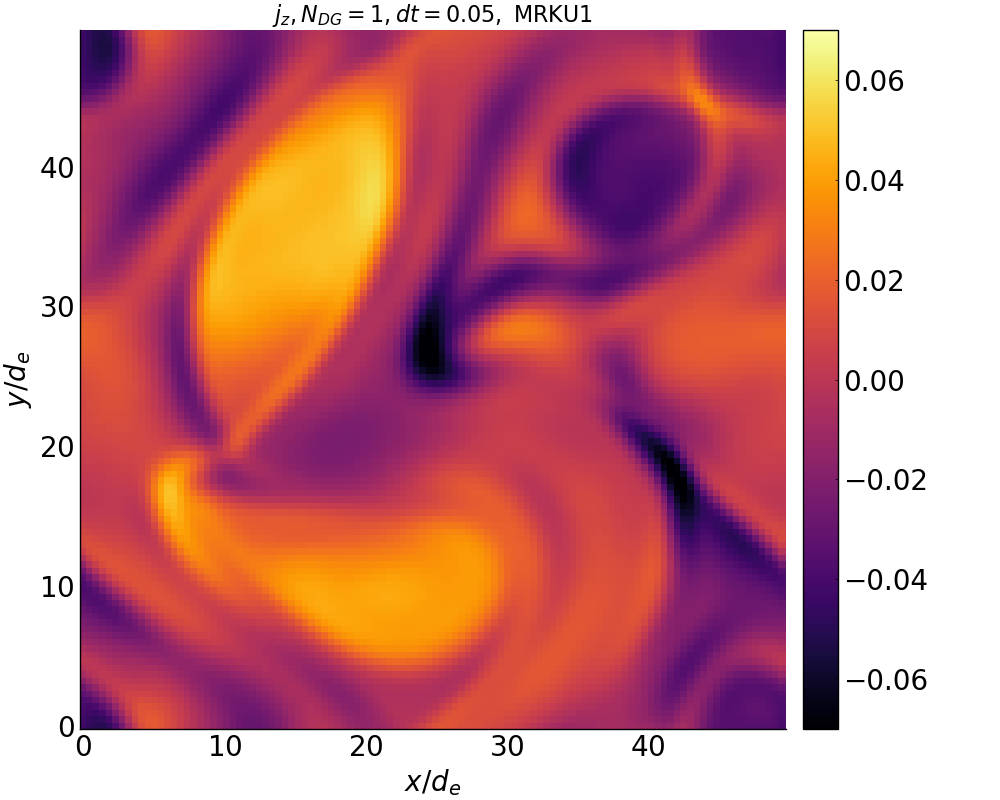}
  \includegraphics[width=0.24\linewidth]{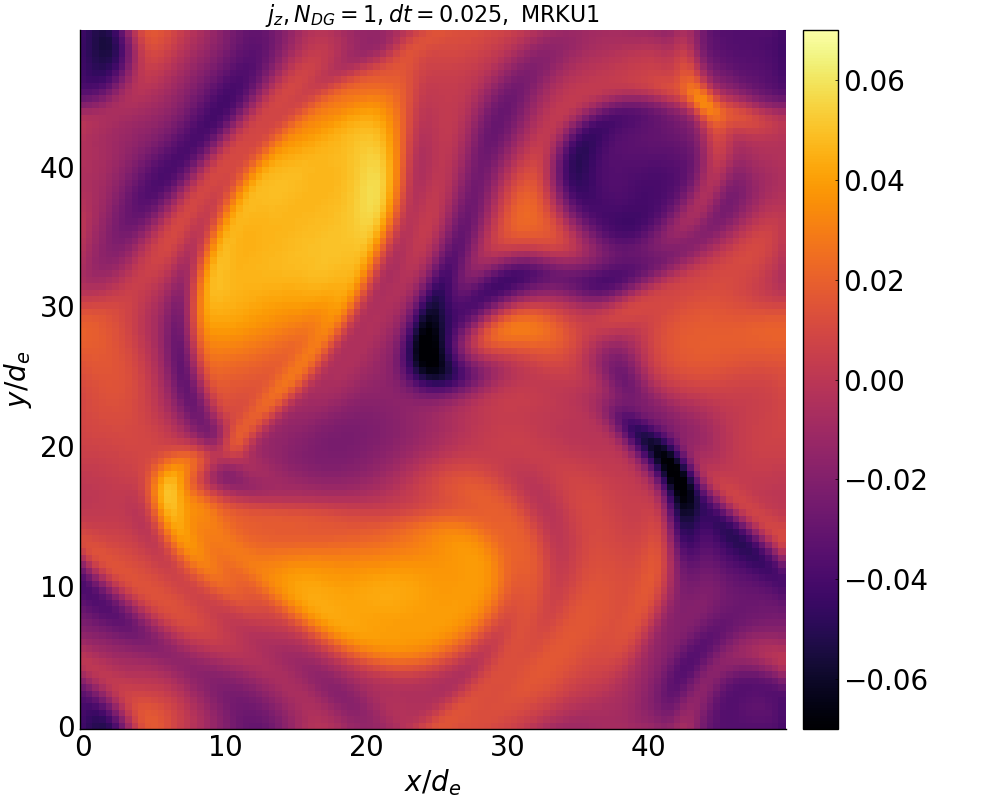}
   \includegraphics[width=0.24\linewidth]{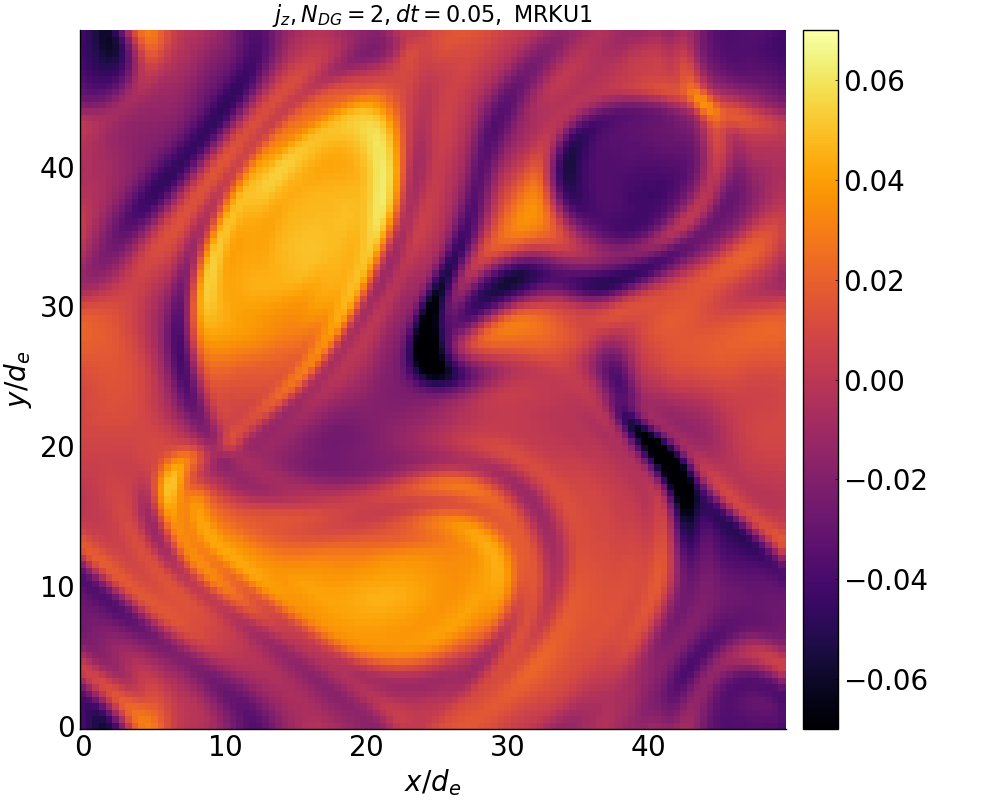}
  \includegraphics[width=0.24\linewidth]{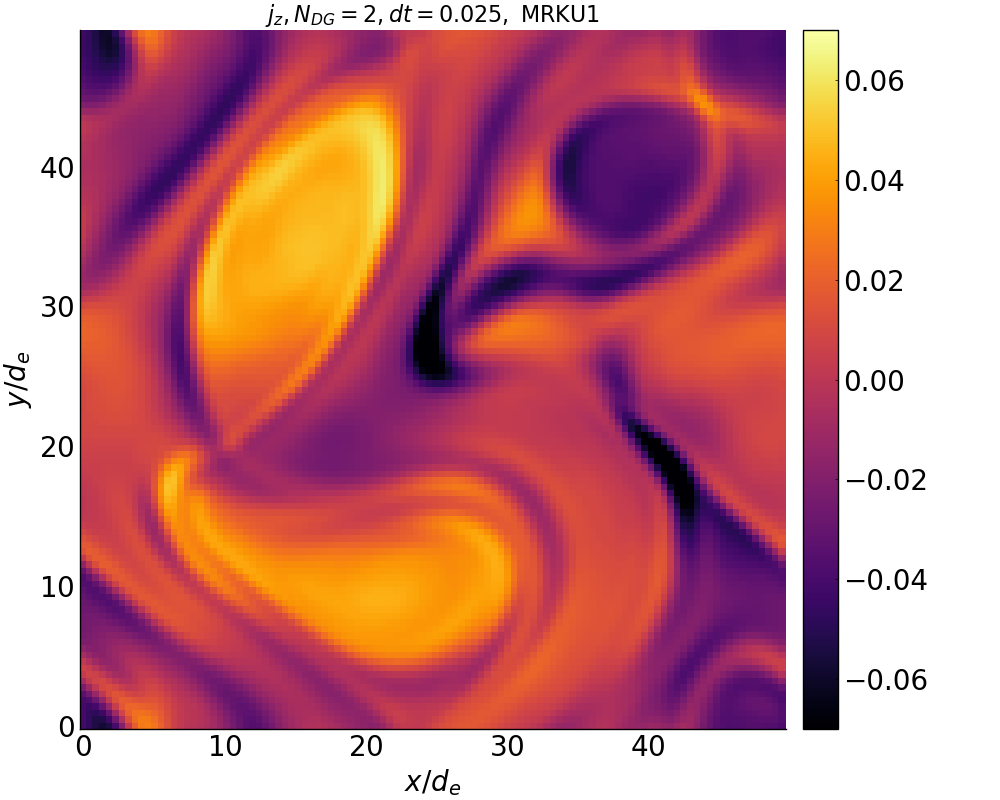}
   \caption{Orszag-Tang vortex benchmark: out of plane current at $t = 1000$
   for
   IMC  (first row), 
   RKC  (second row), 
   MRKC (third row), 
   IMU  (fourth row), 
   RKU (fifth row),
   MRKU0 (sixth row),
   MRKU1 (seventh row)
   with 
   $N_{DG} = 1$ (first and second column),
   $N_{DG} = 2$ (third and fourth column),
   $\Delta t = 0.05$ (first and third column),
   \rev{$\Delta t = 0.025$} (second and fourth column).}
  \label{fig:OT:jz}
\end{figure}

Similarly to the previous numerical tests, we study the seven 
methods summarized in Table~\ref{tab:methods}.
%
As a qualitative measure of the performances of the proposed methods,
we report in Fig.~\ref{fig:OT:jz} the
plasma current along the $z$-axis at the final simulation time $t=1000$.
All runs with $N_{DG}=2$ (third and fourth columns) show 
visually indistinguishable solutions. 
We note that further increasing the temporal, spatial, or velocity resolution
does not lead to visible changes in the results (more refined runs are not reported here).
All runs with $N_{DG}=1$ (first and second columns), except for the MRKU0 method,
produce qualitatively similar plots, but slightly smeared due to the
spatial discretization error.
The MRKU0 method with $N_{DG}=1$ (sixth row, first and second columns)
produce results polluted by high amplitude oscillations, which make the results unreliable.
The nature of those oscillations could be related to the unphysical modes
which were present in the X-mode test in the previous section, 
see the bottom left panel of Fig.~\ref{fig:xmode:u}.

\begin{figure}[H]
  \centering
  \includegraphics[width=0.31\linewidth]{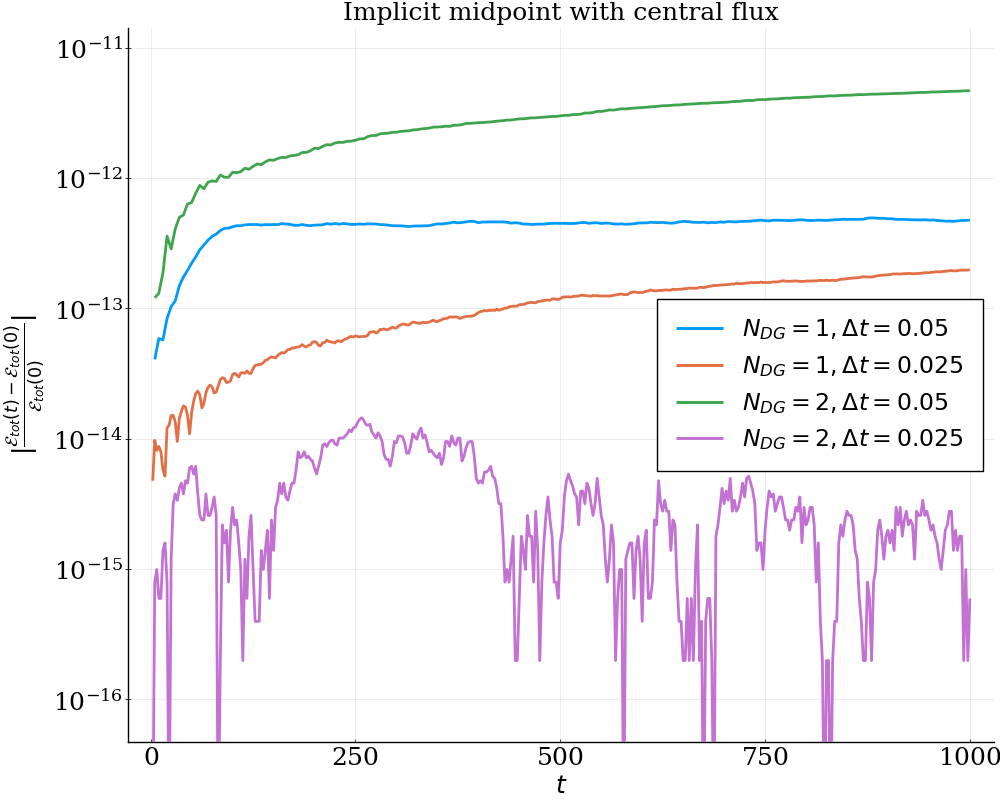}
  \includegraphics[width=0.31\linewidth]{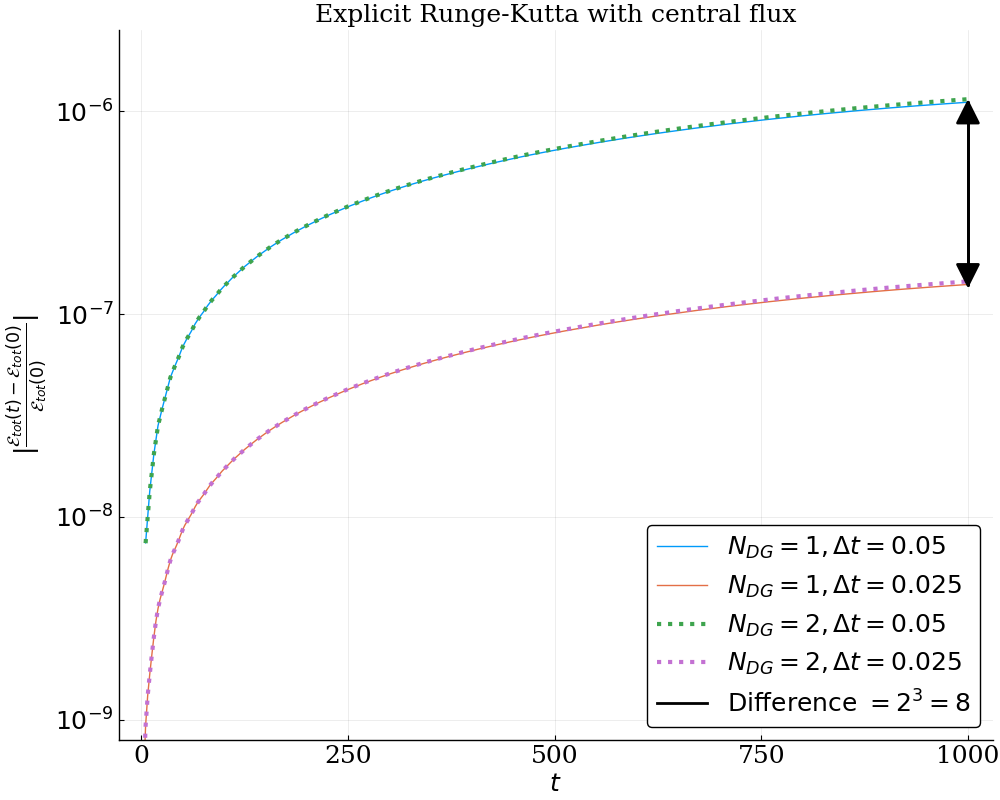}
  \includegraphics[width=0.31\linewidth]{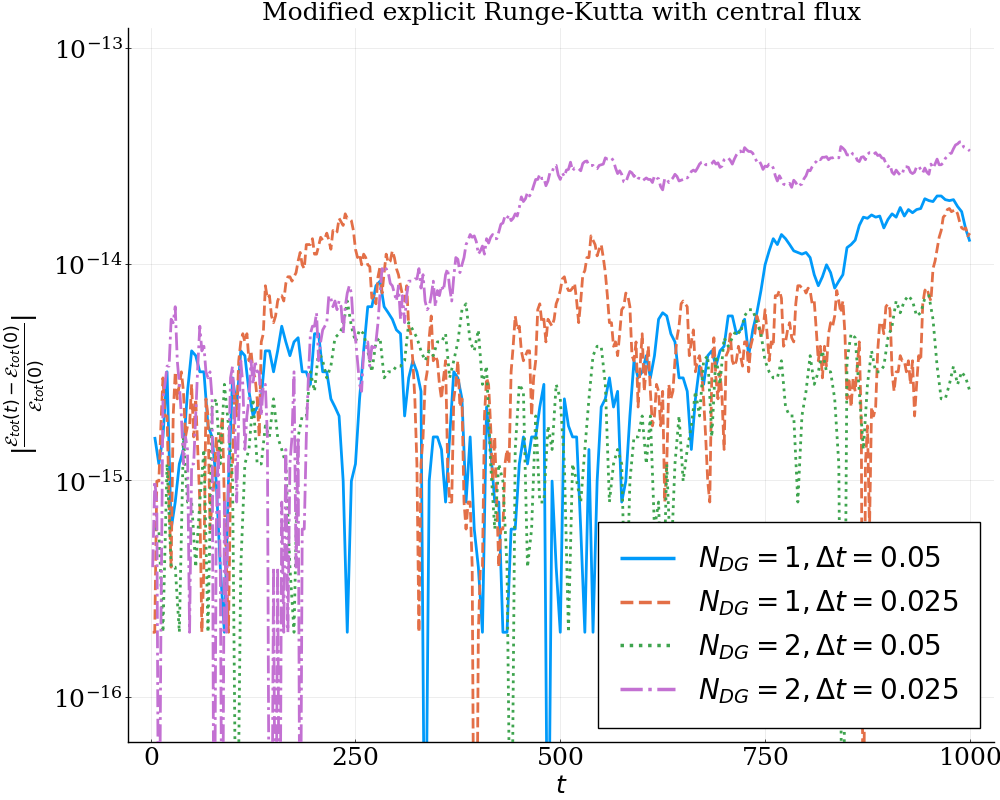}
  \caption{Orszag-Tang vortex benchmark: time evolution of relative error 
  in energy conservation for IMC (left), RKC (middle), and MRKC (right)
  for different spatial and temporal resolution.}
  \label{fig:OT:central:energy}
\end{figure}

\begin{figure}[H]
  \centering
  \includegraphics[width=0.31\linewidth]{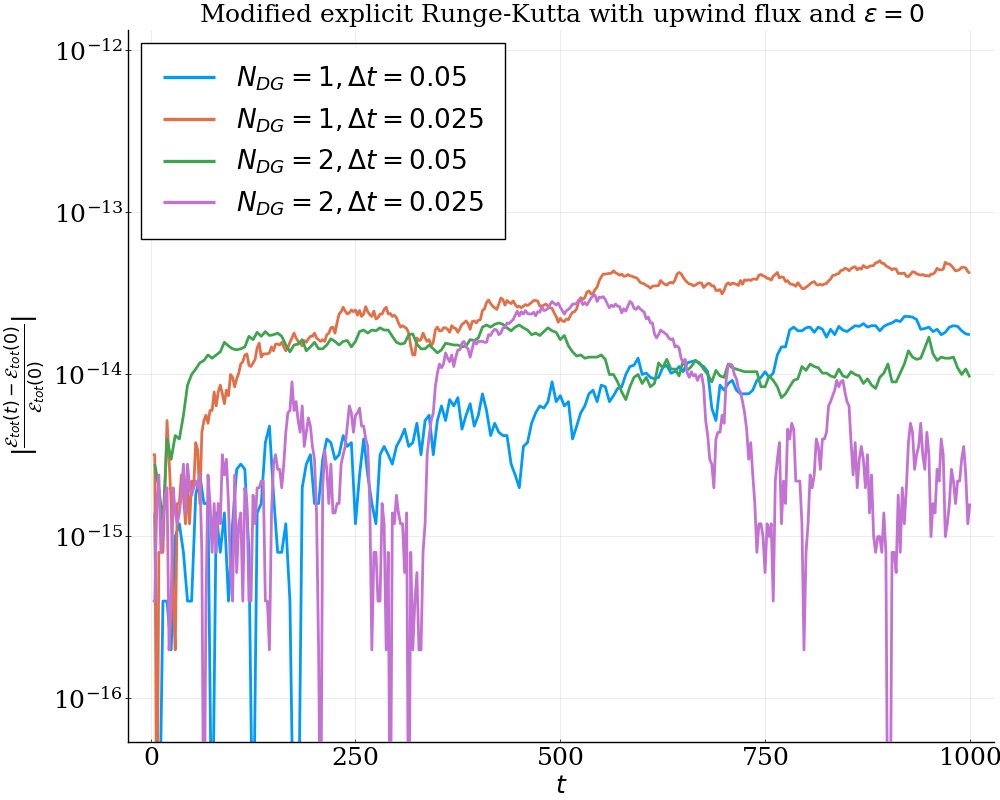}
  \includegraphics[width=0.31\linewidth]{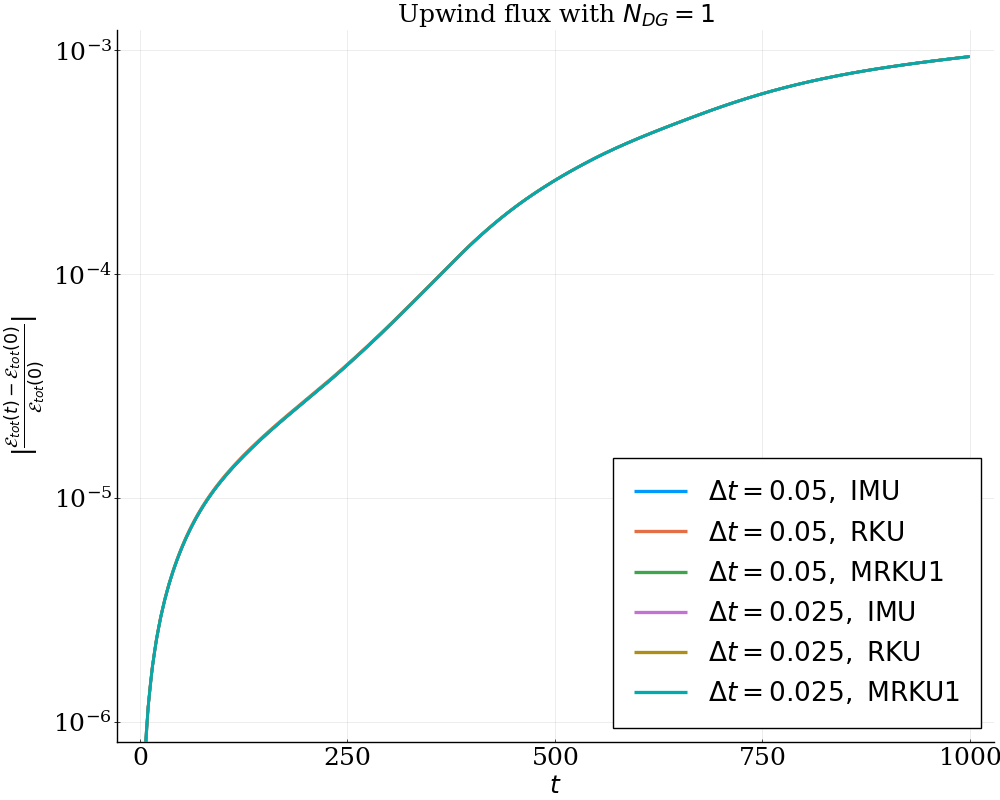}
  \includegraphics[width=0.31\linewidth]{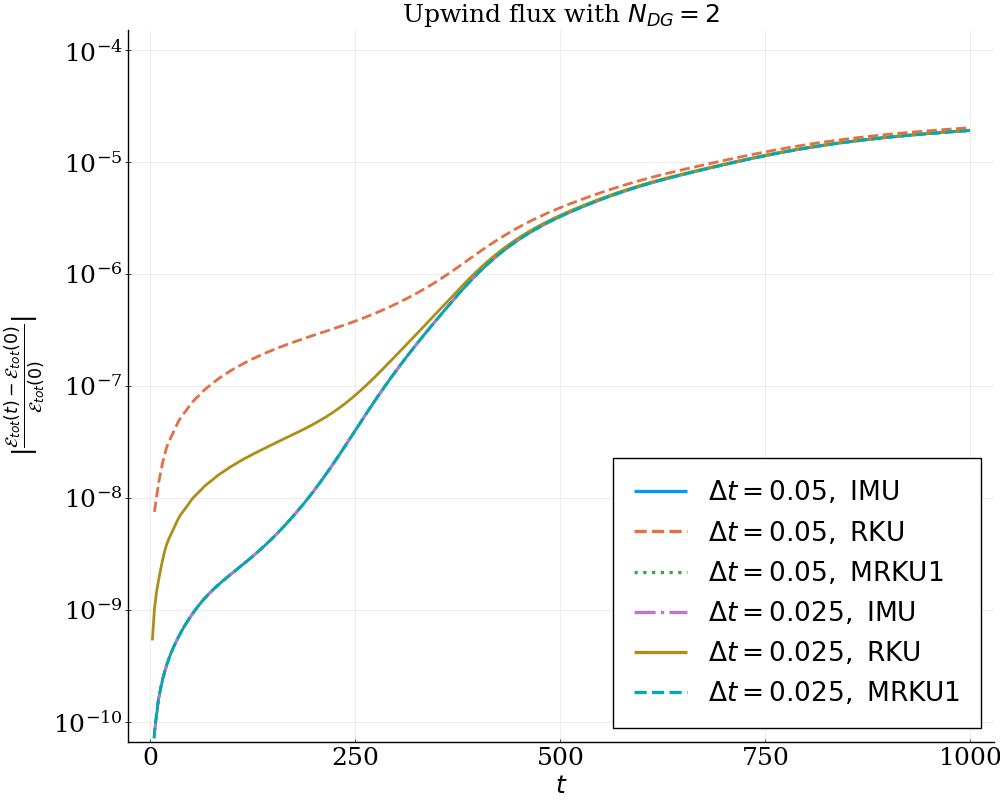}
  \caption{Orszag-Tang vortex benchmark: time evolution of relative error 
  in energy conservation for 
  MRKU0 (left), 
      IMU, RKU, MRKU1 with $N_{DG}=1$ (middle),
  and IMU, RKU, MRKU1 with $N_{DG}=2$ (right),
  for different resolution.}
  \label{fig:OT:upwind:energy}
\end{figure}

Concerning the energy conservation in the Orszag-Tang vortex test,
we distinguish the case of central fluxes in the spatial approximation of Maxwell's equations (Fig.~\ref{fig:OT:central:energy}) and
upwind fluxes (Fig.~\ref{fig:OT:upwind:energy}).
The results for central fluxes are similar to the results obtained with the same numerical methods in the whistler instability test
in Fig.~\ref{fig:whistler_x:central:energy}.
In the left panel of Fig.~\ref{fig:OT:central:energy}, the error in the conservation of energy for the IMC scheme is shown:
energy is conserved up to the
nonlinear solver tolerance.
The middle panel illustrates the performances of the RKC scheme,
where the error in the energy conservation
depends on the accuracy of the temporal discretization
(the error is independent on $N_{DG}$ and scales as the Runge-Kutta order
producing a reduction of the error by a factor of $\sim8 = 2^3$ shown 
with the two sided black arrow).
The right panel shows the results obtained with the MRKC scheme, which confirm that the conservation of energy is satisfied to machine precision.

The results for the methods based on upwind fluxes are shown in Fig.~\ref{fig:OT:upwind:energy}, and
are similar to the ones obtained in the whistler instability test in Fig.~\ref{fig:whistler_x:upwind:energy}.
The MRKU0 method (left panel) conserves
the energy to machine precision although the plot of the current in Fig.~\ref{fig:OT:jz}
shows unsatisfactory results (at least for $N_{DG}=1$).
All other upwind methods, i.e., IMU, RKU, MRKU1 dissipate energy.
The middle panel of Fig.~\ref{fig:OT:upwind:energy}
has results with $N_{DG}=1$ and the right panel has $N_{DG}=2$.
As before, the effect of the temporal discretization is small compared to the spatial discretization error for $N_{DG}=1$ and the curves associated with different methods and different $\Delta t$ coincide. For $N_{DG}=2$ (right panel), the energy error in some of the methods differ for $t<400$. In general, RKU has an energy error that is one to two orders of magnitude larger than that of IMU and MRKU1, showing that in this case the correction for upwind based Runge-Kutta methods is important. The correct third order scaling of RKU is also recovered, as shown by comparing the curves with $\Delta t=0.05$ and $\Delta t = 0.025$. The energy errors for IMU and MRKU1 essentially coincide, irrespective of $\Delta t$. For $t > 400$, the energy error becomes dominated by the spatial discretization error and all the curves essentially overlap (cf. Fig. \ref{fig:whistler_x:upwind:energy}).

\begin{figure}[H]
  \centering
  \includegraphics[width=0.31\linewidth]{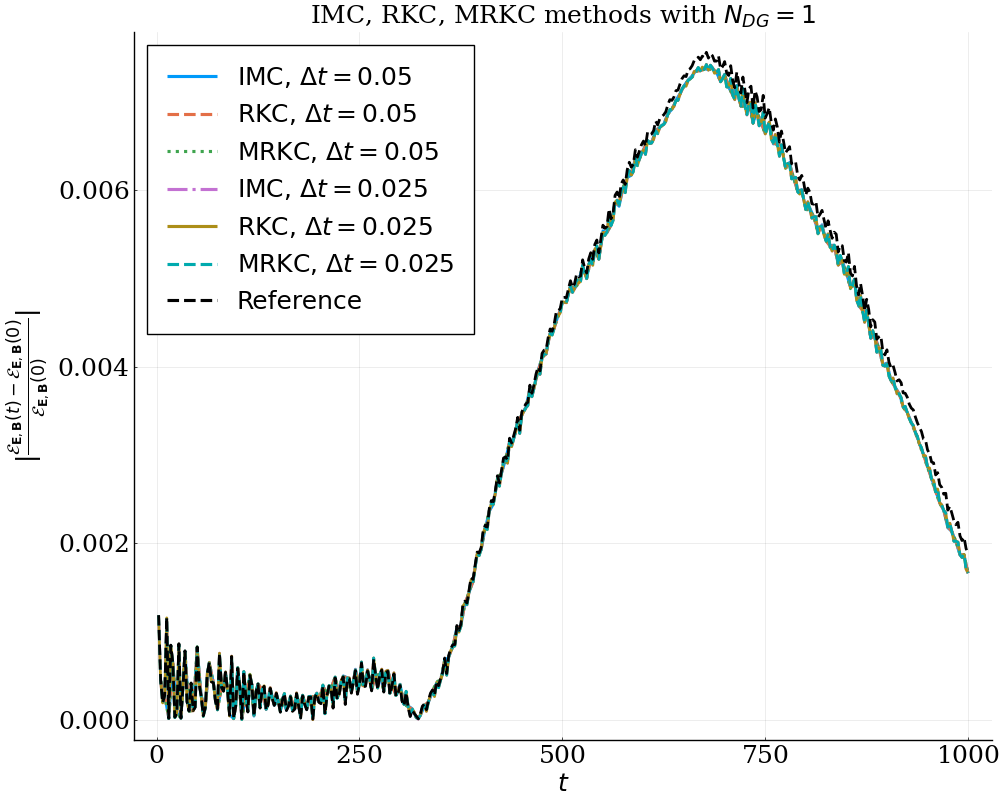}
  \includegraphics[width=0.31\linewidth]{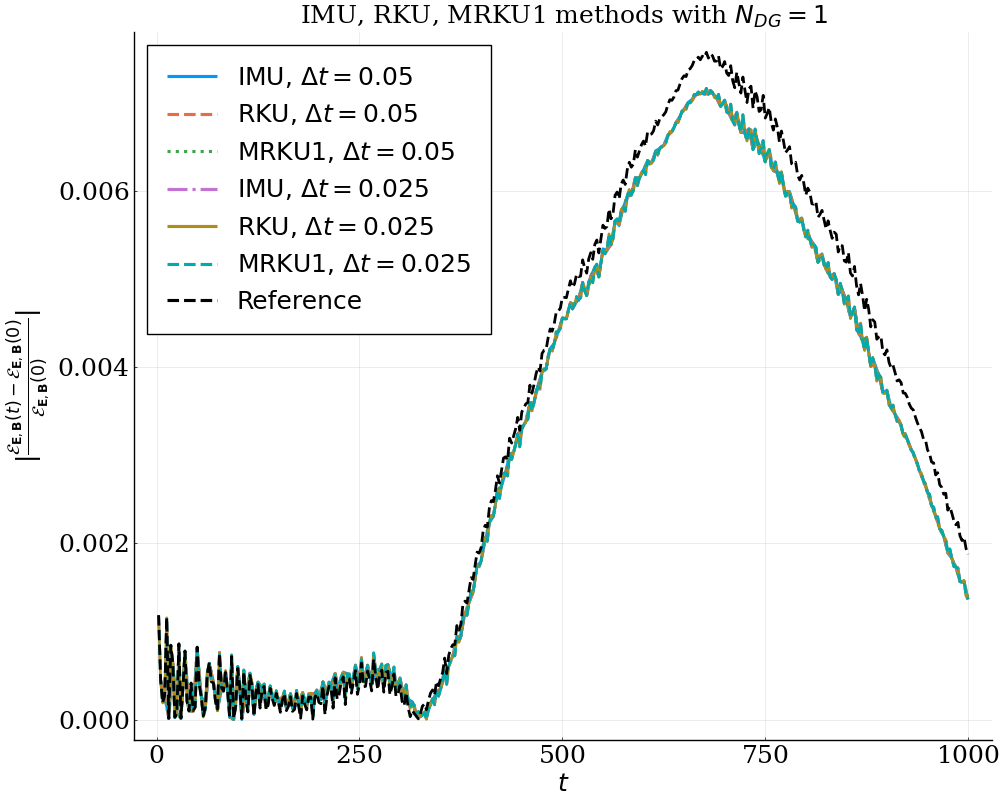}
  \includegraphics[width=0.31\linewidth]{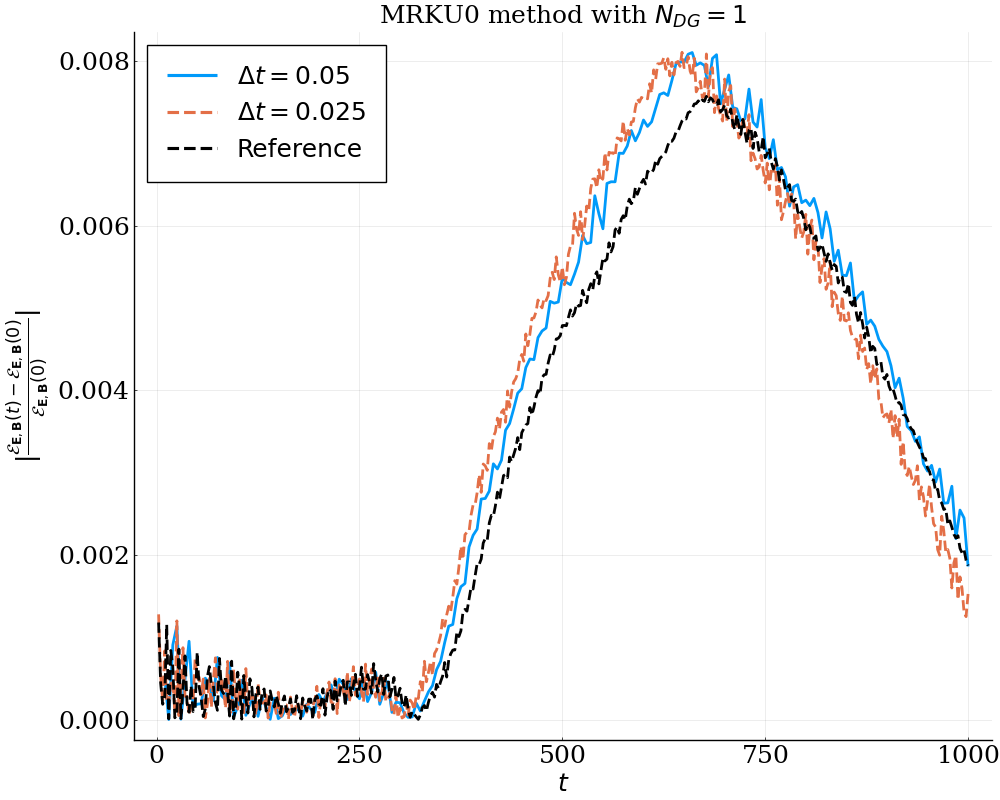}
  \caption{Orszag-Tang vortex benchmark: time evolution of relative change 
  in electromagnetic energy $\mathcal{E}_{\mathbf{E},\mathbf{B}}$ 
  for $N_{DG}=1$
  and 
  IMC, RKC, MRKC  (left),
  IMU, RKU, MRKU1 (middle),
  MRKU0 (right),
  for different $\Delta t$
  against reference solution (IMC with $N_{DG}=2$, $\Delta t = 0.025$).}
  \label{fig:OT:EBenergy}
\end{figure}

Finally, we look at how the electromagnetic energy  $\calE_{\ELE}$ 
defined in (\ref{eq:calEBE}) changes over time for various methods. 
The evolution of the electromagnetic energy for runs with $N_{DG}=2$
looks essentially identical for all methods and time steps
(not shown here).
Therefore, we use a run with $N_{DG}=2$, $\Delta t = 0.025$
and the IMC method as a reference in the following results.
The left panel in Fig.~\ref{fig:OT:EBenergy} shows
the evolution of $\calE_{\ELE}$ for runs with central flux
and $N_{DG}=1$. 
One can see that the central flux based methods reproduce the correct dynamics
of the electromagnetic energy (as seen by comparing them against the more resolved reference solution).
At the same time, 
IMU, RKU, and MRKU1 with $N_{DG}=1$ 
consistently dissipate electromagnetic energy
for different time step sizes,
as shown in the middle panel in Fig.~\ref{fig:OT:EBenergy}.
The right panel in Fig.~\ref{fig:OT:EBenergy}
shows the evolution of $\calE_{\ELE}$ for the MRKU0 method with $N_{DG}=1$.
This method conserves the total energy but it overestimates the
energy stored in the electromagnetic field. 
This additional energy is probably related to the
unphysical small scale fluctuations observed in the
current plots in Fig.~\ref{fig:OT:jz}.

\medskip

%% file: results_Mom.tex
\subsection{Notes on momentum conservation and positivity}
\label{sec:mom}
\begin{figure}[H]
  \centering
  \includegraphics[width=0.29\linewidth]{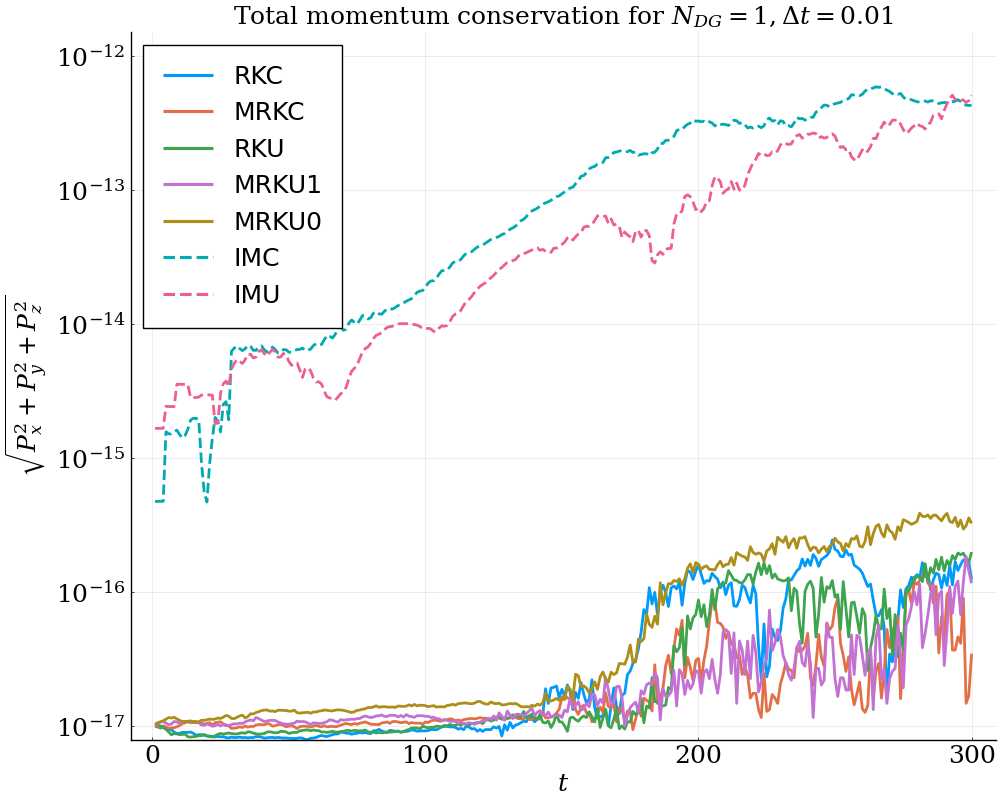}
  \includegraphics[width=0.35\linewidth]{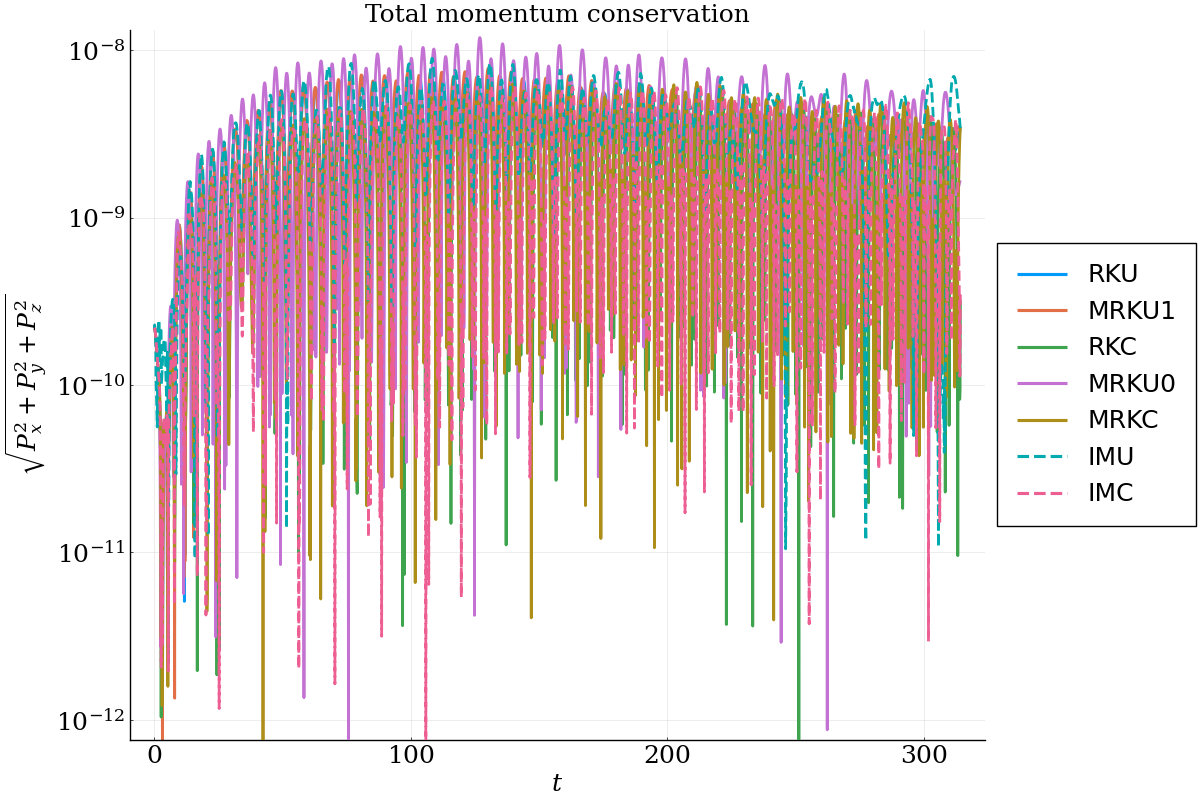}
  \includegraphics[width=0.29\linewidth]{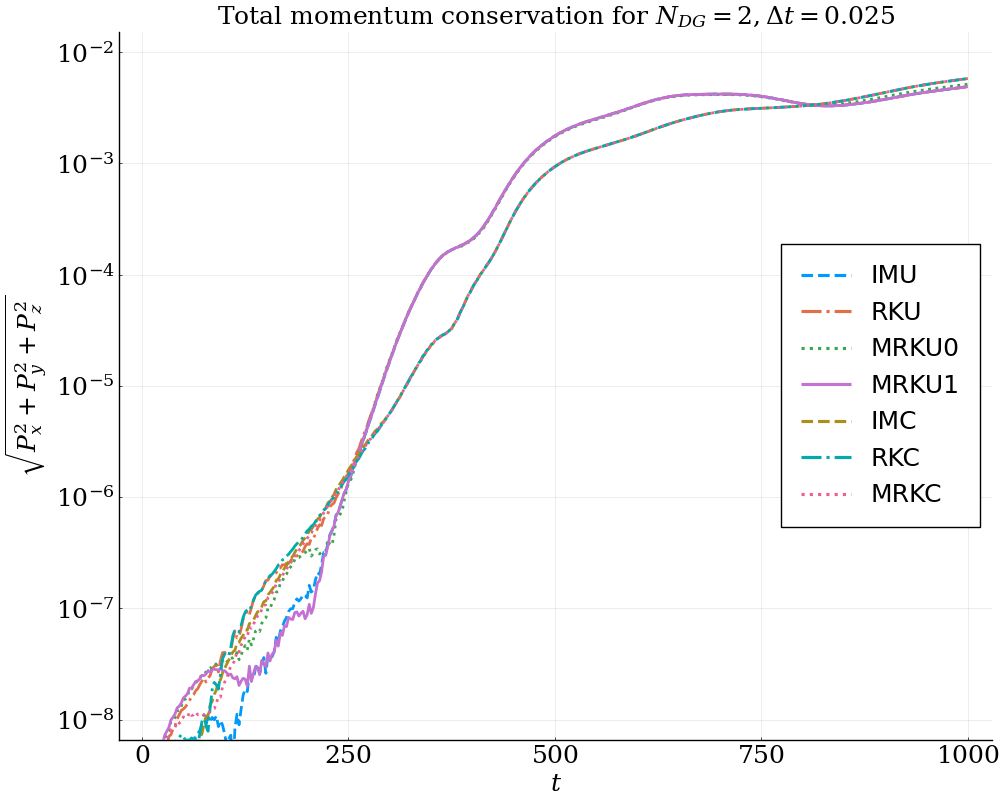}
  \caption{Time evolution of total momentum for different time integration strategies for:
  whistler instability test with $N_{DG} = 1$ and $\Delta t = 0.01$ (left);
  high frequency electromagnetic waves test (middle);
  OT test with $N_{DG} = 2$ and $\Delta t = 0.025$ (right).
}
  \label{fig:momentum}
\end{figure}

The total momentum is a conserved quantity of the Vlasov-Maxwell system at the continuum and for periodic boundary conditions.
However, as shown in \cite[Theorem 5.2]{koshkarov2021},
the Hermite-DG discretization of the problem
does not exactly conserve the discrete total momentum.
If the solution fields are sufficiently regular, the error in the conservation of the discrete  momentum decreases according to the order
of the semi-discrete approximation. In this section we show that, once we introduce a discretization in time, such violation of the momentum conservation is independent of the temporal integrator.
We consider the discrete total momentum defined as
\begin{align*}
    P = \sqrt{P_x^2 +P_y^2 +P_z^2},
\end{align*}
where
\begin{align*}
    P_x = \sum_{I,l} 
    \left (   \int_I\left |  \varphiIl(\xv)   \right |^2 d\xv \right)
    \left ( \left ( \frac{\omega_{ce}}{\omega_{pe}} \right )^2 \left ( E_y^{I,l} B_z^{I,l} - B_y^{I,l}E_z^{I,l} \right )
    + \sum_s m^s \alpha_x^s \alpha_y^s \alpha_z^s \left ( u_x^s C^{s,I,l}_{0,0,0} + \frac{\alpha_x^s}{\sqrt{2}} C^{s,I,l}_{1,0,0}\right ) \right ),
\end{align*}
and $P_y$ and $P_z$ are defined similarly.
Figure~\ref{fig:momentum} shows the evolution of the discrete total momentum (the initial momentum is zero) for different test cases and temporal integrators.
The left panel shows results for the whistler instability test with $N_{DG} = 1$ and $\Delta t = 0.01$, 
where the error in the explicit methods is of the order of machine precision and in the implicit method is of the order of the tolerance of the nonlinear solver.
The middle panel shows results for the high frequency electromagnetic waves test, where the  error is essentially independent
of the time integration strategy.
The right panel shows results for the OT test with $N_{DG} = 2$ and $\Delta t = 0.025$,
where the error is also independent of the time integration strategy.

Another important issue, common to spectral methods, is that the positivity of the distribution function is not guaranteed in the method considered.
As an example, in the most difficult test studied in the paper (the OT test), the reconstructed distribution function has regions with negative values for sufficiently large velocities, whose magnitude at late times and at $v > 3v_T$ (where $v_T$ is the thermal velocity) is $\sim 2$\% of the (positive) peak value. The error associated with the lack of positivity of the distribution function is fairly insensitive to the time discretization strategy.
\par 
We note that the algorithm presented here can be considered as a reduced kinetic model where the distribution function is typically represented with only a few polynomials per velocity direction ($4$ in the OT test). The negativity of the reconstructed distribution function at some values of the velocity has no obvious consequences for the algorithm, as long as it does not introduce significant errors in the first $N_v$ moments of the distribution. Additionally, adding more Hermite moments mitigates the positivity problem effectively.


%% file: conclusion.tex
\section{Conclusions}
\label{sec:conclusions}

We presented the analysis of the conservation properties of the fully-discrete
numerical approximation of the Vlasov-Maxwell equations
provided by the spectral-discontinuous Galerkin method proposed in
\cite{koshkarov2021}. Here the velocity space of the Vlasov equation is expanded in asymmetrically weighted Hermite functions while the spatial part of the Vlasov equation and Maxwell's equations are discretized with the DG method.
The number of particles is conserved in the semi-discrete formulation of the method,
while the total energy is conserved if Maxwell's equations are discretized with central fluxes
and dissipates when upwind fluxes are used. 
Such conservation and dissipation properties are maintained at the fully discrete level
by (implicit) Gauss-Legendre temporal integrators but not by explicit Runge-Kutta schemes. 
Hence, we have
adopted \emph{modified} RK schemes, known as
incremental direction techniques (IDT) \cite{Calvo-HernandezAbreu-Montijano-Randez:2006} and relaxation Runge-Kutta methods (RRK) \cite{Ketcheson19}, and analyzed them in the context of the Hermite-DG discretization of the Vlasov-Maxwell equations. These are
explicit Runge-Kutta methods that maintain the properties of the semi-discrete system by applying a correction step to the regular Runge-Kutta algorithm. The idea behind the method is to tie the correction step to the conservation of a particular quantity (the total energy in this case).

For upwind fluxes in Maxwell's equations, which provide numerical energy dissipation associated with the DG jumps at the cell interfaces, one might be tempted to apply the Runge-Kutta correction step in a way that forces exact conservation of the total energy (i.e. the conservation property of the system at the continuum level and not at the semi-discrete level). Indeed, we have tested this possibility numerically [Eq. (\ref{eq:gamma:def}) with $\varepsilon=0$, scheme MRKU0 in Table 1] against an approach where the correction step is chosen to preserve the upwind energy dissipation [Eq. (\ref{eq:gamma:def}) with $\varepsilon=1$, scheme MRKU1 in Table 1; this is equivalent to enforcing energy conservation at the semi-discrete level such that to preserve upwind energy dissipation]. In all our numerical experiments, MRKU0 consistently showed lower accuracy of the overall solution relative to MRKU1 when used in conjunction with linear polynomials ($N_{DG}=1$) for DG. For instance, in the Orszag-Tang decaying turbulence test, Fig. \ref{fig:OT:jz}, one can see spurious oscillations in the MRKU0 solution (sixth row, first two columns) that are not present in any of the other methods. These oscillations disappear when one increases the order of the spatial discretization ($N_{DG}=2$, sixth row, third and fourth columns in Fig. \ref{fig:OT:jz}).
These results suggest that preserving the properties of the semi-discrete formulation of the method at the fully-discrete level and controlling the numerical dissipation through the accuracy of the spatial discretization might provide a better algorithm when upwind fluxes are used in Maxwell's equations.
\revtwo{For central fluxes, the total energy is a conserved quantity of the evolution problem ensuing from the Hermite-DG discretization since the terms associated with the jumps at cell interfaces disappear. In this case, the modified RK scheme ensures energy conservation at the fully-discrete level.}

In closing, we note that explicit time integrators that can conserve energy at the discrete level, such as those discussed in this paper, are very important for plasma physics applications since these methods are much less computationally expensive than the implicit time integrators that also typically can conserve energy. Depending on the application, implicit methods might still be preferable for their overall better numerical stability. Future work might seek the development of energy-conserving semi-implicit methods that combine the advantages of the two classes of schemes and achieve energy conservation in the discrete, while relaxing the computational cost of fully-implicit schemes (as done recently in the area of Particle-In-Cell methods~\citep{lapenta2017}).

%% file: mainRev.bbl
\begin{thebibliography}{10}

\bibitem{Bittencourt:2004}
J.~A. Bittencourt.
\newblock {\em Fundamentals of Plasma Physics}.
\newblock Springer-Verlag New York, 2004.

\bibitem{Bogacki-Shampine:1989}
P.~Bogacki and L.~F. Shampine.
\newblock {A 3 (2) pair of Runge-Kutta formulas}.
\newblock {\em Applied Mathematics Letters}, 2(4):321--325, 1989.

\bibitem{Boyd-Sandserson:2003}
T.~J.~M. Boyd and J.~J. Sanderson.
\newblock {\em The Physics of Plasmas}.
\newblock Cambridge University Press, 2003.

\bibitem{Calvo-HernandezAbreu-Montijano-Randez:2006}
M.~Calvo, D.~Hern\'andez-Abreu, J.~Montijano, and L.~R\'andez.
\newblock On the preservation of invariants by explicit {R}unge--{K}utta
  methods.
\newblock {\em SIAM Journal on Scientific Computing}, 28(3):868--885, 2006.

\bibitem{CPKS22}
M.~Campos~Pinto, K.~Kormann, and E.~Sonnendr\"{u}cker.
\newblock Variational framework for structure-preserving electromagnetic
  particle-in-cell methods.
\newblock {\em J. Sci. Comput.}, 91(2):Paper No. 46, 39, 2022.

\bibitem{cheng76}
C.~Z. Cheng and G.~Knorr.
\newblock The integration of the {V}lasov equation in configuration space.
\newblock {\em Journal of Computational Physics}, 22(3):330 -- 351, 1976.

\bibitem{Cooper87}
G.~J. Cooper.
\newblock Stability of {R}unge-{K}utta methods for trajectory problems.
\newblock {\em IMA J. Numer. Anal.}, 7(1):1--13, 1987.

\bibitem{DV84}
K.~Dekker and J.~G. Verwer.
\newblock {\em Stability of {R}unge-{K}utta methods for stiff nonlinear
  differential equations}, volume~2 of {\em CWI Monographs}.
\newblock North-Holland Publishing Co., Amsterdam, 1984.

\bibitem{DBM02}
N.~{Del Buono} and C.~Mastroserio.
\newblock Explicit methods based on a class of four stage fourth order
  {R}unge–{K}utta methods for preserving quadratic laws.
\newblock {\em Journal of Computational and Applied Mathematics},
  140(1):231--243, 2002.
\newblock Int. Congress on Computational and Applied Mathematics 2000.

\bibitem{Delzanno:2015}
G.~L. Delzanno.
\newblock Multi-dimensional, fully-implicit, spectral method for the
  {V}lasov-{M}axwell equations with exact conservation laws in discrete form.
\newblock {\em Journal of Computational Physics}, 301:338--356, 2015.

\bibitem{East91}
James~W. Eastwood.
\newblock The virtual particle electromagnetic particle-mesh method.
\newblock {\em Computer Physics Communications}, 64(2):252--266, 1991.

\bibitem{ESF98}
E.~Eich-Soellner and C.~F\"{u}hrer.
\newblock {\em Numerical methods in multibody dynamics}.
\newblock European Consortium for Mathematics in Industry. B. G. Teubner,
  Stuttgart, 1998.

\bibitem{Funaro:2008}
D.~Funaro.
\newblock {\em Polynomial approximation of differential equations}, volume~8.
\newblock Springer Science \& Business Media, 2008.

\bibitem{Gary:2005}
S.~P. Gary.
\newblock {\em Theory of space plasma microinstabilities}.
\newblock Cambridge University Press, 2005.

\bibitem{Glassey:1996}
R.~T. Glassey.
\newblock {\em The {C}auchy problem in kinetic theory}, volume~52.
\newblock SIAM, 1996.

\bibitem{goldston1995}
R.~J. Goldston and P.~H. Rutherford.
\newblock {\em Introduction to Plasma Physics}.
\newblock CRC Press, 1995.

\bibitem{Grad:1949}
H.~Grad.
\newblock On the kinetic theory of rarefied gases.
\newblock {\em Communications on Pure and Applied Mathematics}, 2(4):331--407,
  1949.

\bibitem{GQ05}
V.~Grimm and G.~R.~W. Quispel.
\newblock Geometric integration methods that preserve {L}yapunov functions.
\newblock {\em BIT}, 45(4):709--723, 2005.

\bibitem{HLW06}
E.~Hairer, C.~Lubich, and G.~Wanner.
\newblock {\em Geometric numerical integration. Structure-preserving algorithms
  for ordinary differential equations}, volume~31 of {\em Springer Series in
  Computational Mathematics}.
\newblock Springer, Heidelberg, 2010.

\bibitem{HNW93}
E.~Hairer, S.~P. N{\o}rsett, and G.~Wanner.
\newblock {\em Solving ordinary differential equations {I}. Nonstiff problems},
  volume~8 of {\em Springer Series in Computational Mathematics}.
\newblock Springer-Verlag, Berlin, second edition, 1993.

\bibitem{Hairer-Wanner:1996}
E.~Hairer and G.~Wanner.
\newblock {\em Solving ordinary differential equations {II}. Stiff and
  differential-algebraic problems}, volume~14 of {\em Springer Series in
  Computational Mathematics}.
\newblock Springer-Verlag, Berlin, 2010.

\bibitem{He15}
Y.~He, Y.~Sun, J.~Liu, and H.~Qin.
\newblock Volume-preserving algorithms for charged particle dynamics.
\newblock {\em Journal of Computational Physics}, 281:135--147, 2015.

\bibitem{hesthaven2004high}
J.~S. Hesthaven and T.~Warburton.
\newblock High--order nodal discontinuous {G}alerkin methods for the {M}axwell
  eigenvalue problem.
\newblock {\em Philosophical Transactions of the Royal Society of London.
  Series A: Mathematical, Physical and Engineering Sciences},
  362(1816):493--524, 2004.

\bibitem{IserlesBook}
A.~Iserles.
\newblock {\em A first course in the numerical analysis of differential
  equations}.
\newblock Cambridge Texts in Applied Mathematics. Cambridge University Press,
  Cambridge, 1996.

\bibitem{joyce71}
G.~Joyce, G.~Knorr, and H.~K. Meier.
\newblock Numerical integration methods of the {V}lasov equation.
\newblock {\em Journal of Computational Physics}, 8(1):53--63, 1971.

\bibitem{Juno-Hakim-TenBarge-Shi-Dorland:2018}
J.~Juno, A.~Hakim, J.~TenBarge, E.~Shi, and W.~Dorland.
\newblock Discontinuous {G}alerkin algorithms for fully kinetic plasmas.
\newblock {\em Journal of Computational Physics}, 353:110--147, 2018.

\bibitem{Ketcheson19}
D.~I. Ketcheson.
\newblock Relaxation {R}unge-{K}utta methods: conservation and stability for
  inner-product norms.
\newblock {\em SIAM J. Numer. Anal.}, 57(6):2850--2870, 2019.

\bibitem{KS21}
K.~Kormann and E.~Sonnendr\"ucker.
\newblock Energy-conserving time propagation for a structure-preserving
  particle-in-cell vlasov–maxwell solver.
\newblock {\em Journal of Computational Physics}, 425:109890, 2021.

\bibitem{koshkarov2021}
O.~Koshkarov, G.~Manzini, G.~L. Delzanno, C.~Pagliantini, and V.~Roytershteyn.
\newblock The multi-dimensional {H}ermite-discontinuous {G}alerkin method for
  the {V}lasov--{M}axwell equations.
\newblock {\em Computer Physics Communications}, 264:107866, 2021.

\bibitem{GEMPIC}
M.~Kraus, K.~Kormann, P.~J. Morrison, and E.~Sonnendr\"ucker.
\newblock Gempic: geometric electromagnetic particle-in-cell methods.
\newblock {\em Journal of Plasma Physics}, 83(4):905830401, 2017.

\bibitem{lapenta2017}
G.~Lapenta.
\newblock Exactly energy conserving semi-implicit particle in cell formulation.
\newblock {\em Journal of Computational Physics}, 334:349--366, 2017.

\bibitem{Le70}
H.~R. Lewis.
\newblock Energy-conserving numerical approximations for vlasov plasmas.
\newblock {\em Journal of Computational Physics}, 6(1):136--141, 1970.

\bibitem{Le72}
H.~R. Lewis.
\newblock Variational algorithms for numerical simulation of collisionless
  plasma with point particles including electromagnetic interactions.
\newblock {\em Journal of Computational Physics}, 10(3):400--419, 1972.

\bibitem{Manzini-Delzanno-Markidis-Vencels:2016}
G.~Manzini, G.~L. Delzanno, J.~Vencels, and S.~Markidis.
\newblock A {L}egendre-{F}ourier spectral method with exact conservation laws
  for the {V}lasov-{P}oisson system.
\newblock {\em Journal of Computational Physics}, 317:82--107, 2016.

\bibitem{Orszag:1979}
S.~A. Orszag and C.-M. Tang.
\newblock Small-scale structure of two-dimensional magnetohydrodynamic
  turbulence.
\newblock {\em Journal of Fluid Mechanics}, 90(1):129--143, 1979.

\bibitem{Parashar2009}
T.~N. Parashar, M.~A. Shay, P.~A. Cassak, and W.~H. Matthaeus.
\newblock Kinetic dissipation and anisotropic heating in a turbulent
  collisionless plasma.
\newblock {\em Physics of Plasmas}, 16(3):032310, 2009.

\bibitem{sarmany2007dispersion}
D.~S{\'a}rm{\'a}ny, M.~A. Botchev, and J.~J.~W. van~der Vegt.
\newblock Dispersion and dissipation error in high-order {R}unge-{K}utta
  discontinuous {G}alerkin discretisations of the {M}axwell equations.
\newblock {\em Journal of Scientific Computing}, 33(1):47--74, 2007.

\bibitem{Schumer-Holloway:1998}
J.~W. Schumer and J.~P. Holloway.
\newblock Vlasov simulations using velocity-scaled {H}ermite representations.
\newblock {\em J. Comput. Phys.}, 144(2):626--661, 1998.

\bibitem{Tao16}
M.~Tao.
\newblock Explicit high-order symplectic integrators for charged particles in
  general electromagnetic fields.
\newblock {\em Journal of Computational Physics}, 327:245--251, 2016.

\bibitem{Wan2012}
M.~Wan, W.H. Matthaeus, H.~Karimabadi, V.~Roytershteyn, M.~Shay, P.~Wu,
  W.~Daughton, B.~Loring, and S.C. Chapman.
\newblock {Intermittent Dissipation at Kinetic Scales in Collisionless Plasma
  Turbulence}.
\newblock {\em Physical Review Letters}, 109(19):195001, nov 2012.

\end{thebibliography}
